\documentclass[a4paper,12pt]{amsart} 
\usepackage{amsmath,amssymb}
\usepackage{todonotes}
\usepackage{geometry}
\usepackage{comment}
\usepackage{algorithm}
\usepackage{algpseudocode}
\usepackage{subcaption}
\usepackage{hyperref}
\usepackage{booktabs}
\usepackage{adjustbox}

\numberwithin{figure}{section}
\numberwithin{table}{section}
\numberwithin{equation}{section}
\newtheorem{theorem}{Theorem}[section]

\newtheorem{lemma}[theorem]{Lemma}

\newtheorem{assumption}[theorem]{Assumption}
\newtheorem{remark}[theorem]{Remark}
\newtheorem{definition}[theorem]{Definition}

\newcommand{\rset}{\mathbb{R}}
\newcommand{\Bw}{\bar f}

\newcounter{testcase}
\newcommand{\Testcase}[1]{
  \refstepcounter{testcase}
  \par\medskip
  \noindent\textbf{Test~\thetestcase.} #1\par\medskip
}

\begin{document}
\title[]{Convergence for adaptive resampling of random Fourier features} 
\begin{abstract}

The machine learning random Fourier feature method for data in high dimension is 
computationally and theoretically attractive since the optimization is based on a convex standard least squares problem and independent sampling of Fourier frequencies.
The challenge is to sample the Fourier frequencies well.
This work proves convergence of a data adaptive method based on resampling
the frequencies asymptotically optimally, as the number
of nodes and amount of data tend to infinity. Numerical results based on resampling and adaptive random walk steps together with
approximations of the least squares problem by conjugate gradient iterations 
confirm the analysis for regression and classification problems.

\end{abstract}
\author[Huang]{Xin Huang}
\address{Department of Mathematics and Mathematical Statistics, Ume{\aa} University, 901 87 Ume{\aa}, SWEDEN}
\email{xin.huang@umu.se}

\author[Kammonen]{Aku Kammonen}
\address{Computer, Electrical and Mathematical Sciences and Engineering,
4700 King Abdullah University of Science and Technology (KAUST),
Thuwal 23955-6900, Kingdom of Saudi Arabia}
\email{akujaakko.kammonen@kaust.edu.sa}

\author[Pandey]{Anamika Pandey}
\address{Chair of Mathematics for Uncertainty Quantification, RWTH Aachen University, 52062 Aachen, Germany} 
\email{pandey@uq.rwth-aachen.de}

\author[Sandberg]{Mattias Sandberg}
\address{Institutionen f\"or Matematik, Kungl. Tekniska H\"ogskolan, 100 44 Stockholm, Sweden}
\email{msandb@kth.se}

\author[von Schwerin]{Erik~von~Schwerin}
\address{Computer, Electrical and Mathematical Sciences and Engineering,
4700 King Abdullah University of Science and Technology (KAUST),
Thuwal 23955-6900, Kingdom of Saudi Arabia}
\email{erik.vonschwerin@kaust.edu.sa}

\author[Szepessy]{Anders Szepessy}
\address{Institutionen f\"or Matematik, Kungl. Tekniska H\"ogskolan, 100 44 Stockholm, Sweden}
\email{szepessy@kth.se}

\author[Tempone]{Ra\'{u}l~Tempone}
\address{Computer, Electrical and Mathematical Sciences and Engineering,
4700 King Abdullah University of Science and Technology (KAUST),
Thuwal 23955-6900, Kingdom of Saudi Arabia; Chair of Mathematics for Uncertainty Quantification, RWTH Aachen University, 52062 Aachen, Germany; Alexander von Humboldt Professor in Mathematics for Uncertainty Quantification, RWTH Aachen University, 52062 Aachen, Germany.}
\email{raul.tempone@kaust.edu.sa}


\maketitle
\tableofcontents

\section{Introduction to resampling random features}
\label{sec:intro}
A basic problem in AI and machine learning is to find a neural network function that approximately reproduces the mapping in provided data \cite{bach_book,understand_ml}. The aim in this work is to prove that an adaptive version of the random Fourier feature method, based on resampling, generates a neural network function with an asymptotically optimal expected least squares error, provided the initial data sampling is sufficiently close. 

The learning of neural network functions has two steps:  choosing a set of neural network functions and defining an algorithm to find a neural network function that in some sense approximates the given data well. More precisely
given data $\{(x_j,y_j)\in \mathbb R^d\times \mathbb R\, : \, j=1,\ldots, J\}$, where it is expected that $y_j=f(x_j) + \mbox{noise}$,
the basic learning problem is to find a neural network function that approximates $f$.
A standard machine learning formulation to determine a neural network approximation for $f$ is to minimize the empirical loss
\begin{equation}\label{goal}
\min_{\zeta\in \mathcal N_K} \big({J^{-1}\sum_{j=1}^J {|y_j-\zeta(x_j)|^2}}
 +{\lambda}\sum_{k=1}^K|\hat\zeta_k|^2\big)\,,
\end{equation}
where
\begin{equation*}\label{NN}
\mathcal N_K:=\Big\{\zeta( x)=\sum_{k=1}^K\hat\zeta_k\sigma(\omega_k\cdot x+\omega_0) 
 \Big\}
\end{equation*}
is a set of neural network functions with one hidden layer, based on the parameters $(\hat\zeta_k,\omega_k)\in\mathbb R\times\mathbb R^d$, for $k=1,\ldots,K$, and $\omega_0\in\mathbb R$,  including $\lambda\ge 0$ as a Tikhonov regularization parameter.
The activation function $\sigma:\mathbb R\to \mathbb R$ is assumed to be locally bounded and not a polynomial \cite{activation}. Common choices are the sigmoid function $\sigma(z)=1/(1+e^{-z})$ and the rectified linear unit $\sigma(z)=\max(0,z)$.
The non-convex minimization \eqref{goal} is usually approximated by iterations based on the stochastic gradient method or related variants \cite{bach_book,understand_ml}. 
A main challenge in machine learning, both in practice  and theoretically, is to obtain reliable convergence of the iterations \cite{bach_book,understand_ml, weinan_ICM, belkin}. 

The aim here is to study convergence for a different but related machine learning problem, namely a formulation based on random Fourier features $\zeta\in \mathcal N_K$. %  adapting to the data. 
The random Fourier feature method, 
 using the Fourier activation function 
$\sigma(\omega\cdot x)=e^{{\rm i}\omega\cdot x}$, was 
derived from the kernel method in \cite{rahimi_recht} and takes the following form 
\begin{equation}\label{rffm}
\min_{\hat\zeta\in \mathbb C^K} \big({J^{-1}\sum_{j=1}^J {|y_j-\zeta(x_j)|^2}}
 +{\lambda}\underbrace{\sum_{k=1}^K|\hat\zeta_k|^2}_{=:\|\hat\zeta\|_2^2}\big)
\end{equation}
where the frequency parameters $\omega_1,\ldots, \omega_K$ are random independent samples from a probability distribution $p:\mathbb R^d\to [0,\infty)$ and the amplitudes $\hat\zeta_k\in\mathbb C$. The machine learning random Fourier feature method  is computationally and theoretically attractive since the optimization is based on a convex standard least squares problem and independent sampling of Fourier frequencies, applicable to data in high dimension. The challenge is to sample the Fourier frequencies well to obtain a reliable method providing small generalization error. 
In the case of infinite amount of data, $J=\infty$, and $y_j=f(x_j)$   the generalization error (i.e., the testing error) and training error become equal and satisfy %for the Fourier activation function
\begin{equation}\label{Rff}
\min_{\zeta\in\mathcal N_K}\mathbb E_x[|f(x)-\zeta(x)|^2] +\lambda\|\hat\zeta\|_2^2
\le \mathbb E_\omega\big[\min_{\hat\zeta\in\mathbb C^K}\big(\mathbb E_x[|f(x)-\zeta(x)|^2]+\lambda\|\hat\zeta\|_2^2\big)\big]\le C_p K^{-1}\,,
\end{equation}
for a certain constant $C_p$ given in \eqref{eta2}. It is proved by first using the fact that the minimum over $\omega$ is bounded by its average and then applying Monte Carlo approximation of the Fourier representation of $f$ together with the independence of $\omega_k$ sampled from $p$, see \cite{barron,barron2,bach_book,weinan_ICM, ajpma} and Section \ref{Monte_Carlo}. The generalization error for the random Fourier feature method therefore also bounds the generalization error for the neural network minimization problem \eqref{goal}. 

The aim here is to prove convergence of a variant of the random Fourier feature method with independent samples that adapts the sampling distribution $p$ to the data and thereby minimizes  the constant $C_p$  in the generalization error bound \eqref{Rff}.  An adaptive random Fourier feature method, using the Metropolis method, was introduced in \cite{ajpma}. The subsequent work \cite{Aku2} formulates an adaptive Fourier feature method that includes resampling, 
substantially improving  
the computational results. 

This work formulates a resampling method, following \cite{Aku2}, %that achieves asymptotically optimal convergence rates 
 for the random Fourier feature method  that adapts the feature sampling distribution $p$ to the data, aiming to minimize the constant $C_p$ in \eqref{Rff}, and analyzes its convergence properties. The analysis
 uses four main ideas to prove convergence for the adaptive resampling method obtaining  a nearly optimal feature sampling distribution, namely to 
\begin{enumerate}
\item%[(0)] 
derive the minimal rate constant $C_{p}$,
\item %[(1)] 
use {\it Fourier series} to obtain a discrete set of basis functions,
\item%[(2)] 
establish {\it equal amplitudes} for the same basis function using Tikhonov regularization,
\item%[(3)] 
{\it cut off small amplitudes} to have bounded sums of amplitudes, using the generalization error estimate and a decay of the Fourier coefficients for periodic functions  $f:[-L,L]^d\to \mathbb C$ with bounded derivatives of order $\ell > d$.
\end{enumerate}
The following three sections formulate the random Fourier feature method, derive the equal amplitude property, and describe the resampling procedure.
\subsection{Formulation of the random Fourier feature method}
Assume $f:\mathbb R^d/(2L\mathbb Z)^d\to \mathbb C$ is component-wise $2L$-periodic 
with {\it the Fourier series}
\[
f(x)=\sum_{n\in \mathbb Z^d} \hat f(\omega_n)e^{{\rm i}\omega_n\cdot x}\,,
\]
where $\omega_n:=\frac{\pi}{L}n$ for $n\in \mathbb Z^d$
and $\mathbb T^d:=[-L,L]^d$ for some $L>0$\,.
Suppose we have data points $\{(x_j, y_j)\}_{j=1}^J$
where $\{x_j\in\mathbb R^d\,,\ j=1,\ldots,J \}$ are independent
samples from a density
$\rho:\mathbb T^d\to [0,\infty)$
with noisy values\[y_j=f(x_j)+\xi_j\,,\]
and the perturbations $\{\xi_j\in\mathbb R\}_{j=1}^J$  are independent identically distributed random variables with mean zero and bounded variance.

Consider the random Fourier feature problem
\begin{equation}\label{min}
    \min_{\hat\beta\in\mathbb C^K}\big(\frac{1}{J}\sum_{j=1}^J|\beta(x_j)- y_j|^2 +\lambda\sum_{k=1}^K|\hat\beta_k|^2\big)
   %+\lambda_2\sum_{k=1}^K|\hat\beta_k|\big)
\end{equation}
for the neural network
\[
\beta(x):= \sum_{k=1}^K \hat\beta_k e^{{\rm i}\nu_k\cdot x}\,,
\]
with  a positive Tikhonov regularization parameter $\lambda$,
where 
$\{\nu_k\in \frac{\pi}{L}\mathbb Z^d \,,\ k=1,\ldots,K\}$
are independent samples from a distribution $p:\frac{\pi}{L}\mathbb Z^d\to [0,\infty)$.

The core of the problem of sampling frequencies adaptively well 
is that the rate constant is explicitly known, namely
$C_p=\sum_{n\in\mathbb Z^d}\frac{|\hat f(\omega_n)|^2}{p(\omega_n)}$, as shown in \eqref{b_gen}, 
which is minimized by using the optimal density 
$p_*(\omega_n)=\frac{|\hat f(\omega_n)|}{\sum_{n\in\mathbb Z^d}|\hat f(\omega_n)|}$, see \eqref{p_opt}.
The main idea in this work is  that the known amplitudes $|\hat\beta_k|$, 
together with the equal amplitude property and frequency cutoff explained in the following sections, can be used to approximate the unknown Fourier amplitudes $|\hat f(\omega_n)|$ in the optimal density $p_*$.

\subsection{The equal amplitude property }
Let 
\[
\bar\beta_n:=\sum_{\{k:\,\nu_k=\omega_n\}}\hat\beta_k
\]
be the sum of the neural network amplitudes corresponding to the same basis function $e^{{\rm i}\omega_n\cdot x}$.
The optimization problem \eqref{min} implies
that the amplitudes $\hat\beta_k$,  corresponding to the same
basis function $e^{{\rm i}\omega_n\cdot x}$, satisfy
\begin{equation}\label{equal}\hat\beta_k=\frac{\bar\beta_n}{m_n}
\end{equation}
 with $m_n:=\#\{k:\nu_k=\omega_n\}$\,,
since the regularization term 
\[
\lambda\sum_{\{k:\nu_k=\omega_n\}}|\hat\beta_k|^2
%+\lambda_2\sum_{\{k:\nu_k=\omega_n\}}|\hat\beta_k|
\,,
\]
with the linear constraint $\sum_{\{k:\nu_k=\omega_n\}}\hat\beta_k=\bar\beta_n$\,,
 is minimized when these { amplitudes $\hat\beta_k$ are equal}.
%by the least squares problem \eqref{min} for $\lambda>0$.
This property, that  the amplitudes $\hat\beta_k$ corresponding to the same
basis function are equal,
yields 
\begin{equation}\label{abs}
\sum_{\{k:\nu_k=\omega_n\}}|\hat\beta_k|=|\bar\beta_n|\,.
\end{equation}

\subsection{The cutoff }
Using a resampling method without a cutoff, based on that the frequency $\nu_k$ has the probability
\[ \frac{|\hat\beta_k|}{\sum_{k=1}^K|\hat\beta_k|}\,,
\]
implies therefore  by \eqref{abs} %,  including multiples of $\nu_k$, 
that we in fact sample the lattice point $\omega_n=\nu_k$ with the probability
\[
\frac{|\bar\beta_n|}{\sum_{\{n:|\bar\beta_n|\ge 0\}} |\bar\beta_n|}\,.
\]
We will use the corresponding distribution
\begin{equation}\label{p_resam}
\bar p_n:=\frac{|\bar\beta_n|}{\sum_{\{n:|\bar\beta_n|\ge\epsilon\}} |\bar\beta_n|}
% p_n = \frac{|\bar\beta_n|}{\sum_{n\in\mathbb Z^d}|\bar\beta_n|}\,.
\end{equation}
where amplitudes smaller than $\epsilon>0$ are cut off, as motivated in Section \ref{resamp_sec}.

\subsection{The resampling procedure and cutoff }
In the case of an infinite amount of data, $J=\infty$, the random Fourier feature method has the generalization error
\begin{equation}\label{eta2}
\mathbb E_p\big[\mathbb E_x[|\beta(x)-f(x)|^2] +\lambda\sum_{k=1}^K|\hat\beta_k|^2\big]\le K^{-1}(1+\lambda)\sum_{n\in\mathbb Z^d} \frac{|\hat f(\omega_n)|^2}{p(\omega_n)}=:K^{-1}(1+\lambda) C_p(f)\,,  
\end{equation}
which is well known and proved by
Monte Carlo approximation of the Fourier representation of the approximated function in Section \ref{Monte_Carlo}.

Minimization with respect to the probability distribution $p$, for the frequency sampling, implies that the 
minimal rate constant $C_{p_*}(f) = \min_{p}C_p(f)$ 
is obtained for the optimal distribution
\begin{equation}\label{p*1}
p_*(\omega_n):=\frac{|\hat f(\omega_n)|}{\sum_{n\in\mathbb Z^d}|\hat f(\omega_n)|}\,,
\end{equation}
as shown in Section \ref{sec_opt_p}.
Section \ref{resamp_sec}  formulates a resampling method that picks independent frequency samples, with replacement, from the distribution
\begin{equation}\label{pn_def}
\bar p(\omega_n,\bar\beta)= \frac{|\bar\beta_n|}{\sum_{\{n:|\bar\beta_n|\ge\epsilon\}} |\bar\beta_n|}\,,
\end{equation}
using a certain cutoff $\epsilon$ for the amplitudes.
 In addition, cut  
 frequencies are sampled with a certain small proportion $q_\epsilon$.
 The cutoff $\epsilon$ and proportion $q_\epsilon$ depend on the number of nodes $K$ and the amount of data $J$. % in the case $J<\infty$.
 % and the decay of the Fourier coefficients $|\hat f(\omega_n)|$\,.

\subsection{Computational demonstrations}\label{Sec_1_5}
The main inspiration of the work here is the adaptive random feature Algorithm \ref{alg_FT} constructed in \cite{Aku2}.
\begin{algorithm}[!ht]
\caption{Adaptive resampling with random walk}
\label{alg_FT}
\footnotesize{
\begin{algorithmic}
\State {\bfseries Input:} $\{(x_j, y_j)\}_{j=1}^J$ $~\{\textrm{data}\}$
\State {\bfseries Output:} $x\mapsto\sum_{k=1}^K\hat\beta_k e^{{\rm i}\omega_k\cdot x}$
\State Choose a number of resampling iterations $N$, random walk step size $\delta$ and  Tikhonov parameter $\lambda$
\State  Set $(\omega_1,\ldots,\omega_K) =0$
%\gets \mbox{standard normal in $\rset^{Kd}$}$ 
%\Comment{initialization}
\For{$n = 1$ {\bfseries to} $N$}
\State $\zeta \gets \mbox{standard normal   in $\rset^{Kd}$}$
\State $\omega \gets 
\mbox{$\omega + \delta \zeta$} $ $\{\textrm{random walk}\}$
\State $\hat{\beta} \gets$ solution of the least squares problem \eqref{min} using $\omega$
\State $\omega \gets$ $K$ independent resamples from $\{\omega_1,\ldots, \omega_K\}$   with probability $|\hat\beta_k|/\sum_{\ell}|\hat\beta_\ell|$ for $\omega_k$ $\{\textrm{resampling}\}$
\EndFor
\State $\hat{\beta} \gets \mbox{solution of the least squares problem \eqref{min} using $\omega$}$
\State $x\mapsto\sum_{k=1}^K\hat\beta_k e^{{\rm i}\omega_k\cdot x}$
\end{algorithmic}}
\end{algorithm}
Algorithm~\ref{alg_FT} is motivated as follows.
The initial frequencies can be far from optimally distributed, leading to a large rate constant $C_p(f)$, and too large generalization error. The aim of the resampling step is
to weight the present frequencies approximately as the optimal distribution $p_*$. However, many frequencies
with large weight
in the optimal distribution $p_*$ are not present initially. 
The purpose of iteratively combining each resampling with a random walk step is to provide new input to the next resampling step,
eventually yielding a nearly optimally distributed set of frequencies. 

\begin{figure}[!htbp]
\centering
\begin{minipage}[t]{0.32\textwidth}
    \centering
    \includegraphics[width=\linewidth]{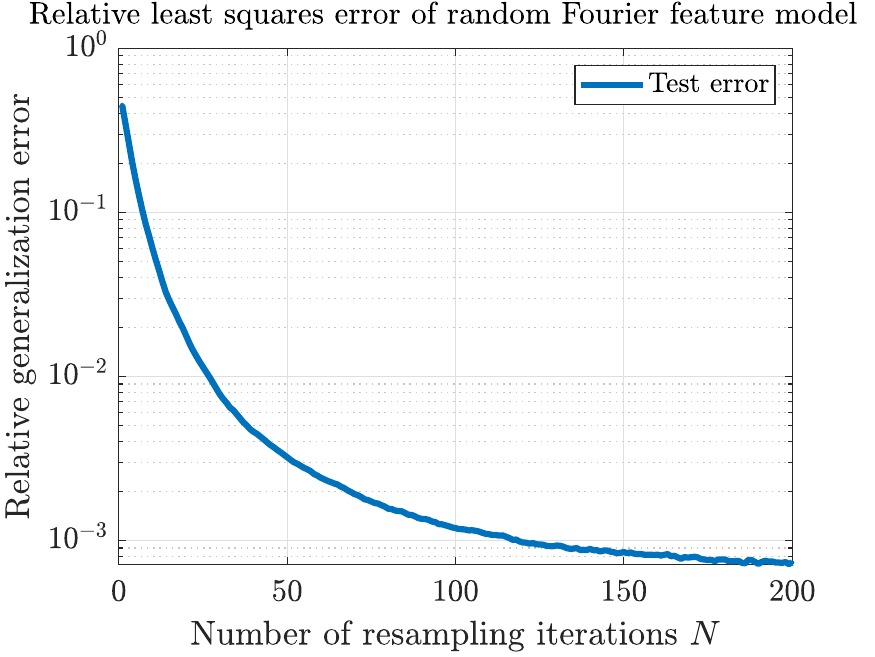}
\end{minipage}%
\hfill
\begin{minipage}[t]{0.32\textwidth}
    \centering
    \includegraphics[width=\linewidth]{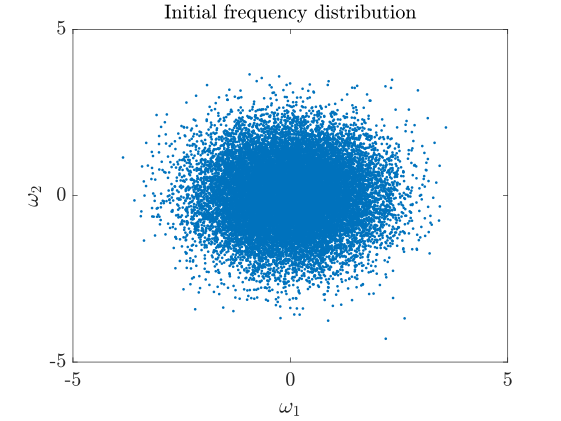}
\end{minipage}%
\hfill
\begin{minipage}[t]{0.32\textwidth}
    \centering
    \includegraphics[width=\linewidth]{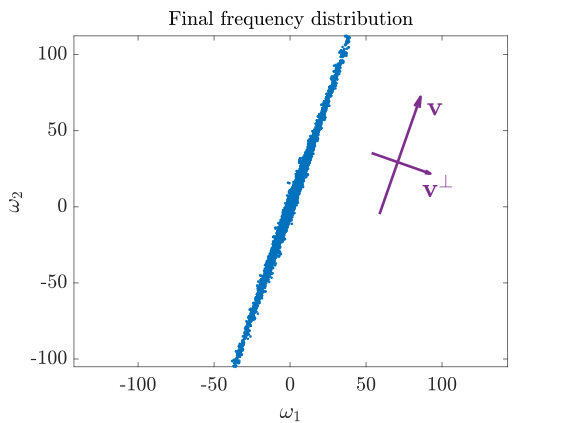}
\end{minipage}
\caption{
Results for the target function $f(x) = e^{-|v\cdot x|/a}e^{-|x|^2/2}$ with parameter $a=0.1$ and a randomly selected direction vector $v=(0.3308,0.9437)$ in dimension $d=2$, using Algorithm~\ref{alg_FT} under non-periodic setting; the training data set size $J=1.5\times 10^4$, frequency sample size $K=1.5J$, random walk step size $\delta=0.5$, and Tikhonov regularization parameter $\lambda=\frac{1}{100}KJ^{-\frac{1}{2}}$. \textbf{Left}: Relative generalization error as a function of resampling iterations. \textbf{Center}: Initial frequency samples drawn from a standard normal distribution. \textbf{Right}: Final frequency samples after 200 resampling steps. 
}
\label{fig:f29}
\end{figure}

\begin{figure}[!htbp]
\centering
\begin{minipage}[t]{0.32\textwidth}
    \centering 
    \includegraphics[width=\linewidth]{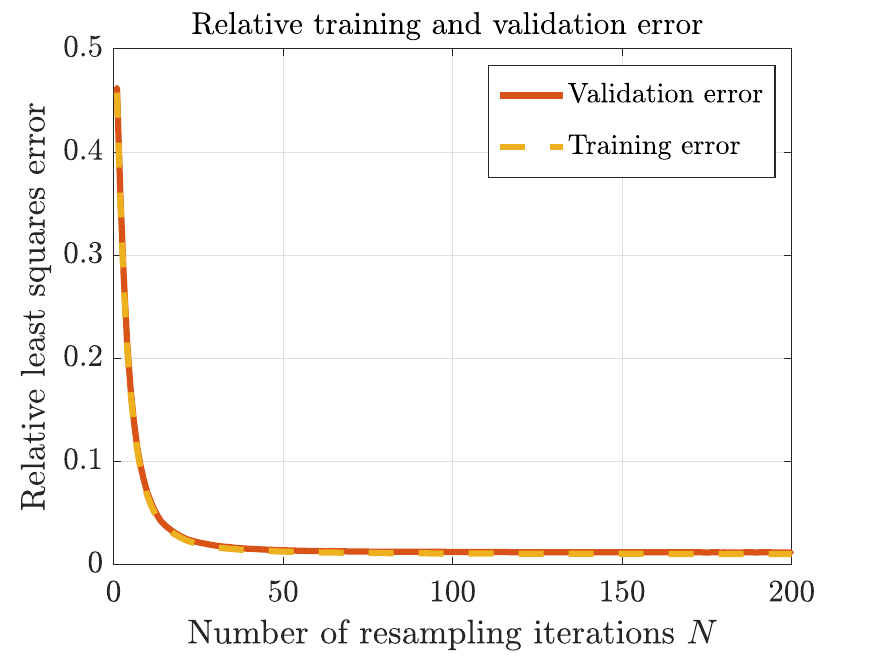}
\end{minipage}
\hfill
\begin{minipage}[t]{0.32\textwidth}
    \centering 
    \includegraphics[width=\linewidth]{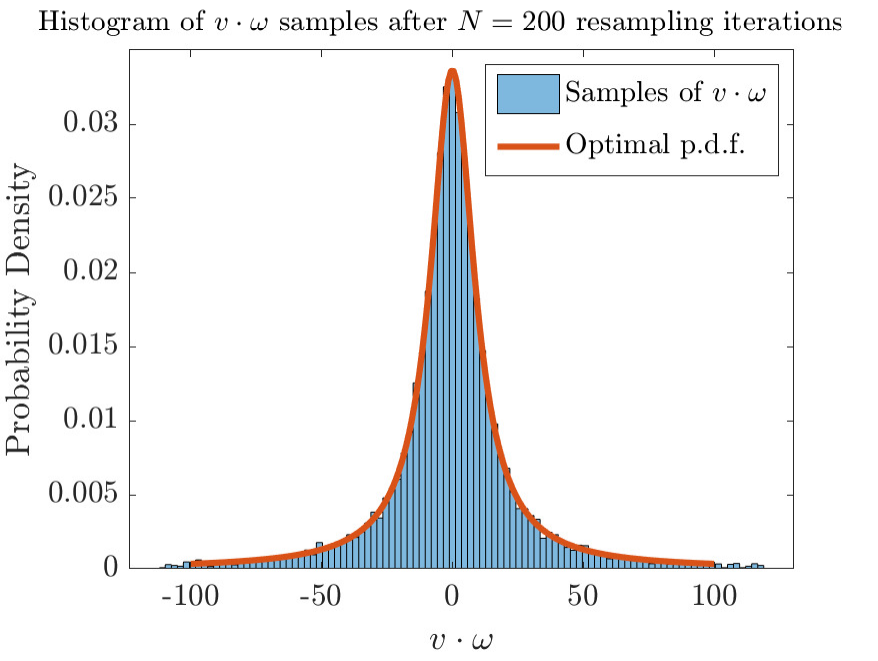}
\end{minipage}
\hfill
\begin{minipage}[t]{0.32\textwidth}
    \centering 
    \includegraphics[width=\linewidth]{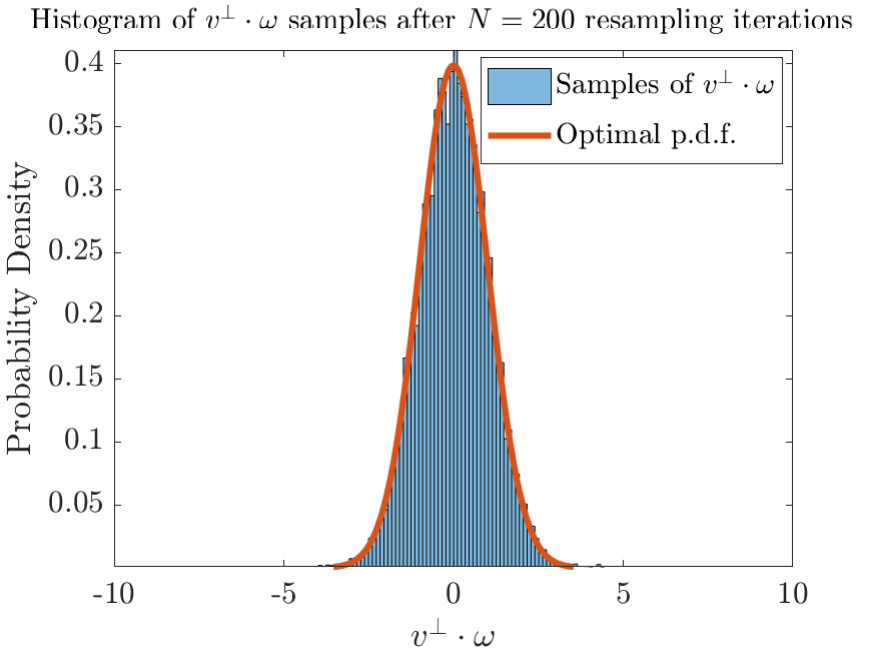}
\end{minipage}
\caption{Frequency distributions after 200 random walk/resampling iterations using Algorithm~\ref{alg_FT} under non-periodic setting.
\textbf{Left}: Relative training and validation error as a function of resampling iterations.
\textbf{Center}: Histogram of projected frequency samples $v\cdot \omega$. 
\textbf{Right}: Histogram of projected frequency samples $v^\perp\cdot \omega$. The plots compare the learned (blue) and optimal (red) frequency probability density function (p.d.f.). 
The hyperparameter settings used here for Algorithm~\ref{alg_FT} are the same as those in Figure~\ref{fig:f29}.
%: projection along \( v \cdot \omega \) (left) and orthogonal projection \( v^\perp \cdot \omega \) (right), showing comparison between learned (blue), and optimal (red) distributions of frequencies. For the implementation of Algorithm~\ref{alg_FT}, the hyperparameter setting is the same as in Figure~\ref{fig:f29}.
}
\label{fig:f30}
\end{figure}

\begin{figure}[!htbp]
\centering
\begin{minipage}[t]{0.32\textwidth}
    \centering
    \includegraphics[width=\linewidth]
    {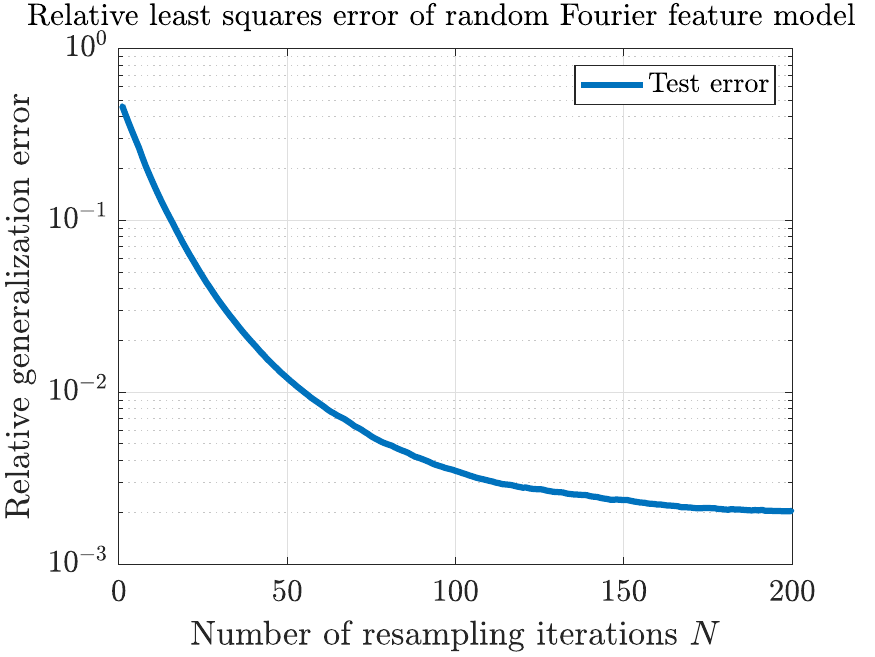}
\end{minipage}%
\hfill
\begin{minipage}[t]{0.32\textwidth}
    \centering
    \includegraphics[width=\linewidth]{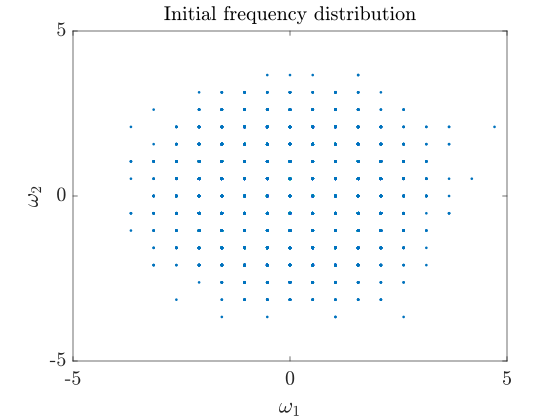}
\end{minipage}%
\hfill
\begin{minipage}[t]{0.32\textwidth}
    \centering
    \includegraphics[width=\linewidth]{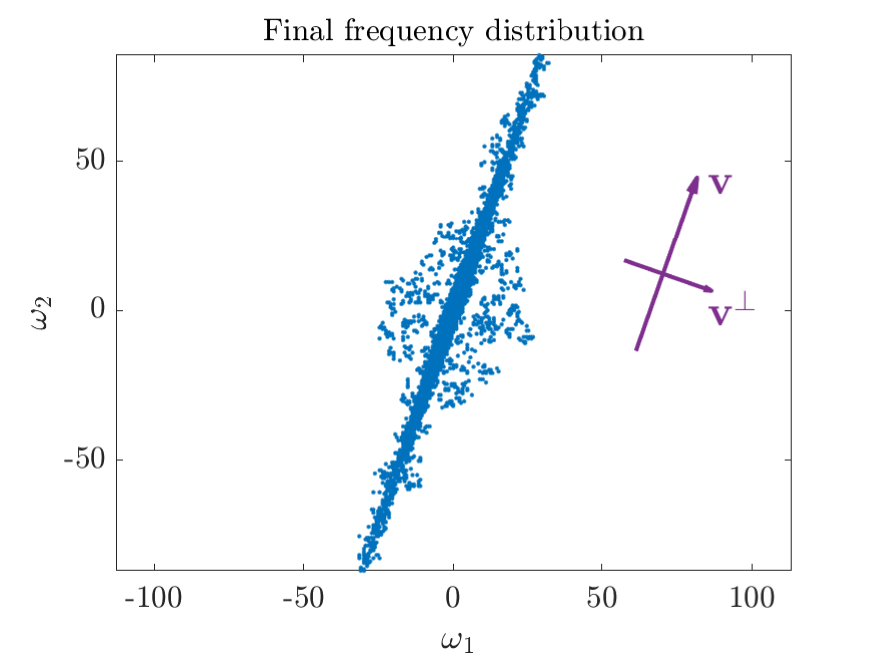}
\end{minipage}
\caption{Results for target function $f(x) = e^{-|v\cdot x|/a}e^{-|x|^2/2}$ with parameter $a=0.1$ in dimension $d=2$, using Algorithm~\ref{alg1} under a periodic setting with period $q = 2L= 12$. \textbf{Left}: Relative generalization error as a function of resampling iterations. \textbf{Center}: Initial frequency samples drawn from a standard normal distribution and projected onto lattice grid points. \textbf{Right}: Final frequency samples after 200 resampling steps. The hyperparameter settings used here for Algorithm~\ref{alg1} are the same as those in Figure~\ref{fig:f29}, except $\delta=0.2$, $\lambda=\frac{1}{500}KJ^{-\frac{1}{2}}$, and cutoff parameter $\epsilon=\frac{1}{200}K^{-\frac{1}{2}}$.
}
\label{fig:f27}
\end{figure}

\begin{figure}[!htbp]
\centering
\begin{minipage}[t]{0.32\textwidth}
    \centering 
    \includegraphics[width=\linewidth]{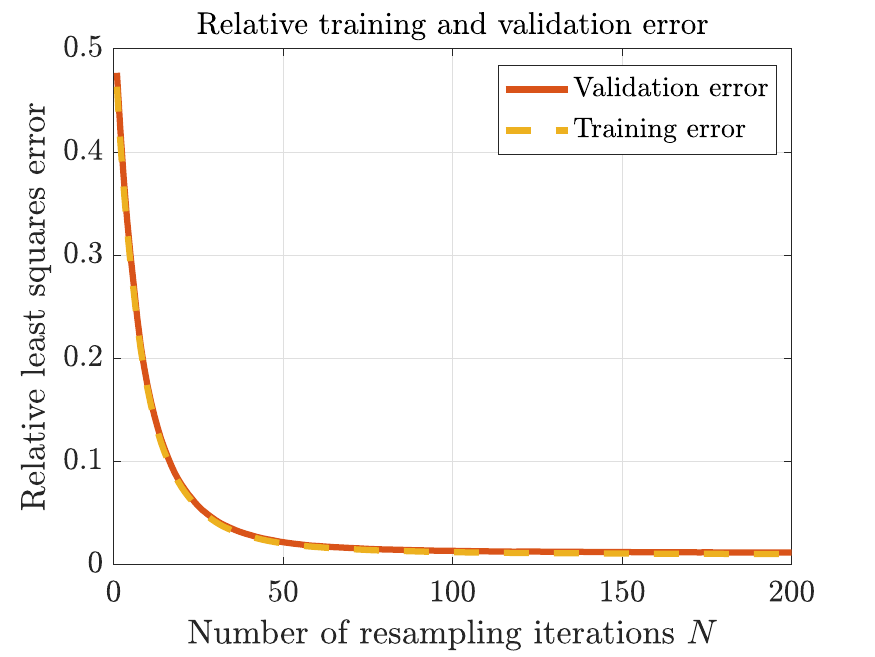}
\end{minipage}
\hfill
\begin{minipage}[t]{0.32\textwidth}
    \centering 
    \includegraphics[width=\linewidth]{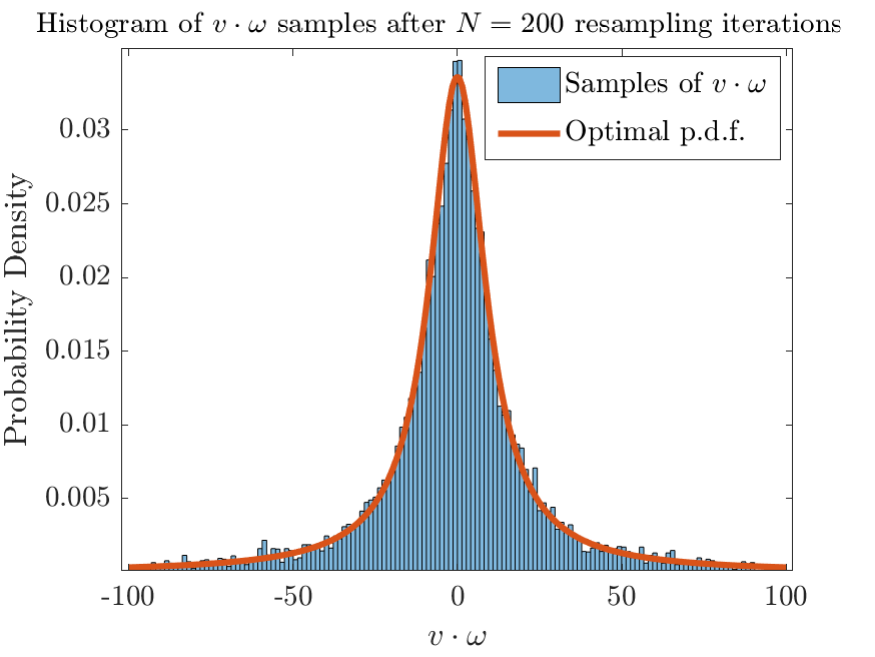}
\end{minipage}
\hfill
\begin{minipage}[t]{0.32\textwidth}
    \centering 
    \includegraphics[width=\linewidth]{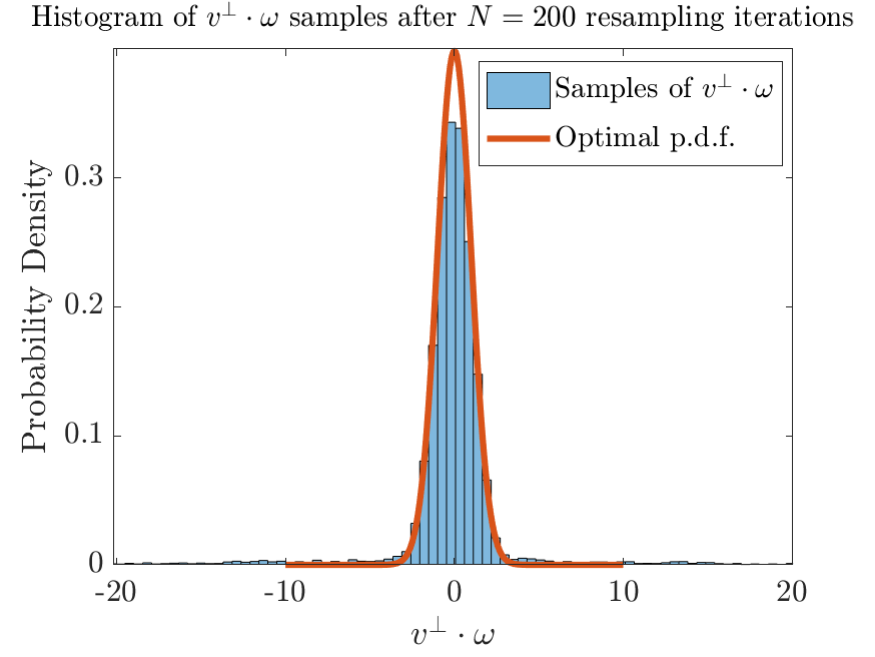}
\end{minipage}
\caption{Frequency distributions after 200 random walk/resampling iterations using Algorithm~\ref{alg1} under a periodic setting with period $q=2L=12$.
\textbf{Left}: Relative training and validation error as a function of resampling iterations.
\textbf{Center}: Histogram of projected frequency samples $v\cdot \omega$. 
\textbf{Right}: Histogram of projected frequency samples $v^\perp\cdot \omega$. The plots compare the learned (blue) and optimal (red) frequency density. 
The hyperparameter settings used here for Algorithm~\ref{alg1} are the same as those in Figure~\ref{fig:f27}.}
\label{fig:f28}
\end{figure}

\begin{figure}[!htbp]
\centering
\begin{minipage}[t]{0.32\textwidth}
    \centering
    \includegraphics[width=\linewidth]
    {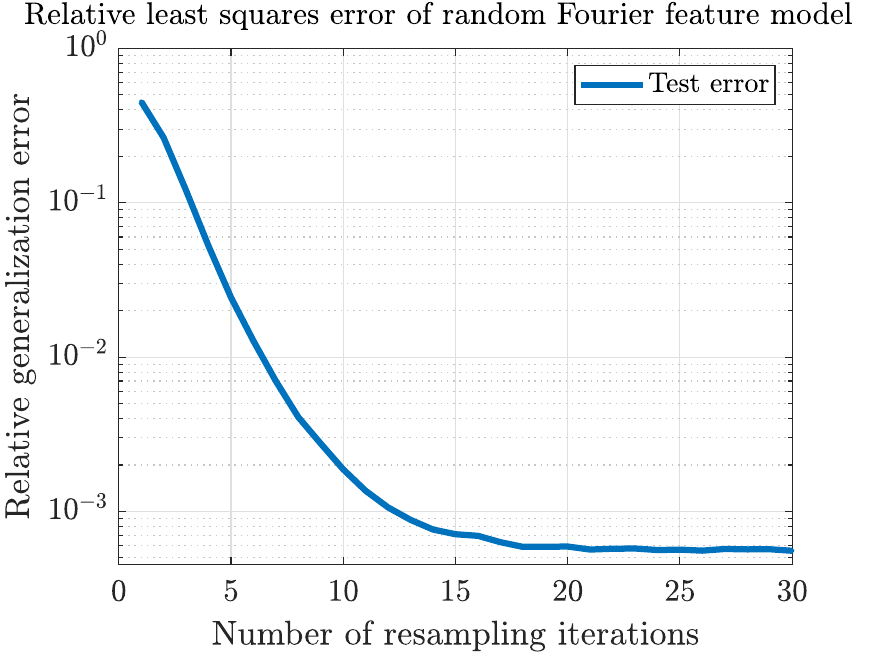}
\end{minipage}%
\hfill
\begin{minipage}[t]{0.32\textwidth}
    \centering
    \includegraphics[width=\linewidth]{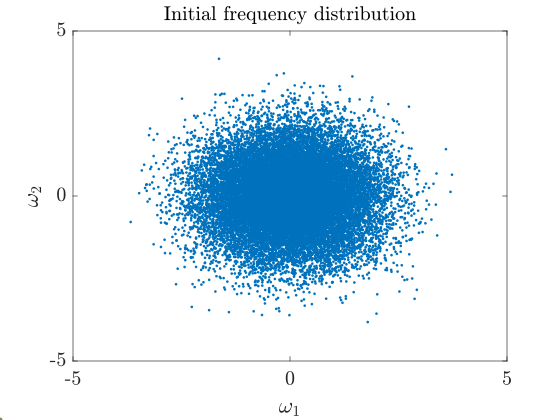}
\end{minipage}%
\hfill
\begin{minipage}[t]{0.32\textwidth}
    \centering
    \includegraphics[width=\linewidth]{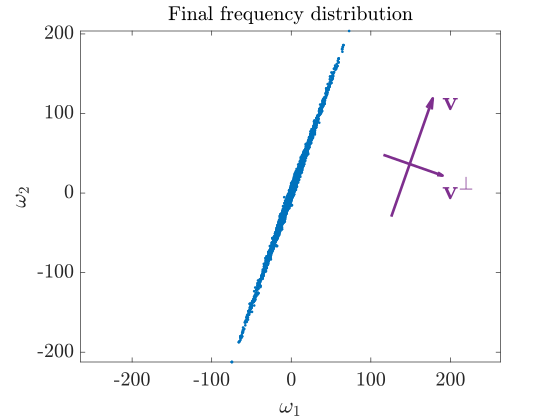}
\end{minipage}
\caption{Results for target function $f(x) = e^{-|v\cdot x|/a}e^{-|x|^2/2}$ with parameter $a=0.1$ in dimension $d=2$, using Algorithm~\ref{alg4} under non-periodic setting. \textbf{Left}: Relative generalization error as a function of resampling iterations. \textbf{Center}: Initial frequency samples drawn from a standard normal distribution. \textbf{Right}: Final frequency samples after 30 resampling steps. The hyperparameter settings used here for Algorithm~\ref{alg4} are the same as those in Figure~\ref{fig:f29}, with an additional parameter $\hat{\epsilon}=1\times 10^{-3}$ added to the diagonal of the empirical covariance matrix of frequency samples for stablization.
}
\label{fig:Sec1_Alg3_1}
\end{figure}

\begin{figure}[!htbp]
\centering
\begin{minipage}[t]{0.32\textwidth}
    \centering 
    \includegraphics[width=\linewidth]{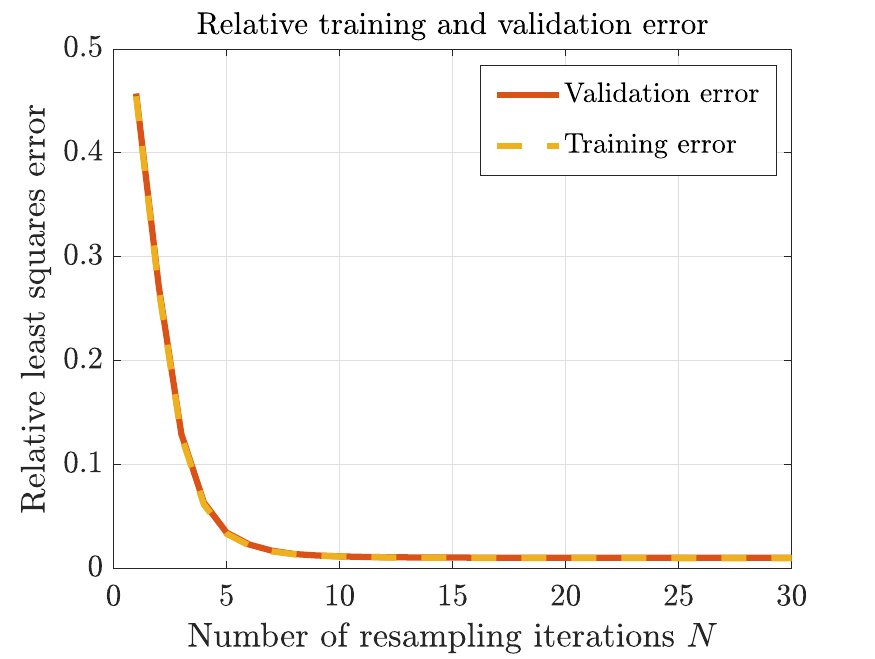}
\end{minipage}
\hfill
\begin{minipage}[t]{0.32\textwidth}
    \centering 
    \includegraphics[width=\linewidth]{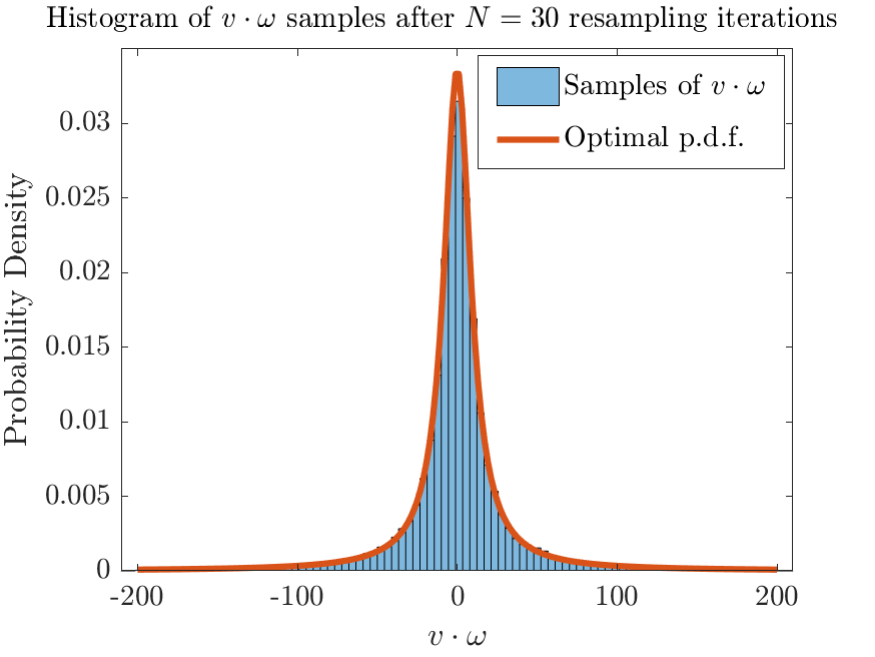}
\end{minipage}
\hfill
\begin{minipage}[t]{0.32\textwidth}
    \centering 
    \includegraphics[width=\linewidth]{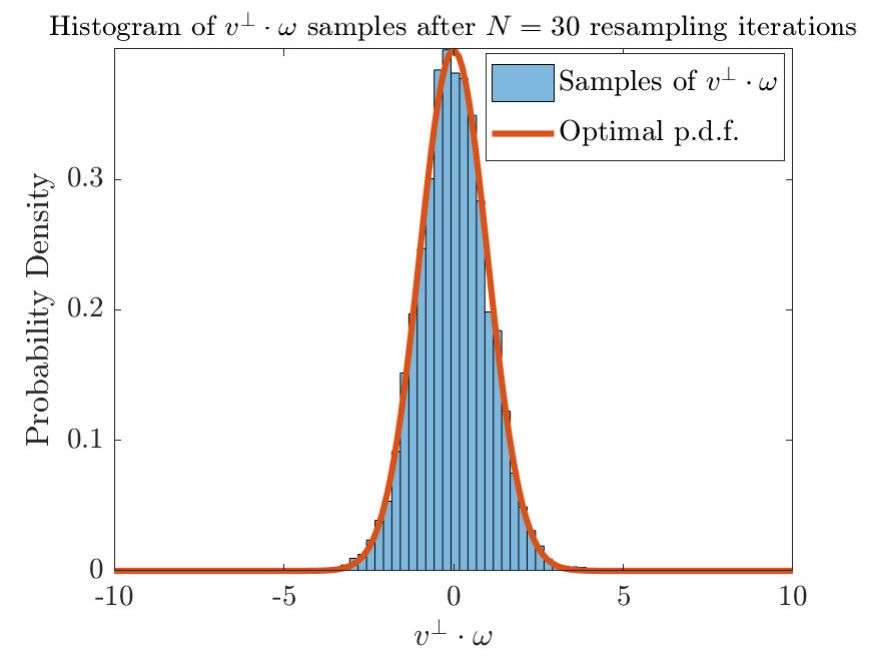}
\end{minipage}
\caption{Frequency distributions after 30 random walk/resampling iterations using Algorithm~\ref{alg4} under non-periodic setting.
\textbf{Left}: Relative training and validation error as a function of resampling iterations.
\textbf{Center}: Histogram of projected frequency samples $v\cdot \omega$. 
\textbf{Right}: Histogram of projected frequency samples $v^\perp\cdot \omega$. The plots compare the learned (blue) and optimal (red) frequency density. 
The hyperparameter settings used here for Algorithm~\ref{alg4} are the same as those in Figure~\ref{fig:Sec1_Alg3_1}.
}
\label{fig:Sec1_Alg3_2}
\end{figure}

\begin{figure}[!htbp]
\centering
\begin{minipage}[t]{0.32\textwidth}
    \centering
    \includegraphics[width=\linewidth]{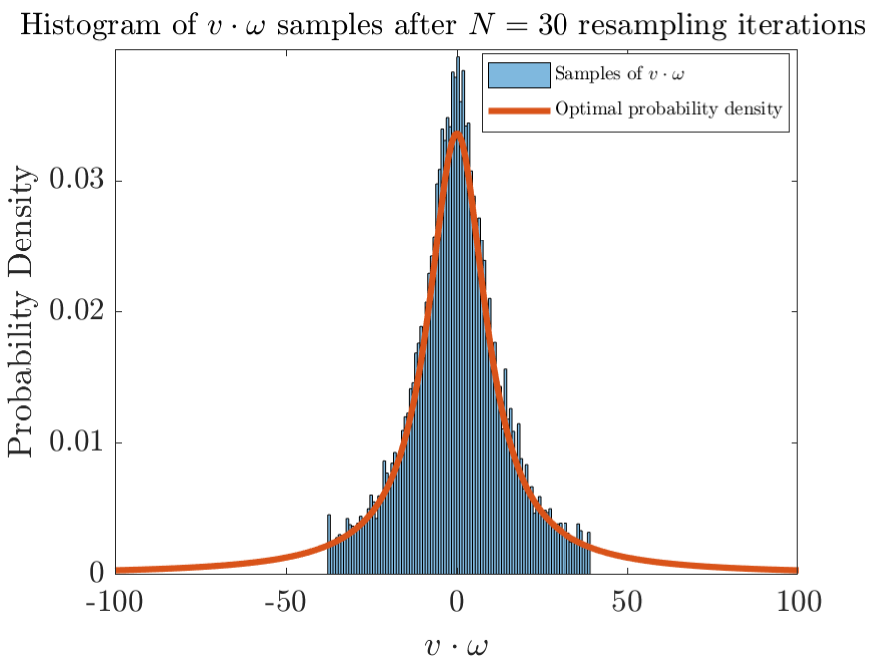}
\end{minipage}%
\hfill
\begin{minipage}[t]{0.32\textwidth}
    \centering
    \includegraphics[width=\linewidth]{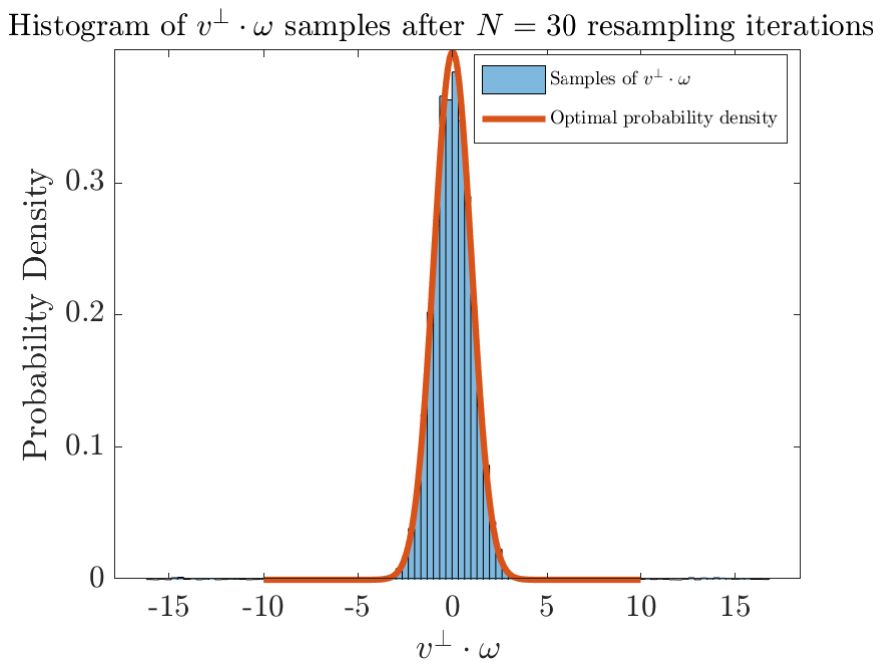}
\end{minipage}%
\hfill
\begin{minipage}[t]{0.32\textwidth}
    \centering
    \includegraphics[width=\linewidth]{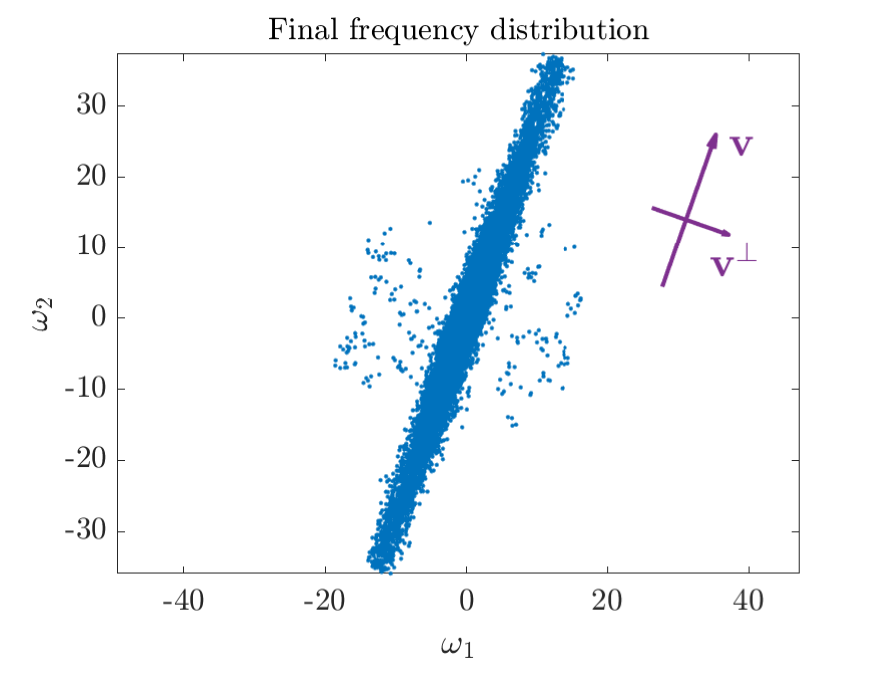}
\end{minipage}
\caption{Results from applying Algorithm \ref{alg_FT} with $30$ random walk/resampling iterations. \textbf{Left}: Histogram of projected frequency samples $v\cdot \omega$ (blue) overlaid with the optimal probability density function (red). \textbf{Center}: Histogram of projected frequency samples $v^\perp\cdot \omega$ (blue) overlaid with the optimal probability density function (red). \textbf{Right}: Final frequency samples after 30 resampling steps. The hyperparameter settings used here for Algorithm~\ref{alg_FT} are the same as those in Figure~\ref{fig:f29}.
}
\label{fig:f31}
\end{figure}

Figure \ref{fig:f30}  illustrates that a nearly optimal frequency
distribution is achieved after $200$ random walk and resampling iterations of Algorithm \ref{alg_FT}.
Figure \ref{fig:f28} shows a similar result for
 the related Algorithm \ref{alg1}, which uses a Fourier series representation instead of the Fourier transform representation in Algorithm \ref{alg_FT}. Figures~\ref{fig:f29} and~\ref{fig:f27} show that as the frequency samples converge toward the optimal distribution through resampling iterations, the  generalization error of the random Fourier feature model decreases.
Figures~\ref{fig:Sec1_Alg3_1} and~\ref{fig:Sec1_Alg3_2} further demonstrate that incorporating the empirical covariance matrix of the frequency samples into Algorithm~\ref{alg4}, a variant of Algorithm~\ref{alg_FT}, improves the efficiency of the resampling process, achieving comparable results with only 30 iterations. In contrast, Figure \ref{fig:f31} shows that
after $30$ resampling iterations with Algorithm~\ref{alg_FT}, the frequency distribution has improved from the initial standard normal distribution, but has not yet captured sufficiently high frequencies  to match the approximation quality observed in Figure~\ref{fig:Sec1_Alg3_2}.

For the numerical tests in Section~\ref{Sec_1_5}, Algorithms~\ref{alg_FT}, \ref{alg1}, and~\ref{alg4} are applied to the target function $f(x)=e^{-|v\cdot x|/a}\,e^{-|x|^2/2}$ (see Figure~\ref{fig:f_zoomed}).
%in Figure \ref{fig:f_zoomed}, 
This target function $f$ is associated with the frequency distribution $p_\ast$, which is fat-tailed in the direction $v$. 
To mimic real-world observations, Gaussian noise with mean zero and standard deviation $2^{-6}$ is added to the training data. The Conjugate Gradient method is used to solve the minimization problem \eqref{min3} with a relative tolerance of $10^{-3}$. % 

Additional numerical results are presented in Section~\ref{sec_num}, including experiments on the MNIST classification problem and a systematic investigation of hyperparameter sensitivity, such as the impact of the standard deviation of the noise variable and the random walk step size $\delta$. In Section~\ref{sec_num}, Algorithms~\ref{alg1} and~\ref{alg4}, which are modifications of Algorithm~\ref{alg_FT}, are described in greater detail and evaluated under various settings. These variants incorporate Fourier series expansions, iterative least squares solutions via the Conjugate Gradient method with zero-initialized frequencies, and adaptive random walks informed by the empirical covariance of the frequency samples. 
The zero-initialized frequencies are introduced to reduce the number of hyperparameters, thereby highlighting the influence of the random walk step size $\delta$ on model performance, while the empirical covariance enables adaptive random walks that accelerate convergence. If the function $f$ has compact support, the Fourier series and Fourier feature  setting based on the Fourier transform representation are directly related, using a sufficiently large periodic cell, as described in Section \ref{sec:FT_setting}.
 \subsection{Paper structure and previous results}
Theorem \ref{theorem1} proves that the resampling method achieves the optimal sampling rate constant $C_{p*}(f)$ asymptotically as $K\to \infty$.
The proof assumes that the function $f$ has bounded derivatives of order $\ell>d$,
the initial probability distribution $p$ is sufficiently close to the optimal distribution and that the amount of data is infinite, $J=\infty$.
In this case the minimal generalization error is smaller than $\mathcal O(K^{-1})$ in \eqref{eta2}, namely $\mathcal O(K^{1-2\ell/d })$, 
which is obtained  using the $K$ frequencies with the largest  Fourier coefficients as shown in Section \ref{small_GE}. 
The minimal generalization error and the optimal random Fourier feature error are close if $\ell/d\approx 1$.
The convergence in Theorem \ref{theorem1} is generalized to a finite amount of data, $J<\infty$, in Section \ref{sec:noise}. The results in Sections \ref{finite-K-section} and \ref{sec:noise} are local in the sense that the initial frequencies are assumed to be sufficiently close to the optimal. A method to handle the global problem of starting with initial frequencies equal to zero and to combine successive random walk and resampling iterations is studied in Section \ref{sec_iter}, based on Algorithm \ref{alg_FT} modified to use the periodic setting. 
Section \ref{other} relaxes the regularity conditions on $f$ to anisotropic decay rates including a few directions with slow decay. 

Computational methods to optimize the generalization error for the minimization %\eqref{goal} with respect to both amplitudes and frequencies are 
 \[ \min_{(\hat\beta,\omega)\in \mathbb C^K\times \rset^{Kd}}
 \big(\frac{1}{J}\sum_{j=1}^J|\beta(x_j)- y_j|^2 +\lambda\sum_{k=1}^K|\hat\beta_k|^2\big)\,
 \]
have been well studied based on the stochastic gradient method,  or variants thereof, %applied to the  optimization problem, 
see \cite{bach_book, understand_ml}. 
%to find both $(\hat\beta,\omega)\in \mathbb C^K\times\frac{\pi}{L}\mathbb Z^{dK}$. 
%This non convex minimization is a theoretical and computational challenge in general. 
The asymptotic behavior of the gradient method and the stochastic gradient method as the number of nodes $K$ and iterations tend to infinity is analyzed in \cite{ belkin2,Bach, mei}.

The alternative to use the data $\{(x_j,y_j)\, :\, j=1,\ldots,J\}$ to optimize the generalization error for the random Fourier feature problem \eqref{min} has received less study. The work \cite{LiY} formulates in their Algorithm 2 a gradient descent method applied to a loss function based on the sum of the mean squared error  for the neural network and an added penalty with respect to the kernel approximation. In \cite{ajpma} instead the explicit form of the optimal distribution \eqref{p*1} is used in the learning by replacing the unknown exact Fourier coefficients $|\hat f(\omega_k)|$ by the corresponding approximate amplitudes $|\hat\beta_k|$, with the aim to sample approximately optimal frequencies using a non linear Metropolis step based on the distribution $|\hat\beta_k|/\sum_{\ell=1}^K|\hat\beta_\ell|$. The work \cite{Aku2} applies resampling to the distribution $|\hat\beta_k|/\sum_{\ell=1}^K|\hat\beta_\ell|$ and show in numerical tests improvements compared to \cite{ajpma}, both in convergence rate and stability with respect to parameter perturbations.
In the present work we use the random Fourier feature algorithm including resampling in \cite{Aku2}, with the modifications that we include a cut-off for small amplitudes, adaptive random walk increments and iterative solutions of the least squares problem \eqref{min}. Adaptive random walks have been used before, for instance in the Metropolis algorithm \cite{sacs}.
%solved iteratively by the conjugate gradient method 

The main difference with the well studied convergence properties of resampling methods for particle filters and related problems
is that the sample average to evaluate here is based on a neural network function $\beta$ with coefficients $\hat\beta_k$ that
depend substantially on the sampling distribution for the frequencies.  Therefore the central limit theorem based techniques to obtain convergence for particle filters, see \cite{delmoral}, is not used here. Instead the proof of Theorem \ref{theorem1} directly estimates the rate constant $C_p(f)$,
using the cutoff and Chebyshev's inequality combined with the generalization error estimate \eqref{eta2} and a decay of the Fourier coefficients. 
Also in the case with finite amount of data including noise
the generalization error is minimized when the rate constant $C_p(f)$ is minimal. The proof
of this property of the generalization error
in Section \ref{sec_general} is new in the sense that it does not use the standard technique based on the Rademacher complexity \cite{bach_book, understand_ml,weinan_understand}. Instead   
we use Tikhonov regularizations $\lambda_1\|\hat\beta\|^2+\lambda_2\|\hat\beta\|^4$
and the orthogonality provided by the random Fourier features, which in a certain sense separates the data dependence and the amplitudes $\hat\beta$, although the amplitudes depend on the data. 

Resampling has also been applied to random Fourier features in \cite{Li}, using a leverage function based on the least squares kernel matrix. The leverage function in \cite{Li} does not depend on the data $y_j$ therefore the method in \cite{Li} is not adaptive in the sense of the random Fourier feature method with resampling presented in this study.

Sections \ref{Monte_Carlo} and \ref{sec_general} derive the generalization error for random Fourier features for infinite and finite amount of data, respectively. 
Section \ref{inf_nodes&data} analyses the resampling method for infinite number of nodes and data.
Sections \ref{finite-K-section} and \ref{sec:noise} prove
for the case of finite number of nodes $K$
in Theorem \ref{theorem1} (with infinite amount of data) and in Equation \eqref{conclude3} (with finite amount of data) that the resampling method achieves asymptotically optimal sampling, assuming that the initial sampling is sufficiently close to the optimal.  % in the case of finite number of nodes $K$ and infinite amount of data $J$.
Section \ref{sec_iter} studies the combination of random walk and resampling steps to iteratively spread initial random frequencies well from an initial distribution that is not close to the optimal.
Section \ref{sec_num} presents numerical results and the adaptive random walks. Section \ref{other} considers anisotropic decay rates for the Fourier coefficients. 

Our main wish for future related work is to  improve the mathematical precision %of the present work would be to formulate 
of the iteration procedure in Section \ref{sec_iter} and the extension to anisotropic decays in Section \ref{other}.
\section{The generalization error for infinite amount of data and the optimal distribution}
\label{Monte_Carlo}
\subsection{The generalization error}
Assume that
$f\in H^k(\mathbb T^d)$  for some $k>d/2$, then by Sobolev's inequality its  Fourier series converges pointwise, see Remark \ref{rem_sobolev}.
Using that $\nu_k$ are independent samples from the distribution $p$, Monte Carlo approximation of the
Fourier series representation yields
\[
\begin{split}
f(x)&=\sum_{n\in\mathbb Z^d}\hat f(\omega_n)e^{{\rm i}\omega_n\cdot x}\\
&\simeq \sum_{k=1}^K \frac{\hat f(\nu_k)}{Kp(\nu_k)}e^{{\rm i}\nu_k\cdot x}\\
&=:\sum_{k=1}^K \frac{\hat\zeta_k}{K}e^{{\rm i}\nu_k\cdot x}=:\zeta(x)\,,
\end{split}
\]
with the property of unbiasedness
\[
\begin{split}
&\mathbb E_\nu[\zeta(x)] = \mathbb E_\nu[K^{-1}\sum_{k=1}^K \frac{\hat f(\nu_k)e^{{\rm i}\nu_k\cdot x}}{p(\nu_k)}]\\
&=K^{-1}\sum_{k=1}^K \sum_{v_1\in \frac{\pi}{L}\mathbb Z^d}\ldots  \sum_{v_K\in \frac{\pi}{L}\mathbb Z^d}
\frac{\hat f(v_k)e^{{\rm i}v_k\cdot x}}{p(v_k)}p(v_1)\ldots p(v_K)\\
&=K^{-1}\sum_{k=1}^K\sum_{v_k\in \frac{\pi}{L}\mathbb Z^d}
\frac{\hat f(v_k)e^{{\rm i}v_k\cdot x}}{p(v_k)}
p(v_k)\sum_{v_1}p(v_1)\ldots
\underbrace{\sum_{v_{k-1}}p(v_{k-1})}_{=1}\sum_{v_{k+1}}p(v_{k+1})
\ldots \sum_{v_{K}}p(v_{K})\\
&=\sum_{v\in \frac{\pi}{L}\mathbb Z^d}\hat f(v)e^{{\rm i}v\cdot x}
%\\&
=f(x)\,.
\end{split}
\]
This no bias implies that the variance has the representation
\begin{align}
\label{var_a}
&\mathbb E_\nu[|\zeta(x)-f(x)|^2]
=
\mathbb E_\nu\Big[K^{-2}\sum_{k=1}^K\sum_{\ell=1}^K
\big(\hat\zeta_ke^{{\rm i}\nu_k\cdot x}-f(x)\big)^*
\underbrace{\big(\hat\zeta_\ell e^{{\rm i}\nu_\ell\cdot x}-f(x)\big)}_{=:g(\nu_\ell)}\Big]
\nonumber\\ \nonumber
&=K^{-2}\sum_{k=1}^K\sum_{\ell=1}^K
\sum_{v_1\in \frac{\pi}{L}\mathbb Z^d}\cdots \sum_{v_K\in \frac{\pi}{L}\mathbb Z^d}
g(v_k)^*g(v_\ell)\,p(v_1)\cdots p(v_K)
\\\nonumber
&=K^{-2}\sum_{k=1}^K
\sum_{v_1\in \frac{\pi}{L}\mathbb Z^d}\cdots \sum_{v_K\in \frac{\pi}{L}\mathbb Z^d}
|g(v_k)|^2\,p(v_1)\cdots p(v_K)
\\\nonumber
&\quad +K^{-2}\sum_{k=1}^K\sum_{\ell\ne k}
\sum_{v_1\in \frac{\pi}{L}\mathbb Z^d}\cdots \sum_{v_K\in \frac{\pi}{L}\mathbb Z^d}
g(v_k)^*g(v_\ell)\,p(v_1)\cdots p(v_K)
\\ \nonumber \displaybreak[3]
&=K^{-2}\sum_{k=1}^K
\Big(\sum_{v_k}|g(v_k)|^2p(v_k)\Big)\times
\underbrace{\sum_{v_1}\cdots \sum_{v_{k-1}}\sum_{v_{k+1}}\cdots \sum_{v_{K}}
p(v_1)\cdots p(v_{k-1})p(v_{k+1})\cdots p(v_K)}_{=1}
\\ 
&\quad +K^{-2}\sum_{k=1}^K\sum_{\ell\ne k}
\underbrace{\Big(\sum_{v_k}g(v_k)p(v_k)\Big)^*}_{=0}
\Big(\sum_{v_\ell}g(v_\ell)p(v_\ell)\Big)\times
\\\nonumber
&\qquad\times \sum_{v_1}\cdots \sum_{v_{k-1}}\sum_{v_{k+1}}\cdots
\sum_{v_{\ell-1}}\sum_{v_{\ell+1}}\cdots \sum_{v_{K}}
p(v_1)\cdots p(v_{k-1})p(v_{k+1})\cdots p(v_{\ell-1})p(v_{\ell+1})\cdots p(v_K)
\\\nonumber
&=K^{-1}\sum_{v\in \frac{\pi}{L}\mathbb Z^d}|g(v)|^2p(v)
\\\nonumber
&=K^{-1}\Big(\sum_{n\in\mathbb Z^d}\frac{|\hat f(\omega_n)|^2}{p(\omega_n)}-|f(x)|^2\Big)
=\mathcal O(K^{-1})\,,
\end{align}
for $x\in \mathbb T^d$.

The variance \eqref{var_a} establishes the generalization error in the case of infinite amount of data $J=\infty$
\begin{equation}\label{b_gen}
\begin{split}
\mathbb E_\nu\big[\min_{\hat\beta\in\mathbb C^K}\mathbb E_x[|\beta(x)-f(x)|^2]+\lambda\sum_{k=1}^K|\hat\beta_k|^2\big]
&\le \mathbb E_\nu\big[\mathbb E_x[|\zeta(x)-f(x)|^2]+\lambda\sum_{k=1}^K|K^{-1}\hat\zeta_k|^2\big]\\
&= \mathbb E_\nu\big[\mathbb E_x[|\zeta(x)-f(x)|^2]\big] +\lambda K^{-2}\sum_{k=1}^K\mathbb E_\nu[|\hat\zeta_k|^2]\\
&\le K^{-1}(1+\lambda)\sum_{n\in\mathbb Z^d} \frac{|\hat f(\omega_n)|^2}{p(\omega_n)}\,.
\end{split}
\end{equation}
\begin{remark}\label{rem_ae}{\rm
The condition $f\in H^{k}(\mathbb T^d)$ can be relaxed to $f\in L^2(\mathbb T^d)$
using that by Carleson's theorem \cite{carleson}
the Fourier series then converges almost everywhere, for sums of frequencies in increasing rectangular and general polygonal sets \cite{fefferman}. Consequently \eqref{b_gen} follows from \eqref{var_a}, which holds almost everywhere, and then integrating over the sampling density $\rho\in L^1(\mathbb T^d)$. If $f\in H^{k}(\mathbb T^d)$ then the result \eqref{b_gen} holds for any probability measure $\rho$.
}
\end{remark}

\subsection{The optimal distribution}\label{sec_opt_p}
The minimization problem 
\begin{equation}\label{min_p8}
\min_p\sum_{n\in\mathbb Z^d} \frac{|\hat f(\omega_n)|^2}{p(\omega_n)}
\end{equation}
over probability densities $p:\frac{\pi}{L}\mathbb Z^d\to (0,\infty)$ is a convex problem
since each term in 
\eqref{min_p8} is a  convex function
$p(\omega_n)\mapsto\frac{|\hat f(\omega_n)|^2}{p(\omega_n)}, \ p(\omega_n)>0$, depending on individual variables $p(\omega_n)$, and the   constraint is linear
$\sum_{n\in\mathbb Z^d}p(\omega_n)=1$. Hence \eqref{min_p8} has a solution.
The Lagrange formulation
\[
\sum_{n\in\mathbb Z^d} \frac{|\hat f(\omega_n)|^2}{p(\omega_n)} + \lambda_*\big(\sum_{n\in\mathbb Z^d}p(\omega_n)-1\big)
\]
with the Lagrange multiplier $\lambda_*\in\rset$
implies $\frac{|\hat f(\omega_n)|^2}{p(\omega_n)^2}=\lambda_*$ for all $n\in\mathbb Z^d$.
Therefore
the probability distribution 
\begin{equation}\label{p_opt}
p_*(\omega_n)
=\frac{|\hat f(\omega_n)|}{\sum_{n\in\mathbb Z^d}|\hat f(\omega_n)|}
\end{equation}
minimizes \eqref{min_p8} and the right hand sides in \eqref{var_a} and \eqref{b_gen}.

\subsection{The Fourier transform setting}\label{sec:FT_setting}
The more common setting with random Fourier feature frequency sampling in $\rset^d$, as in Algorithms \ref{alg_FT} and \ref{alg4}, can be related to the periodic formulation here as follows. Assume $f$ has compact support and that Assumption \ref{assump1} holds. Then the rate constant $C_p$ 
%and the estimates $\mathcal E(C_*,0)$ 
is bounded independent of the
period length $2L$. Therefore, the resampling results in Theorem \ref{theorem1} and Equation \eqref{conclude3} hold uniformly with respect to $L$ so that the resampling Algorithm \ref{alg1} can be based on any lattice $\frac{\pi}{L}\mathbb Z^d$, where $L\to\infty$ corresponds to the Fourier transform setting.

\section{Infinite amount of nodes and data}
\label{inf_nodes&data}
This section shows that the resampling method provides the optimal sampling in the case of infinite number of nodes, $K=\infty$, and infinite amount of data, $J=\infty$.
In the limit with $J=\infty$ and $K=\infty$, the generalization error estimate \eqref{eta2} implies  that
\[
\beta(x)=f(x)=\sum_{n\in \mathbb Z^d} \hat f(\omega_n)e^{{\rm i}\omega_n\cdot x}\,, \quad \mbox{almost surely with respect to $\nu$ and $x$}\,,
\]
so that
\[
\sum_{k=1}^K \hat\beta_k e^{{\rm i}\nu_k\cdot x}
=\sum_{n\in \mathbb Z^d} \hat f(\omega_n)e^{{\rm i}\omega_n\cdot x}\,,
\]
which implies
\begin{equation}\label{f_b}
\hat f(\omega_n)=\frac{1}{(2L)^d} \sum_{k=1}^K\hat\beta_k\int_{\mathbb T^d}
e^{{\rm i}(\nu_k-\omega_n)\cdot x}{\rm d}x
=\sum_{\{k:\,\nu_k=\omega_n\}}\hat\beta_k=\bar\beta_n\,.
\end{equation}

The  resampling distribution $
\bar p_n = 
\frac{|\bar\beta_n|}{\sum_{\{n:|\bar\beta_n|\ge\epsilon\}} |\bar\beta_n|}$
%\frac{|\bar\beta_n|}{\sum_{n\in\mathbb Z^d}|\bar\beta_n|}$ 
in \eqref{p_resam}
 is by \eqref{f_b} for $\epsilon\to 0+$ equal to the optimal distribution
\[
p_*(\omega_n)=\frac{|\hat f(\omega_n)|}{\sum_{n\in \mathbb Z^d}|\hat f(\omega_n)|}\,,
\]
derived in \eqref{p_opt} for  $J=\infty$.
\section{Finite amount of nodes  and the cutoff procedure}
\label{finite-K-section}
This section formulates and proves Theorem \ref{theorem1} on near optimal sampling for the resampling method with a finite number of nodes $K$, in the following steps, presented in the next four subsections:
\begin{itemize}
\item[-] error representation for $\beta-f$,
\item[-] assumptions on $\hat f$,
\item[-] cut off for small amplitudes,
\item[-] generalization error for the resampling.
\end{itemize}

\subsection{Error representation for \texorpdfstring{$\beta-f$}{beta-f}}
In the case of finite $K$ and infinite amount of data, $J=\infty$, the training  and testing errors become the same
and we have the generalization error estimate
\[
\beta(x;\nu)= f(x)+\eta(x,\nu)\,,
\]
where
\begin{equation}\label{eta}
\mathbb E[|\eta|^2 +\lambda\sum_{k=1}^K|\hat\beta_k|^2]\le K^{-1}(1+\lambda)\sum_{n\in\mathbb Z^d} \frac{|\hat f(\omega_n)|^2}{p(\omega_n)}\,,  
\end{equation}
is derived in \eqref{b_gen}.
The optimal distribution that minimizes the rate constant
\[
C_p(f):=\sum_{n\in\mathbb Z^d} \frac{|\hat f(\omega_n)|^2}{p(\omega_n)}
\]
is
$p_*(\omega_n)=\frac{|\hat f(\omega_n)|}{\sum_{n\in\mathbb Z^d}|\hat f(\omega_n)|}$,
which yields 
\begin{equation}\label{Cp}
C_{p_*}(f)=(\sum_{n\in\mathbb Z^d}|\hat f(\omega_n)|)^2\,.
\end{equation}
We obtain for $\bar\beta_n=\sum_{\{ k:\nu_k=\omega_n\}}\hat\beta_k$ the Fourier error representation
\begin{equation}\label{bfn}
 \bar\beta_n = \hat f(\omega_n) + (2L)^{-d}\int_{\mathbb T^d}\eta(x,\nu)e^{-{\rm i}\omega_n\cdot x}{\rm d}x\,,   
\end{equation}
by replacing $f$ by $f+\eta$ in \eqref{f_b}. The coefficient
\[
(2L)^{-d}\int_{\mathbb T^d}\eta(x,\nu)e^{-{\rm i}\omega_n\cdot x}{\rm d}x=:\hat\eta_n
\]
satisfies by \eqref{bfn} 
\[
\big||\bar\beta_n|-|\hat f(\omega_n)|\big|\le |\hat\eta_n|\,,
\]
which combined with the generalization error estimate \eqref{eta} 
establishes  the convergence rates
\begin{equation}\label{gen}
\begin{split}
\mathbb E[\big(|\bar\beta_n|-|\hat f(\omega_n)|\big)^2]
&\le \mathbb E[|\hat\eta_n|^2]\\
&\le (2L)^{-d}\mathbb E[\int_{\mathbb T^d}|\eta(x,\nu)|^2 {\rm d}x]\\
&\le K^{-1}(1+\lambda) C_p(f)\,, % \mathcal O(K^{-1})\,.
\end{split}
\end{equation}
and
\begin{equation}\label{gen2}
\begin{split}
\mathbb E[\big||\bar\beta_n|-|\hat f(\omega_n)|\big|]
&\le \mathbb E[|\hat\eta_n|]\\
&\le (\mathbb E[|\hat\eta_n|^2])^{1/2}\\
%&\le \underbrace{(2L)^{-d}\mathbb E[\int_{\mathbb T^d}|\eta(x,\nu)|^2 {\rm d}x]}_{=:\|\eta\|^2}\\
&\le \big(K^{-1}(1+\lambda) C_p(f)\big)^{1/2}\,.
\end{split}
\end{equation}

\subsection{Assumptions on \texorpdfstring{$\hat f $}{\hat f}}
\begin{assumption}[Regularity class for $f$]\label{assu:decay}
Assume that there are positive constants $c\,,C$ and $\ell>d$ such that
the function $f$ has Fourier coefficients with the decay 
\begin{equation}\label{decay}
(\frac{c}{|n|})^\ell\le |\hat f(\omega_n)|\le (\frac{C}{|n|})^\ell\,,
 \mbox{ for $n\in\mathbb Z^d\setminus 0$ \, where $|n|:=(n_1^2+\ldots n_d^2)^{1/2}$\,.   }
\end{equation}
%which holds for periodic functions with $L^2-$bounded derivatives of order $\ell $.
\end{assumption}

Now we use polar coordinates to obtain  an estimate for the high frequency content of $f,$ namely
\begin{equation}\label{f_N_decay}
\begin{split}
\sum_{|n|\ge N}|\hat f(\omega_n)|^2
&\le 
\min\big(\sigma_d\int_{N-\sqrt d}^\infty C^{2\ell}|n|^{-2\ell +d-1}{\rm d}|n|,\|\hat f\|^2_{\ell^2} \big)\\
&\le 
\min\big(\frac{\sigma_dC^{2\ell} }{(N/2 )^{2\ell -d}(2\ell-d)},\|\hat f\|^2_{\ell^2}\big)\,,
\end{split}
\end{equation}
where $\|\hat f\|^2_{\ell^2}:=\sum_{n\in\mathbb Z^d}|\hat f(\omega_n)|^2 $,  $\sigma_d=d(\pi)^{d/2}/\Gamma(\frac{d}{2} +1)$ is the surface area of the unit ball in dimension $d$ and $\Gamma$ is the Gamma function.
%To bound the right hand side we have also used the left inequality in \eqref{decay}.
%assume that the Fourier coefficients of $f$ also have a lower %bound given by
%\begin{equation} \label{lower}
%(\frac{c}{|n|})^\ell\le |\hat f(\omega_n)|\,,
%\end{equation}
%for  some $c\le C$.
Furthermore, with this bound and the left inequality in \eqref{decay} we obtain for any $\theta>0$ the following cutoff estimates,
\begin{equation}\label{theta}
\begin{split}
&\sum_{\{n:|\hat f(\omega_n)|\le \theta\}}|\hat f(\omega_n)|^2\le \min\big(\frac{\sigma_dC^{2\ell}}{(c/2)^{2\ell-d}(2\ell-d)(1-d^{1/2}\theta^{1/\ell}c^{-1})^{2\ell-d}} \theta^{2-d/\ell},\|\hat f\|^2_{\ell^2}\big)\,,\\
&\sum_{\{n:|\hat f(\omega_n)|\le \theta\}}|\hat f(\omega_n)|\le \min\big(\frac{\sigma_dC^{\ell}}{(c/2)^{\ell-d}(\ell-d)(1-d^{1/2}\theta^{1/\ell}c^{-1})^{\ell-d}} \theta^{1-d/\ell},
\underbrace{\sum_{n\in\mathbb Z^d}|\hat f(\omega_n)|}_{=:\|\hat f\|_{\ell^1}}\big)\,.
\end{split}
\end{equation}

Later, for instance in Theorem \ref{theorem1}, we will use the cutoff $\epsilon$, satisfying
\begin{equation}\label{epsilon}
c_1\le \epsilon K^{1/(3-d/\ell)}\le C_1\,,
\end{equation}
for constants $c_1$ and $C_1$,
in the resampling distribution  $\frac{|\bar\beta_n|}{\sum_{\{n:|\bar \beta_n|\ge\epsilon}|\bar\beta_n|}$, where this precise choice  is motivated below by optimizing error estimates.
We will also use a domain \[\{n\in\mathbb Z^d: |n|<N_*\}\] for  resolved Fourier modes, connected with the cutoff by Assumption \ref{assu:decay}, taking
\[
\frac{\epsilon}{2}=(\frac{C}{{N_*}})^\ell\,.
\]
Note that for the resolved Fourier modes, we have the following counting estimate,
\[
\sum_{\{n\in\mathbb Z^d: |n|+\sqrt d<N_*\}}1 \le \frac{\sigma_d}{d} N_*^d=\frac{\sigma_d}{d} \frac{2^{d/\ell}C^d}{\epsilon^{d/\ell}}\,,
%=\frac{2^{d/\ell}\sigma_d}{d} \frac{C^d\gamma^{-d/\ell}}{\kappa_f^{d/(2\ell)}}\,,
%=\mathcal O(K^{d/(2\ell)})\,,
\]
and by the assumption $\ell>d$ %we have for $\gamma=K^{1/2-d/\ell}\gg 1$
% and  for the given $\gamma\gg 1$ we choose $K$ sufficiently large so that 
 \begin{equation}\label{N_bound}
 \frac{\sigma_d}{d} N_*^d = \mathcal O(K^{\frac{d/\ell}{3-d/\ell}})\,, %\gamma^{-d/\ell}\mathcal O(K^{d/(2\ell)})\,.
 \end{equation}
 for the choice \eqref{epsilon}.
 We will see that this  condition on resolving the  kept low frequency modes is required for the resampling neural network  to accurately approximate the
 optimal sampling distribution  and in particular
 we will use the number
of basis functions, $\sum_{|\hat f(\omega_n)|\ge\epsilon} 1$, related to the cutoff. 
We have by \eqref{N_bound} that %\eqref{epsilon} and \eqref{decay} that
\begin{equation}\label{N_*}
\sum_{|\hat f(\omega_n)|\ge\epsilon/2} 1 \le 
\sigma_d \frac{N_*^d}{d} =\frac{2^{d/\ell}\sigma_d}{d} \frac{C^d}{\epsilon^{d/\ell}}=\mathcal O(K^{\frac{d/\ell}{3-d/\ell}})\,.
%\gamma^{-d/\ell}\mathcal O(K^{d/(2\ell)})\,.
%\mathcal O(K^{d/(2\ell)})\,.
%=\gamma^{-1/4}\mathcal O(K^{1/2})\,.
\end{equation}
%and for the cutoff at $|\bar\beta_n|<\epsilon$ and $\epsilon = K^{-1/2}$
\begin{remark}\label{rem_sobolev}{\rm
The assumption $\ell>d$ implies by Sobolev's inequality that the Fourier series for $f$ converges in $L^\infty(\mathbb T^d)$ as follows.
By Sobolev's inequality there is a constant $C$ such that
\[
\|v\|_{L^\infty(\mathbb T^d)}\le C\|v\|_{H^{k}(\mathbb T^d)}
\]
 holds for $k>d/2$  and $v\in \mathcal C(\mathbb T^d)$, where $H^{k}(\mathbb T^d)$ is the Sobolev space of functions with $k$ derivatives in $L^2(\mathbb T^d)$.
 Sobolev's inequality is proved  using the Fourier series representation 
 \[v(x)=\sum_{n\in \mathbb Z^d} \frac{\hat v(\omega_n)(1+|\omega_n|^k)}{(1+|\omega_n|^k)} e^{{\rm i}\omega_n\cdot x}\] and applying Cauchy's inequality for differentiable functions $v$ that are dense in $H^{k}(\mathbb T^d)$, see \cite{rauch}.
 Convergence in this Sobolev space follows as in \eqref{f_N_decay}, using Parseval's identity provided $2(\ell-k)-d>0$. 
Consequently the convergence of the Fourier series requires 
$
\ell>d/2+d/2=d \,.
$
}
\end{remark}
To approximately sample from the optimal distribution $p_*$ is advantageous, since the corresponding approximation error $K^{-1}C_{p_*}\ll K^{-1}C_{p}$ can be much smaller
than a random Fourier feature approximation based on a distribution $p$ with faster frequency decay than $\hat f$. 
%standard normal distributed 
%so that $C_p\gg C_{p_*}$ implying large approximation error, $K^{-1}C_p\gg K^{-1}C_{p_*}$.
%compared to using the optimal distribution $p_*$. 
We note that although the decay suggested for the tail of $\hat f$ in Assumption \ref{assu:decay} is sufficiently fast, the corresponding optimal constant $C_{p_*}(f)$ can be arbitrarily large if $C$ is large. The size of the constant $C$ matters, as $C_{p_*}=\mathcal O(C^\ell)$ can also become large for $\ell>d\gg 1$. The anisotropic setting with a few directions, with large constants $C$
having small decay powers $\ell$, %combined with possible slow decay with small powers $\ell$, 
and several remaining directions with small constants $C$ are discussed in Section \ref{other}.

\subsection{Cutoff for small amplitudes and resampling}\label{resamp_sec}
We will study the generalization error for the random feature problem based on a nonstandard choice of sampling distribution.

Recall that the optimal frequency sampling distribution $p_*$ is proportional to $|\hat f|$, so ideally, given $p$ we would like to compute amplitudes $|\bar \beta_n|$ that are accurate approximations for the Fourier amplitudes of $f$ and then use them to sample the frequencies in a nearly optimal way. However, since the total number of frequencies $K$ is finite, we will have a poor resolution of the amplitudes corresponding to high frequency. Consequently, we will not be able to trust these high-frequency amplitudes to guide our sampling in any meaningful way, so their information will be discarded by the sampling cutoff procedure described below. Thus, in the absence of further information on the high frequencies, we use the given sampling distribution $p$ to sample the remaining frequencies.

This work applies a conditional stochastic resampling distribution over the set of frequencies with the following structure:

{\it Input:} A probability distribution $p:\frac{\pi}{L}\mathbb Z^d\to [0,1]$,  an amplitude cutoff value $\epsilon>0$, and a cutoff ratio $q_\epsilon$ satisfying $K^{-1}\ll q_\epsilon\ll 1$.
%$0<q_\epsilon<1$, 

{\it First step:}
Sample  $\{\nu_k\}_{k=1}^K$ frequencies independent distributed from $p$ and solve the regularized least squares problem \eqref{min} for the approximation of $f$ using the sampled frequencies $\nu$, yielding the corresponding approximation $\beta$ with amplitudes $\{\hat\beta_k(\nu)\}$ and $\{\bar\beta_n(\nu)\},$ according to %\eqref{min} and 
\eqref{abs}.

{\it Second step:} Conditioned on $\nu = (\nu_1,\ldots,\nu_K)$, compute two numbers adding to $K$, namely $\bar K:=\lfloor K(1-q_\epsilon)\rfloor$ and $\tilde K:=\lceil Kq_\epsilon\rceil$,  % $\bar K +\tilde K =K$ 
and the resampled conditional probability function $\bar p(\cdot |\nu,\epsilon)$
defined by
\[
\bar p(n|\nu, \epsilon):= \frac{|\bar\beta_n|}{\sum_{\{n:|\bar\beta_n|\ge\epsilon\}} |\bar\beta_n|}\,,
\]

{\it Third step:} 
Sample  $\{\bar \nu_k\}_{k=1}^{\bar K}$ independent frequencies  from $\bar p(\cdot|\nu, \epsilon)$
and sample  independently $\{\tilde \nu_k\}_{k=1}^{\tilde K}$ independent frequencies from $p$. %$\tilde p_\epsilon(\cdot|\nu)$ 

As a result, our resampling random feature approximation solves the regularized least squares problem \eqref{min_alp_b} using the frequencies 
\begin{equation}\label{eq:omegadef}
\bar\omega := (\bar \nu,\tilde\nu)
\end{equation} 
corresponding to a single sample from 
 \begin{equation}\label{eq:barpdef}
    \check p(\bar\nu,\tilde\nu,\nu):= \left(\otimes_{n=1}^{\bar K}\bar p(\bar\nu_n|\nu,\epsilon)\otimes \otimes_{m=1}^{\tilde K} p(\tilde \nu_m) \right) \otimes_{k=1}^K p(\nu_k),
\end{equation}   
using $\nu$ as an auxiliary, seed sample.
%Observe that in \eqref{eq:barpdef} the components of $\nu$ are iid, while the components of $\bar \nu$ and $\tilde \nu$ are only conditionally iid.
%
Along this line, we make the following definition for a random Fourier feature method, based on the sampled frequencies $\omega$, 
%letting \[v(x;\omega)=\sum_{k=1}^{ K}\hat v_ke^{{\rm i} \omega_k\cdot x}\]
%and then solve, based on
 and $J$ data points.
 \begin{definition}\label{assump1}
Assume $\beta(\cdot;\nu)$ solves the random Fourier feature problem \eqref{min}, using $K$ independent frequencies $\{\nu_k\}_{k=1}^K$ sampled from $p$, % with iid frequencies $\nu$ 
and that \[\alpha(x;\bar\omega) := \sum_{k=1}^K\hat{\alpha}_ke^{{\rm i}\bar\omega_k\cdot x}\] has amplitudes that solve, given $J$ data points,
\begin{equation}\label{min_alp_b}
\begin{split}
\hat \alpha :=&\arg\min_{\hat{v}\in\mathbb C^K}\big(J^{-1}\sum_{j=1}^J|\sum_{k=1}^K\hat{v}_ke^{{\rm i}\bar\omega_k\cdot x_j}- f(x_j)|^2 +\lambda\sum_{k=1}^K|\hat{v}_k|^2\big)\,,%\\
%&v(x;\omega) = \sum_{k=1}^K\hat{v}_ke^{{\rm i}\omega_k\cdot x}\,.
\end{split}
\end{equation}
using frequencies $\bar\omega$ defined in \eqref{eq:omegadef}. 
\end{definition}
The proof of Theorem \ref{theorem1} is based on splitting $f$ into two parts
\begin{equation}\label{f_split}
f(x)=\sum_{\{n:|\bar\beta_n|>\epsilon\}} \hat f(\omega_n)e^{{\rm i} \omega_n\cdot x}
+\sum_{\{n:|\bar\beta_n|\le \epsilon\}} \hat f(\omega_n)e^{{\rm i} \omega_n\cdot x}
=: f_\epsilon(x) + (f- f_\epsilon)(x)\,.
\end{equation}
Note that the above decomposition is inherently random. Moreover, in the particular case of infinite data, namely $J=\infty$, it is a deterministic construction given %function of 
the independent sampled frequencies $\nu$ from $p$. % which are $p-$iid. 
The splitting \eqref{f_split} is motivated as follows.

The goal is to compare the resulting rate constant $C_{\check p}$ to the optimal $C_{p_*}$ and to conclude that the use of \eqref{eq:barpdef} is asymptotically optimal as $K\to \infty$. % cf. Theorem \ref{theorem1}.
%to avoid the coupling of iterated sampling densities  $|\bar\beta_n|/\sum_n|\bar\beta_n|$.
To obtain a convergence rate for the feature sampling distribution requires  to estimate %in some sense to estimate
\[
\mathbb E[\big|\sum_n|\bar\beta_n|-\sum_n|\hat f(\omega_n)|\big|]\,,\]
and it is useful to cut off the sum, retaining only $\sum_{\{n:|\bar\beta_n|\ge\epsilon\}}|\bar\beta_n|$
to avoid summing over many high-frequency, small-amplitude terms with relatively large errors.

The generalization error estimate 
\begin{equation}\label{eta1}
\mathbb E[|\eta|]\le (\mathbb E[|\eta|^2])^{1/2}\le (\frac{(1+\lambda)C_p(f)}{K})^{1/2}=:\mathcal E_f^{1/2}%\mathcal O(K^{-1/2})
\end{equation}
%we have  by the generalization error estimate \eqref{gen} 
%for $|\bar\beta_n|>\epsilon$ 
and \eqref{gen2} imply that
\[
|\hat f(\omega_n)| - \mathcal E_f^{1/2}\le \mathbb E[|\bar\beta_n|] \le |\hat f(\omega_n)| + \mathcal E_f^{1/2}\,,
%\lesssim 2\mathbb E[|\hat f(\omega_n)|]
\]
which shows that small amplitudes can have a relatively large error contribution
so that the expected sum
$\mathbb E[\sum_{\{n:|\bar\beta_n|\le\epsilon\}}|\bar\beta_n|]$ 
% $\mathbb{E}\!\left[\sum_{n:\,\allowbreak |\bar\beta_n|\le \varepsilon} |\bar\beta_n|\right]$
may not approximate $\mathbb E[\sum_{\{n:|\bar\beta_n|\le\epsilon\}}|\hat f(\omega_n)|]$.
Therefore we make the splitting \eqref{f_split}, where the resampled frequencies $\bar\nu$, based on $\bar p(\cdot|\nu,\epsilon)$, are used to estimate the part $f_\epsilon$ while the frequencies $\tilde\nu$ sampled from $p$ are used to estimate the cut part $f-f_\epsilon.$ This approximation for the split representation of $f$ also determines the choice of the cutoff value $\epsilon$ and the cutoff ratio $q_\epsilon$ by properly balancing the resulting error terms. Namely, for the part approximating $f_\epsilon$, the resampled frequencies have a distribution close to the optimal, while the cut part $f-f_\epsilon$ is small due to the decay \eqref{decay} of the Fourier coefficients, so that only a few samples $Kq_\epsilon$ of the non-optimal density $p$ make the corresponding error contribution asymptotically negligible. 

\begin{remark}{\rm
In the numerical solution of the least squares problem \eqref{min} we obtain $\hat\beta_k$ and we need to determine also $\bar\beta_n=\sum_{\{k:\nu_k=\omega_n\}}\hat\beta_k$ to perform the cutoff. The naive idea to base this search for equal frequencies directly on $\nu_k$ seems computationally infeasible for large $K$. Instead, we sort $|\hat\beta_k|$ and identify the sets of equal amplitudes to determine the sets $\{k:\nu_k=\omega_n\}$ of equal frequencies. This sorting has a computational cost of the order $\mathcal O(K\log K)$, which is negligible compared to solving the least squares problem. } 
\end{remark}

The next step is to estimate the generalization error for the resampling method in Definition \ref{assump1}. %$\bar\alpha$. %$\alpha$ and $\alpha'$. %, which both use the sampling distribution $\bar p$.
%which happens in the resampling method. Namely,

\subsection{The generalization error for the resampling method}
In this section we consider the case of infinite amount of data $J=\infty$ and prove
\begin{theorem}\label{theorem1}
Let $f$ satisfy Assumption \ref{assu:decay}.
Suppose the neural networks $\beta(\cdot;\nu)$ and $\alpha(\cdot;\bar\omega)$ satisfy Definition \ref{assump1} with $\lambda$ a positive constant and $J=\infty$, there is a constant $C''$ such that
 \begin{equation}\label{C''}
 \sup_{n\in\mathbb Z^d}\frac{|\hat f(\omega_n)|}{p(\omega_n)}\le C''\,,
 \end{equation}
and the cutoff $\epsilon$ satisfies $c_1\le \epsilon K^{\frac{1}{3-d/\ell}}\le C_1$, for positive constants $c_1$ and $C_1$. %is $\epsilon=\mathcal{O}(K^{-\frac{1}{3-d/\ell}})$  %for some $\gamma\gg 1$
%the rate constant $C_p$ is finite,
Then the resampled random Fourier feature method for $\alpha(\cdot;\bar\omega)$, based on Definition \ref{assump1},
has for $c_2\le q_\epsilon K^{\frac{1-d/\ell}{6-2d/\ell}}\le C_2$, with positive constants $c_2$ and  $C_2$,
%q_\epsilon= %\mathcal{O}(K^{-\frac{1-d/\ell}{6-2d/\ell}})$ 
the asymptotically optimal generalization error  
\begin{equation}\label{conclude}
\begin{split}
\mathbb E_{\check p}\big[
%\min_{\hat{\bar \alpha}\in\mathbb C^K} 
\mathbb E_x[|\alpha(x;\bar\omega)-f(x)|^2]\big] & \le
\mathbb E_{\check p}\big[\min_{\hat v\in\mathbb C^{K}}\big(\mathbb E_x[|\sum_{k=1}^{ K}\hat v_ke^{{\rm i} \bar\omega_k\cdot x}- f(x)|^2] +\lambda\sum_{k=1}^{ K}|\hat v_k|^2\big)\big] \\
%\mathbb E\big[
%\min_{\hat{\bar \alpha}\in\mathbb C^K} 
%\mathbb E_x[|\alpha(x;\omega)-f(x)|^2]\big]
&\le  \frac{(1+\lambda)}{K}\big(\sum_{n\in\mathbb Z^d}|\hat{f}(\omega_n)|\big)^2
\big(1+C''\mathcal O(q_\epsilon^{1/2})\big) 
%\kappa_K(1+\gamma^{-1/4})(1+\gamma^{-1/2})\\
%\\&\quad 
 %+K^{-1}\mathcal O(q_\epsilon^{1/2})
 \,.\\
\end{split}
\end{equation}
\end{theorem}

We see that
%,  for functions $f$ with the fast decay \eqref{decay} of its tail Fourier coefficients and the sufficiently resolved low frequency modes \eqref{N_bound}, 
the generalization error
for the resampling method approaches the nearly optimal rate. That is, also if we start with a distribution $p$ for determining $\beta$
where \[
C_p\gg C_{p_*}=(\sum_{n\in\mathbb Z^d}|\hat f(\omega_n)|)^2\,, \]
and $C''\gg 1$, the resampling asymptotically yields the  optimal rate constant 
\begin{equation}\label{C''est}
  C_{\check p}=\big(1+C''\mathcal O(q_\epsilon^{1/2})\big)C_{p_*}\,.
\end{equation}
% in this case.
\begin{remark}{\rm
   Condition \eqref{C''} means that the initial sampling  distribution needs to be sufficiently close to the optimal
and consequently Theorem \ref{theorem1} is a  local result. 
However, Assumption \ref{assu:decay} allows us to construct a sampling distribution that satisfies \eqref{C''}, that is, by taking
\[p(\omega_n) \propto (1+|n|)^{-\ell}\]
so \eqref{C''} is not restrictive when the decay rate $\ell$ is known as given in Assumption \ref{assu:decay}. Note in passing that the value of $\ell$ also plays a crucial role in determining acceptable ranges for $\epsilon$ and $q_\epsilon$ when applying Theorem \ref{theorem1}.
The question of finding an approximation to the optimal sampling distribution from an initial distribution not satisfying \eqref{C''}, that is $C''=\infty$, %where $C''q_\epsilon^{1/2}\gg 1$, 
is studied in Section \ref{sec_iter} by performing iterated resampling and random walk steps. The extension of \eqref{conclude} to finite amount of data is in \eqref{conclude3}.}
\end{remark}

The proof of the theorem is based on the following lemmas.
\begin{lemma}[Error split]\label{lemma1}
  For any $\delta>0$, 
 we have
 \begin{equation}\label{gen_split}
\begin{split}
&\mathbb E_{\check p}\big[\min_{\hat v\in\mathbb C^{K}}\big(\mathbb E_x[|\sum_{k=1}^{ K}\hat v_ke^{{\rm i} \bar\omega_k\cdot x}- f(x)|^2 +\lambda\sum_{k=1}^{ K}|\hat v_k|^2]\big)\big]\\
&\le (1+\delta)\mathbb E_{\check p}\big[\varphi_1(\bar \nu,\nu )\big]
 +(1+\delta^{-1}) \mathbb E_{\check p}\big[\varphi_2(\tilde \nu,\nu )\big]\,,\\
\end{split}
\end{equation}

with 
\[\varphi_1(\bar \nu,\nu ) := 
\min_{\hat v\in\mathbb C^{\bar K}}\big(\mathbb E_x[|\sum_{k=1}^{\bar K}\hat v_ke^{{\rm i} \bar\nu_k\cdot x}- f_\epsilon(x)|^2 +\lambda\sum_{k=1}^{\bar K}|\hat v_k|^2]\big)
\]
and
\[\varphi_2(\tilde \nu,\nu ) := \min_{\hat v\in\mathbb C^{\tilde K}}\big(\mathbb E_x[\sum_{k=1}^{\tilde K}\hat v_ke^{{\rm i} \tilde\nu_k\cdot x}- (f-f_\epsilon)(x)|^2 +\lambda\sum_{k=1}^{\tilde K}|\hat v_k|^2]\big)\,.
\]
\end{lemma}
\begin{proof}
Consider a generic $\hat v\in\mathbb C^{K}$ corresponding to frequencies $\bar\omega$ and split it accordingly to the partition of frequencies \eqref{eq:omegadef} into 
$\hat v = (\hat u, \hat z)$, with $\hat u\in\mathbb C^{\bar K}$ and $\hat z\in\mathbb C^{\tilde K}$.
Then by the arithmetic-geometric mean inequality $2|ab|\le a^2/\delta+b^2\delta,$ %Young's inequality implies, 
for any $\delta>0,$ and using the splitting of $f$ \eqref{f_split},
\[
\begin{split}
&\mathbb E_x[|\sum_{k=1}^{ K}\hat v_ke^{{\rm i} \bar\omega_k\cdot x}- f(x)|^2 +\lambda\sum_{k=1}^{ K}|\hat v_k|^2]\\
& \le
(1+\delta) \mathbb E_x[|\sum_{k=1}^{ \bar K}\hat u_ke^{{\rm i} \bar \nu_k\cdot x}- f_\epsilon(x)|^2 +\lambda\sum_{k=1}^{ \bar K}|\hat u_k|^2] \\
& \quad +
(1+\delta^{-1})\mathbb E_x[|\sum_{k=1}^{ \tilde K}\hat z_ke^{{\rm i} \tilde \nu_k\cdot x}- (f-f_\epsilon)(x)|^2 +\lambda\sum_{k=1}^{\tilde K}|\hat z_k|^2]\,.
\\
\end{split}
\]
%To conclude, observe that the minimum over $\hat v$ is bounded by the sum of the two minima over $\hat u_k$ and $\hat z_k$.
Take minimum in the above over $\hat v$ and observe that the right hand side of the inequality splits into the sum of two minima. 
Then use the expected value with respect to $\check p$ to obtain \eqref{gen_split}.

%    Prove it using Young's inequality and carefully invoking the conditional sampling procedure.
\end{proof}
\begin{lemma}[Estimate for $f_\epsilon$]\label{lemma2}
We have
\[
\begin{split}
\mathbb E_{\check p}\big[\varphi_1(\bar \nu,\nu ) \big] &=
\mathbb E_{\check p}\Big[\ \min_{\hat v\in\mathbb C^{\bar K}}\big(\mathbb E_x[|\sum_{k=1}^{\bar K}\hat v_ke^{{\rm i} \bar\nu_k\cdot x}- f_\epsilon(x)|^2 +\lambda\sum_{k=1}^{\bar K}|\hat v_k|^2]\big)\Big] \\
&\le \big(K(1-q_\epsilon)\big)^{-1}(1+\lambda)\Big[\big(\sum_{n\in \mathbb Z^d}|\hat f(\omega_n)|\big)^2\big(1+(C'')^{1/2}\mathcal O(K^{-\frac{1-d/\ell}{6-2d/\ell}})\big)\\
&\quad +\sum_{n\in \mathbb Z^d}|\hat f(\omega_n)|\big(
(C'')^{1/2}\mathcal O(K^{-\frac{3(1-d/\ell)}{6-2d/\ell}})
+C''\mathcal O(\lambda^{-1/2}K^{-\frac{1-d/\ell}{6-2d/\ell}})\big)\Big]\,.
\end{split}
\]
\end{lemma}
\begin{proof}
We note that, by the construction of the sample measure $\check p$ and conditional on the sampled frequencies $\nu$ which determine $f_\epsilon$, we can apply the generalization error estimate for the random features based on $\bar \nu,$ namely
\begin{equation}\label{namely}
\begin{aligned}
\mathbb E_{\check p}\Big[\varphi_1(\bar \nu,\nu ) \Big]  &\le \big(K(1-q_\epsilon)\big)^{-1}(1+\lambda) \mathbb E_{ p}\big[\sum_{\{n: |\bar\beta_n| > \epsilon\}} \frac{|\hat{ f_\epsilon}(\omega_n)|^2}{\bar p(n|\nu,\epsilon)}\big]\\
&= \big(K(1-q_\epsilon)\big)^{-1}(1+\lambda) \mathbb E_p\big[\big(\sum_{\{n: |\bar\beta_n| > \epsilon\}} \frac{|\hat{f_\epsilon }(\omega_n)|^2}{|\bar\beta_n|}\big) \big(\sum_{\{n: |\bar\beta_n| > \epsilon\}}|\bar\beta_n| \big)\big]\,.
\end{aligned}
\end{equation}
%% Here goes the rest of the proof
The Fourier coefficient for the perturbation $\eta=\beta-f$ has the bound
\[
|\hat\eta_n|=(2L)^{-d}|\int_{\mathbb T^d} \eta(x,\nu)e^{{\rm i}\omega_n\cdot x}{\rm d}x|
\le (2L)^{-d}\int_{\mathbb T^d} |\eta(x,\nu)|{\rm d}x=:\|\eta\|\,,
\]
and in  the right hand side  of \eqref{namely} we make the splitting
\begin{equation}\label{splitt1}
\begin{split}
\sum_{\{n: |\bar\beta_n| > \epsilon\}} \frac{|\hat{ f_\epsilon}(\omega_n)|^2}{\bar p_n}
&=\Big(\sum_{\{n: |\bar\beta_n| > \epsilon\}} \frac{|\hat{f_\epsilon }(\omega_n)|^2}{|\bar\beta_n|}\Big) \Big(\sum_{\{n: |\bar\beta_n| > \epsilon\}}|\bar\beta_n| \Big)\\
&= \big(\mathbf 1_{\|\eta\| < \epsilon/2} + \mathbf 1_{\|\eta\|  \ge \epsilon/2}\big)
\Big(\sum_{\{n: |\bar\beta_n| > \epsilon\}} \frac{|\hat{f_\epsilon }(\omega_n)|^2}{|\bar\beta_n|}\Big) \Big(\sum_{\{n: |\bar\beta_n| > \epsilon\}}|\bar\beta_n| \Big)
\,.
\end{split}
\end{equation}
For the first term in the right hand side of \eqref{splitt1} we have by \eqref{N_bound} and \eqref{eta1} 
\begin{equation}\label{bound21}
\begin{split}
&\mathbb E_p\big[
\mathbf 1_{\|\eta\| < \epsilon/2} 
\big(\sum_{\{n: |\bar\beta_n| > \epsilon\}} \frac{|\hat{f_\epsilon }(\omega_n)|^2}{|\bar\beta_n|}\big) \big(\sum_{\{n: |\bar\beta_n| > \epsilon\}}|\bar\beta_n| \big)
\big]\\
&\le 
\mathbb E_p\big[\mathbf 1_{\|\eta\| < \epsilon/2}  \sum_{\{n: |\bar\beta_n| > \epsilon\}} |\hat{f_\epsilon }(\omega_n)|(1+\frac{\|\eta\|}{\epsilon})
(\sum_{\{n: |\hat f(\omega_n)| > \epsilon/2\}} |\hat{f_\epsilon }(\omega_n)| +\|\eta\|\sum_{\{n: |\hat f(\omega_n)| > \epsilon/2\}} 1)\big]\\
&\le 
\mathbb E_p\big[\mathbf 1_{\|\eta\| < \epsilon/2}  \sum_{\{n: |\hat f(\omega_n)| > \epsilon/2\}} |\hat{f_\epsilon }(\omega_n)|(1+\frac{\|\eta\|}{\epsilon})
(\sum_{\{n: |\hat f(\omega_n)| > \epsilon/2\}} |\hat{f_\epsilon }(\omega_n)| +\|\eta\|\sum_{\{n: |\hat f(\omega_n)| > \epsilon/2\}} 1)\big]\\
&\le (\sum_{n\in \mathbb Z^d}|\hat f(\omega_n)|)^2
\mathbb E_p[\mathbf 1_{\|\eta\| < \epsilon/2} ( 1 +\frac{\|\eta\|}{\epsilon})]  
+ \big(\sum_{n\in \mathbb Z^d}|\hat f(\omega_n)| \big)\frac{2^{d/\ell}\sigma_d}{d} N_*^d\mathbb E_p[\mathbf 1_{\|\eta\| < \epsilon/2}\|\eta\|(1+\frac{\|\eta\|}{\epsilon})]\\
&\le \big(\sum_{n\in \mathbb Z^d}|\hat f(\omega_n)|\big)^2
\big(1+\mathcal O(\mathcal E_f^{1/2}\epsilon^{-1})\big)
+\sum_{n\in \mathbb Z^d}|\hat f(\omega_n)|\mathcal O(\mathcal E_f^{1/2}\epsilon^{-d/\ell})\\
&\le \big(\sum_{n\in \mathbb Z^d}|\hat f(\omega_n)|\big)^2\big(1+\sqrt{C''}\mathcal O(K^{-\frac{1-d/\ell}{6-2d/\ell}})\big)
+\sum_{n\in \mathbb Z^d}|\hat f(\omega_n)|\sqrt{C''}\mathcal O(K^{-\frac{3(1-d/\ell)}{6-2d/\ell}})\,,
\end{split}
\end{equation}
and by Chebyshev's inequality we obtain for the second term in the right hand side of \eqref{splitt1}
\begin{equation}\label{cheby}
\begin{split}
&\mathbb E_p\big[
\mathbf 1_{\|\eta\|  \ge \epsilon/2}
\big(\sum_{\{n: |\bar\beta_n| > \epsilon\}} \frac{|\hat{f_\epsilon }(\omega_n)|^2}{|\bar\beta_n|}\big) \big(\sum_{\{n: |\bar\beta_n| > \epsilon\}}|\bar\beta_n| \big)\big] \\
&\le \mathbb E_p\big[ \mathbf 1_{\|\eta\|  \ge \epsilon/2}
\sum_{\{n: |\bar\beta_n| > \epsilon\}} |\hat{f_\epsilon }(\omega_n)|(1+\frac{\|\eta\|}{\epsilon})\sum_{\{n: |\bar\beta_n| > \epsilon\}}|
\bar\beta_n| \big] \\
&\le
\sum_{n\in \mathbb Z^d}|\hat f(\omega_n)|\big(
\mathbb E_p[\sum_{k=1}^K|\hat\beta_k|\mathbf 1_{\|\eta\|  \ge \epsilon/2}]
+\mathbb E_p[\sum_{k=1}^K|\hat\beta_k|\frac{\|\eta\|}{\epsilon}]\big)\\
&\le \sum_{n\in \mathbb Z^d}|\hat f(\omega_n)|
\Big( \big(\mathbb E_p[(\sum_{k=1}^K|\hat\beta_k|)^2] \big)^{1/2}\big((\mathbb E_p[\mathbf 1_{\|\eta\|  \ge \epsilon/2}])^{1/2}
+(\mathbb E_p[\frac{\|\eta\|^2}{\epsilon^2}])^{1/2}\big)\Big)\\
&\le \sum_{n\in \mathbb Z^d}|\hat f(\omega_n)|
\Big(\big(K\mathbb E_p [\sum_{k=1}^K|\hat\beta_k|^2] \big)^{1/2}\big( \frac{2\mathcal E_{f}^{1/2}}{\epsilon }
+\frac{\mathcal E_{f}^{1/2}}{\epsilon}\big)\Big)\\
%&\le  \sum_{n\in \mathbb Z^d}|\hat f(\omega_n)|\mathcal O(\lambda^{-1/2}C''K^{-1/2}\epsilon^{-1})\\
&\le  \sum_{n\in \mathbb Z^d}|\hat f(\omega_n)|C''\mathcal O(\lambda^{-1/2}K^{-\frac{1-d/\ell}{6-2d/\ell}})\,,
\end{split}
\end{equation}
using also that 
$
\sum_{n\in\mathbb Z^d}|\bar\beta_n|=\sum_{k=1}^K|\hat \beta_k|
$
and $\mathbb E[\sum_{k=1}^K|\hat \beta_k|^2]\le \lambda^{-1}\mathcal E_{f}$ obtained from \eqref{eta}.
The combination of \eqref{splitt1}, \eqref{bound21} and \eqref{cheby} yields
\begin{equation}\label{step1_eq}
\begin{split}
\mathbb E_p[\sum_{\{n: |\bar\beta_n| > \epsilon\}} \frac{|\hat{ f_\epsilon}(\omega_n)|^2}{\bar p_n}]
&\le  (\sum_{n\in \mathbb Z^d}|\hat f(\omega_n)|)^2\big(1+(C'')^{1/2}\mathcal O(K^{-\frac{1-d/\ell}{6-2d/\ell}})\big)\\
&\quad +\sum_{n\in \mathbb Z^d}|\hat f(\omega_n)|(C'')^{1/2}\mathcal O(K^{-\frac{3(1-d/\ell)}{6-2d/\ell}})\\
&\quad + \sum_{n\in \mathbb Z^d}|\hat f(\omega_n)|C''\mathcal O(\lambda^{-1/2}K^{-\frac{1-d/\ell}{6-2d/\ell}})\,,
\end{split}
\end{equation}
which together with \eqref{namely} proves the lemma.
\end{proof}

\begin{lemma} [Estimate for $f-f_\epsilon$]\label{lemma3}
There holds
\[
\begin{split}
\mathbb E_{\check p}\Big[\varphi_2(\tilde\nu,\nu ) \Big] &=
\mathbb E_{\check p}\big[\min_{\hat v\in\mathbb C^{\tilde K}}\big(\mathbb E_x[\sum_{k=1}^{\tilde K}\hat v_ke^{{\rm i} \tilde\nu_k\cdot x}- (f-f_\epsilon)(x)|^2 +\lambda\sum_{k=1}^{\tilde K}|\hat v_k|^2]\big)\,\big]\\
&\le
C''\mathcal O(K^{-1})q_\epsilon^{-1} \epsilon^{1-d/\ell}\big(1+C''\mathcal O(K^{-(1-d/\ell)/(3-d/\ell})\big)^{(1-d/\ell)/2}\,.
\end{split}
\]
\end{lemma}
\begin{proof}
    The second term in \eqref{gen_split} has by \eqref{eta},\eqref{C''}, \eqref{theta} and Chebyshev's inequality the bound
\begin{equation}\label{remainder2}
\begin{split}
%&\mathbb E_{ p,\epsilon}\big[\min_{\hat\alpha'\in \mathbb C^K}\big(\mathbb E_x[|\alpha'(x)- (f-f_\epsilon)(x)|^2] 
%+\lambda\sum_{k=1}^K|\hat\alpha'_k|^2\big)\big]\\
\mathbb E_{\check p}\Big[\varphi_2(\tilde\nu,\nu ) \Big] &\le\big(Kq_\epsilon\big)^{-1}(1+\lambda) \mathbb E_p\big[\sum_{\{n: |\bar\beta_n| \le \epsilon\}} \frac{|\hat{ f}(\omega_n)|^2}{ p(\omega_n)}\big]\\
&\le \big(Kq_\epsilon\big)^{-1}(1+\lambda) \mathbb E_p\big[C''\sum_{\{n: |\hat f(\omega_n)| \le \epsilon+\|\eta\|\} } |\hat{ f}(\omega_n)|\big]\\
&\le \big(Kq_\epsilon\big)^{-1}(1+\lambda) \Big(\mathbb E_p[(\epsilon+\|\eta\|)^{1-d/\ell}]\frac{C''\sigma_d
(c/2)^{d-\ell} C^\ell}{\ell-d} 
+ \|\hat f\|_{\ell^1}\mathbb E_p[1_{\|\eta\|>(\frac{c}{2d^{1/2}})^\ell}]\Big)\\
&\le \big(Kq_\epsilon\big)^{-1}(1+\lambda) 
\Big((\mathbb E_p[(\epsilon+\|\eta\|)^2])^{(1-d/\ell)/2} \frac{C''\sigma_d(c/2)^{d-\ell} C^\ell}{\ell-d}\\
&\quad +\|\hat f\|_{\ell^1}\mathbb E_p[\|\eta\|]/(\frac{c}{2d^{1/2}})^{\ell}\Big)\\
&=C''\mathcal O(K^{-1})q_\epsilon^{-1}(\epsilon^2+\mathcal E_f)^{(1-d/\ell)/2}\\
&=C''\mathcal O(K^{-1})q_\epsilon^{-1} \epsilon^{1-d/\ell}\big(1+C''\mathcal O(K^{-(1-d/\ell)/(3-d/\ell})\big)^{(1-d/\ell)/2}\,.\\
 \end{split}
\end{equation}
\end{proof}

Now we proceed with the proof of Theorem \ref{theorem1}.
\begin{proof}    
Using  Lemmas \ref{lemma1}, \ref{lemma2} and \ref{lemma3}, we conclude by %\eqref{gen_split},\eqref{gen_split2}, \eqref{step1_eq}, \eqref{remainder2}, 
the assumption on  the  decay \eqref{decay} of the Fourier coefficients and the resolution property, \eqref{N_bound}, that  for $\delta=q_\epsilon^{1/2}$, $ \lambda=1$, $J=\infty$ and \[\mathcal E_f/\epsilon^2=C''\mathcal O(K^{-(1-d/\ell)/(3-d/\ell)})\ll 1,\]  the neural network based on resampling has the generalization error
\begin{equation}\label{conclude_4}
\begin{split}
\mathbb E_{\check p}\big[ %\min_{\hat{\bar \alpha}\in\mathbb C^K} 
\mathbb E_x[|\alpha(x;\bar\omega)-f(x)|^2]\big]
%\mathbb E_{\bar p}\big[\min_{\hat \alpha} \mathbb E_x[|\alpha(x)-f(x)|^2]\big]
&\le  K^{-1}(1+\lambda)(\sum_{n\in\mathbb Z^d}|\hat{f}(\omega_n)|)^2
(1+\delta)\big(1+\mathcal O(q_\epsilon+\frac{\mathcal E_f^{1/2}}{\epsilon})\big) \\
%\kappa_K(1+\gamma^{-1/4})(1+\gamma^{-1/2})\\
%\\&\quad 
 &\quad +(1+\delta^{-1})K^{-1}C''\mathcal O( \frac{\epsilon^{1-d/\ell}}{q_\epsilon})\\
 &= K^{-1}(1+\lambda)(\sum_{n\in\mathbb Z^d}|\hat{f}(\omega_n)|)^2 +K^{-1}C''\mathcal O(q_\epsilon^{1/2})
 \,,\\
\end{split}
\end{equation}
where the remainder in the last inequality is obtained by first estimating
\[
\begin{split}
&\delta +q_\epsilon + \epsilon^{-1}\mathcal E_f^{1/2}+ (1+\delta^{-1})C'' q_\epsilon^{-1} \epsilon^{1-d/\ell}\\
    &\le C'' \Big(\delta +\big(q_\epsilon + (1+\lambda)^{1/2}\epsilon^{-1} K^{-1/2}\big)+(1+\delta^{-1})q_\epsilon^{-1}\epsilon^{1-d/\ell}\Big)
\end{split}
\]
and then letting $q_\epsilon^2=\epsilon^{1-d/\ell}$ together with $K^{-1/2}\epsilon^{-1}=q_\epsilon=
\epsilon^{(1-d/\ell)/2}$ and $\delta=q_\epsilon/\delta$.
The leading order term $ K^{-1}(1+\lambda)(\sum_{n\in\mathbb Z^d}|\hat{f}(\omega_n)|)^2=K^{-1}(1+\lambda)C_{p_*}$ in \eqref{conclude_4} is the optimal  rate constant and dominates for large values $K\gg 1$.
\end{proof}
\section{Finite number of nodes \texorpdfstring{$K$}{K} and finite amount of data \texorpdfstring{$J$}{J}}
\label{sec:noise}
%\begin{remark}
In the case with finite amount of data, $\{(x_j,y_j)\, j=1,\ldots,J\}$ and $y_j=f(x_j)+\xi_j$ including independent noise $\xi_j,\ j=1,\ldots, J,$ 
with mean $\mathbb E[\xi_j]=0$ and  variance $\mathbb E[|\xi_j|^2]=s^2$,  we use in  Theorem \ref{thm_noise} for given regularization parameters $\lambda_1>0$ and $\lambda_2>0$
%the work \cite{HPSS} uses
the minimization problem for $\beta(x;\nu)=\sum_{k=1}^K\hat \beta_k e^{{\rm i}\nu_k\cdot x}$
\begin{equation}\label{min3}
    \min_{\hat\beta\in\mathbb C^K}\big(\frac{1}{J}\sum_{j=1}^J|\beta(x_j)- y_j|^2 +\lambda_1\sum_{k=1}^K|\hat\beta_k|^2
    +\lambda_2(\sum_{k=1}^K|\hat\beta_k|^2)^2\big)\,,
   %+\lambda_2\sum_{k=1}^K|\hat\beta_k|\big)
\end{equation}
with  independent samples $\nu_1,\ldots,\nu_K$ from $p$,
%derives the generalization error in the case without noise, i.e. $s=0$. We extend in Theorem \ref{thm_noise}
%the analysis to include noise with and 
and prove for any $\delta>0$ the generalization error bound
%to derive the generalization error estimate for any $\delta>0$, see Theorem ,
\begin{equation}\label{gen_Jxi2}
\begin{split}
&\mathbb E_ {\{x_j,\xi_j\}}\Big[\mathbb E_p \big[\mathbb E_x[ |f(x)-\beta(x)|^2\,\big|\, \{x_j,\xi_j,\nu_k\}] \big] \Big]\\
%&=\mathcal O\big(\frac{1}{K}+\lambda +(1+\frac{1}{\lambda K})(\sqrt{\frac{\log\frac{K}{\delta}}{J}})+\delta\big)\\
%&=\mathcal O\big(\frac{1}{K}+\frac{1+|s|}{\sqrt J} \big)\,,
&\le (1+3\delta/2)\big((\frac{1}{K}+\frac{\lambda_1}{K})C_p +\frac{\lambda_2}{K^2}(C^2_p+\frac{C_p'}{K})\big)
+\frac{Kc_1}{2\delta\lambda_1 J} +\frac{K^2}{2\delta\lambda_2J}\,, \\
%&\le\frac{(1+3\delta/2)C_p}{K}(1+\lambda) + (\frac{(1+3\delta/2) c_1C_p}{2J\delta})^{1/2}\\
%&\quad +\big((1+3\delta/2)(C_p^2+\frac{C_p'}{K})\big)^{1/2} \big(\frac{1}{K}+\frac{1}{(J\delta)^{1/2}}\big)\\
%%+(\frac{(1+3\delta/2)(C_p^2+\frac{(C'')^2C_p}{K})}{J\delta})^{1/2}\\
%&=:\mathcal E(C_p,C'')
\end{split}
\end{equation}
 where
 \[
 \begin{split}
 c_1 &:=\mathbb E_x[ |f(x)|^2]+s^2\,,\\
 C_p' &:= \sum_{n\in\mathbb Z^d} \frac{|\hat f(\omega_n)|^4}{p(\omega_n)^3}\le C_p(C'')^2\,,
 \end{split}
 \]
 and $\mathbb E_{\{x_j,\xi_j\}}[\ldots]$
denotes the expected value with respect to the data set. 
 The choice 
 \[
 \begin{split}
 \lambda_1&=\max\big(K(\frac{(1+3\delta/2)C_pc_1}{2J\delta})^{1/2},\lambda\big)\,,\\ 
 \lambda_2 &=K\max\Big(K\big(\frac{1}{(1+3\delta/2)(C_p^2+K^{-1}C'_p)2J\delta}\big)^{1/2},\lambda\big((1+3\delta/2)(C_p^2+K^{-1}C'_p)2\delta\big)^{-1/2}\Big)\,,
 \end{split}
 \] 
 where $\lambda<2$ is a positive constant,
 implies
 \begin{equation}\label{Ec}
 \begin{split}
& \mathbb E_ {\{x_j,\xi_j\}}\Big[\mathbb E_p \big[\mathbb E_x[ |f(x)-\beta(x)|^2
\,\big|\, \{x_j,\xi_j,\nu_k\}] \big] \Big]\\
 &\le\frac{(1+3\delta/2)C_p}{K}(1+\lambda) + (\frac{2(1+3\delta/2) c_1C_p}{J\delta})^{1/2}\\
&\quad +\big((1+3\delta/2)(C_p^2+c_K)\big)^{1/2} \big(\frac{\lambda}{K}+\frac{2^{1/2}}{(J\delta)^{1/2}}\big)\\
%+(\frac{(1+3\delta/2)(C_p^2+\frac{(C'')^2C_p}{K})}{J\delta})^{1/2}\\
&=:\mathcal E(C_p,c_K)
 \end{split}
 \end{equation}
 and $c_K:=K^{-1}C_p'$.
 The work \cite{HPSS} derives the analogous generalization error in the case without noise.

We have roughly 
\[
\mathcal E(C_p,c_K)\approx \max\big((1+2\lambda)\frac{C_p}{K}, \frac{3^{1/2}C_p+(3c_1C_p)^{1/2}}{J^{1/2}}\big)\,,
\]
by choosing $\delta$ small if $J^{1/2}\gg K$ and $\delta$ large if $J^{1/2}\ll K$.
We observe that the generalization error bound \eqref{gen_Jxi2} is minimized  by using the distribution $p_*$ that minimizes $C_p$.
We note that also with the additional regularization term $\lambda_2(\sum_{k=1}^K|\hat\beta_k|^2)^2$
the equal amplitude property \eqref{equal}  and \eqref{abs} hold. Therefore,  the proof of Theorem \ref{theorem1} can
be applied to the optimization \eqref{min3} with the generalization error bound \eqref{gen_Jxi2} replacing $K^{-1}(1+\lambda)C_p$.  
Consequently, the generalization error estimate \eqref{gen_Jxi2}
% with the additional asymptotically negligible error term $\mathcal O(J^{-1/2})$, 
can  replace \eqref{eta}  %\eqref{resamp} 
for finite amount of data $J$  and we obtain as in Theorem \ref{theorem1},
 by choosing $\epsilon := \mathcal E(C_p,c_K)^{\frac{1}{4}}$ and 
 $q_\epsilon:=\mathcal E(C_p,c_K)^{\frac{1-d/\ell}{8}}$, 
 %and replacing $C''/K$ in $\mathcal E(C_{p_*})$ by $\|\hat f(\omega_n)\|_{\ell^\infty}/(\epsilon K)$ for $\bar\alpha$ 
%we obtain as in Theorem \ref{theorem1}
that for any $K$ and $J$  the resampled generalization error for $\alpha$ in Definition \ref{assump1} has the asymptotically minimal bound 
\begin{equation}\label{conclude3}
\begin{split}
\mathbb E\big[
%\min_{\hat{\bar \alpha}\in\mathbb C^K}
\mathbb E_x[|\alpha(x)-f(x)|^2\,\big|\, \{x_j,\xi_j,\nu_k\}]\big]
&\le  \mathcal E\big(C_{p_*}, 0\big)\big(1+C''\mathcal O(q_\epsilon^{1/2})\big) 
%\kappa_K(1+\gamma^{-1/4})(1+\gamma^{-1/2})\\
%\\&\quad 
% +K^{-1}\mathcal O(q_\epsilon^{1/2}) +\mathcal O(J^{-1/2})
\,,\\
\end{split}
\end{equation}
as $\mathcal E(C_p,c_K)\to 0+$. Here  the expected value $\mathbb E[\ldots]$ is based on the data set and  the sampling  with respect to $\check p$ in \eqref{eq:barpdef}. 
%$p$ for the frequencies $\omega_k$ and the resampling
%with respect $\bar p_n$ for the frequencies  $\nu_\ell$.

\begin{proof}[Proof of \eqref{conclude3}]
The aim is to estimate  $\mathbb E[C_{\bar p}(f_\epsilon)]$ and $\mathbb E[C'_{\bar p}(f_\epsilon)/K]$
in the resampling setting and then replace $C_p$ with  $\mathbb E[C_{\bar p}(f_\epsilon)]$
and $C'_p/K$ by zero in \eqref{Ec}.
Following the proof of Theorem \ref{theorem1} in Lemma \ref{lemma2} we therefore estimate as in \eqref{step1_eq}
the rate function \[\mathbb E[C_{\bar p}(f_\epsilon)]=\mathbb E[\sum_{\{n: |\bar\beta_n| > \epsilon\}} \frac{|\hat{ f_\epsilon}(\omega_n)|^2}{\bar p_n}]\,,\]
replacing $C_p(f)$ in \eqref{Ec},
and in addition the new quantity
\[\mathbb E[\frac{C'_{\bar p}(f_\epsilon)}{K}]=K^{-1}\mathbb E[\sum_{\{n: |\bar\beta_n| > \epsilon\}} \frac{|\hat{ f_\epsilon}(\omega_n)|^4}{\bar p_n^3}],\]
replacing to $C'_p(f)/K$,
now choosing
 \[
 \begin{split}
 %&\mbox{$K^{-1}$ by $\mathcal E(C_p,c_K)$,}\\
 &\mbox{ $\epsilon=\mathcal E(C_p,c_K)^{1/4}$ , }\\
& \mbox{$q_\epsilon=\lambda^{-1/2}\epsilon^{(1-d/\ell)/2}= \lambda^{-1/2}\mathcal E(C_p,c_K)^{(1-d/\ell)/8}$}\,. 
%by $\lambda^{-1/2}\mathcal E(C_p,c_K)^{-(1-d/\ell)/(6-2d/\ell)}$.}
\end{split}
\]

{\it Estimation of $C_{\bar p}(f_\epsilon)$.} We have as in \eqref{step1_eq} the estimate 
\[\mathbb E[\sum_{\{n: |\bar\beta_n| > \epsilon\}} \frac{|\hat{ f_\epsilon}(\omega_n)|^2}{\bar p_n}]
\le (\sum_{n\in\mathbb Z^d}|\hat f(\omega_n)|)^2 + \mathcal O\big(\frac{\mathcal E(C_p,c_K)^{1/2}}{\epsilon}\big)\,,
\]
%for $\mathbb E_p[\sum_{\{n: |\bar\beta_n| > \epsilon\}} \frac{|\hat{ f_\epsilon}(\omega_n)|^2}{\bar p_n}]$
and use Jensen's inequality to obtain the bound
\begin{equation}\label{E1}
\begin{split}
\mathbb E[(\sum_{\{n: |\bar\beta_n| > \epsilon\}} \frac{|\hat{ f_\epsilon}(\omega_n)|^2}{\bar p_n})^{1/2}] &\le 
\big(\mathbb E[\sum_{\{n: |\bar\beta_n| > \epsilon\}} \frac{|\hat{ f_\epsilon}(\omega_n)|^2}{\bar p_n}]\big)^{1/2} \\
&
\le \sum_{n\in\mathbb Z^d}|\hat f(\omega_n)| + \mathcal O\big((\frac{\mathcal E(C_p,c_K)^{1/2}}{\epsilon})^{1/2}\big)\,.
%\mathcal O\big(\mathcal E(C_p,c_K)^{-(1-d/\ell)/(12-4d/\ell)}\big)\,.\\
%\,.
\end{split}
\end{equation}

{\it Estimation of $C'_{\bar p}(f_\epsilon)/K$.} The estimate of $C'_{\bar p}(f_\epsilon)/K$ uses 
\[
\begin{split}
\mathbb E[\sum_{\{n: |\bar\beta_n| > \epsilon\}} \frac{|\hat{ f_\epsilon}(\omega_n)|^4}{\bar p_n^3}]
& = \mathbb E[\sum_{\{n: |\bar\beta_n| > \epsilon\}} \frac{|\hat{ f_\epsilon}(\omega_n)|^4}{|\bar \beta_n|^3}
(\sum_{\{n: |\bar\beta_n| > \epsilon\}}|\bar\beta_n| )^3]\\
&\le\sum_{n\in\mathbb Z^d} \frac{|\hat{ f}(\omega_n)|^4}{\epsilon^3}
 \mathbb E[(\sum_{\{n: |\bar\beta_n| > \epsilon\}}|\bar\beta_n| )^3]
\end{split}
\]
together with Cauchy's and Jensen's inequalities for the bounds
\[
\begin{split}
 \mathbb E[(\sum_{\{n: |\bar\beta_n| > \epsilon\}}|\bar\beta_n| )^3] 
 &\le \mathbb E[(\sum_{k=1}^K|\hat\beta_k| )^3] \\
 &\le  \mathbb E[(K\sum_{k=1}^K|\hat\beta_k|^2 )^{3/2}]\\
 &\le  K^{3/2}\big(\mathbb E[(\sum_{k=1}^K|\hat\beta_k|^2)^2 ]\big)^{3/4}\\
 &\le \frac{K^{3/2}}{\lambda_2^{3/4}}\big(2\mathcal E(C_p, c_K)\big)^{3/4}\,,
\end{split}
\]
%obtained by Jensen's inequality 
where the last inequality is based on the training  error bound, proved in \eqref{train},
\[
\begin{split}
&\mathbb E_{\{x_j,\xi_j\}}\Big[\mathbb E_p\big[J^{-1}\sum_{j=1}^J|f(x_j)-\beta(x_j)|^2 +\lambda_1\sum_{k=1}^K|\hat\beta_k|^2
+\lambda_2(\sum_{k=1}^K|\hat\beta_k|^2)^2\big]\Big]
\le 2\mathcal E(C_p,c_K)\,.
%(1+\frac{3\delta}{2})\mathcal E_t + \frac{KC_3 s^2}{2\delta\lambda_1 J}\\
%&=\mathcal O(\frac{1}{K} + \frac{\lambda_1}{K} + \frac{Ks^2}{\lambda_1 J} + \frac{\lambda_2}{K^2})\,.
\end{split}
\]
%it is this step that uses 
By the construction of $\lambda_2$ we have $\lambda_2\gtrsim K\lambda>0$, so that
% \[
% \mathbb E_ {\{x_j,\xi_j\}}\Big[\mathbb E_\omega \big[\mathbb E_x[ |f(x)-\beta(x)|^2
%+\lambda_1\sum_{k=1}^K |\hat\beta_k|^2+\lambda_2(\sum_{k=1}^K |\hat\beta_k|^2)^2
%\,\big|\, \{x_j,\xi_j,\omega_k\}] \big] \Big]\le 2\mathcal E(C_p,c_K)\,.
%\]
the bound for the remainder $K^{-1}C'_{\bar p}$ in \eqref{gen_Jxi2} becomes
\begin{equation}\label{E2}
\begin{split}
K^{-1}\mathbb E[\sum_{\{n: |\bar\beta_n| > \epsilon\}} \frac{|\hat{ f_\epsilon}(\omega_n)|^4}{\bar p_n^3}] &\le
K^{-1}\epsilon^{-3}K^{3/2}\lambda_2^{-3/4}\big(2\mathcal E(C_p, c_K)\big)^{3/4} \\
&=2^{3/4}
K^{1/2}\lambda_2^{-3/4}(\frac{\mathcal E(C_p, c_K)^{1/4}}{\epsilon})^3\\
%&=2^{3/4}\frac{1}{K^{1/4}\lambda^{3/4}}\mathcal E(C_p, c_K)^{ -\frac{1+d/\ell}{4(1-d/(3\ell)}}\\
&=\mathcal O(K^{-1/4})\,.
\end{split}
\end{equation}
using the assumption $\epsilon=\mathcal E(C_p,c_K)^{1/4}$.

{\it The generalization error.} We obtain as in \eqref{conclude_4}, now using
\eqref{E1} and \eqref{E2}, the generalization error
\begin{equation*}
\begin{split}
\mathbb E\big[
%\min_{\hat{\bar \alpha}\in\mathbb C^K} 
\mathbb E_x[|\alpha(x)-f(x)|^2]\big]
&\le  \mathcal E(C_{p_*},0)
(1+\delta)\big(1+\mathcal O(q_\epsilon+\frac{\mathcal E(C_{p},c_K)^{1/2}}{\epsilon}+K^{-1/4})\big) \\
 &\quad +(1+\delta^{-1})\mathcal E(C_{p},c_K)\mathcal O( \frac{\epsilon^{1-d/\ell}}{q_\epsilon})\\
 &=  \mathcal E\big(C_{p_*}, 0\big)\big(1+C''\mathcal O(q_\epsilon^{1/2})\big) 
 \,,\\
\end{split}
\end{equation*}
%
%We obtain the remainder terms
%\[
%\mathcal O\big(K^{-1/4} 
%+ \frac{\mathcal E(C_p, c_K)^{1/2}}{\epsilon} + q_\epsilon +C''\frac{\epsilon^{1-d/\ell}}{q_\epsilon}\big)\,,
%\]
%which replace the remainder
%$\mathcal O(\frac{\mathcal E_f^{1/2}}{\epsilon} + q_\epsilon +C''\frac{\epsilon^{1-d/\ell}}{q_\epsilon})$
%in \eqref{conclude_4}
%and proves that we asymptotically can replace $C_p$ by $C_{p*}$ and $c_K$ by zero.
which proves \eqref{conclude3}  also using the assumption $\mathcal E(C_{p},c_K)^{1/2}/\epsilon\ll 1$.
%using also the estimates in Step 2 of Theorem \ref{theorem1}.

\end{proof}

\section{The resampling iterations}
\label{sec_iter}
This section motivates using   iterations for the frequencies 
$\nu[j]=(\nu_1,\ldots,\nu_K)[j]$, with $j=1,2,3,\ldots$, by combining resampling steps and  random walk steps.
We may do the random walk steps by projecting a normally distributed random variable in $\mathbb R^d$ to the lattice $\frac{\pi}{L}\mathbb Z^d$.
%The resampling will contribute to the drift and the random walk increments to the diffusion. 
Sections \ref{finite-K-section} and \ref{sec:noise} demonstrate that the resampling method can accurately approximate the optimal sampling distribution $p_*$ provided the relevant Fourier modes are present in the set of available frequencies and  condition \eqref{C''} holds. Now we study how to iteratively obtain such a resolving set of frequencies. In particular we will investigate how the variance of the random walk relates to the resampling error term $\mathcal O(q_\epsilon^{1/2})$  in the rate constant bound \eqref{conclude} and \eqref{conclude3}. 

We start, for instance, with $K$ independent  normally distributed samples $\nu_1,\ldots, \nu_K$,
with standard deviation $\delta$ each projected to the nearest lattice point. In the case $\delta=1$
these frequencies will by \eqref{C''} only resolve Fourier modes with low frequencies, namely the frequencies $\omega_n$ satisfying
\[
\frac{e^{-|\omega_n|^2/2}}{(2\pi)^{d/2}}\ge \frac{|\hat f(\omega_n)|}{C''}=\mathcal O(\epsilon)\,,
%\mathcal O(K^{-1/(3-d/\ell)})\,,
\]
which implies $|\omega_n|\le (\log \epsilon^{-1})^{1/2}+\mathcal O(1)$\,.
%$|\omega_n|\le (\frac{2}{3-d/\ell}\log K)^{1/2} +\mathcal O(1)$\,.
Theorem \ref{theorem1} and the splitting of $f$ in \eqref{gen_split} shows that the cut function 
\begin{equation}\label{f1}
f_{\epsilon_1}(x):=
\sum_{|\omega_n|\le (\log \epsilon^{-1})^{1/2}+c}\hat f(\omega_n)e^{{\rm i}\omega_n\cdot x}
%\sum_{|\omega_n|\le(\frac{2}{3-d/\ell}\log K)^{1/2} +c}\hat f(\omega_n)e^{{\rm i}\omega_n\cdot x}
\end{equation}
will be near optimally approximated by the resampling method using the initial standard normal distributed variables.
The cutoff with the function $f_{\epsilon_1}$ introduces a substantial approximation error of the order 
$\mathcal O\big(\log(\epsilon^{-1})\big)$\,.
%$\mathcal O((\log K)^{d/2-\ell})$.
Therefore the frequencies need to be updated. 
A straightforward  method is to do random walk steps based on the resampled frequencies.
%An alternative is to combine the random walk with an adaptive Metropolis step, where the Metropolis test makes a decision of the random walk step will be accepted or rejected. 
 Then iterative resampling and random walk steps are continued until the approximation of the function $f$ does not improve.

Assume the probability to be in state $\omega_m$ is $q(j,\omega_m)$, after the resampling in iteration $j$, then the  probability 
after the random walk  in $\frac{\pi}{L}\mathbb Z^d$ satisfies %the Kolmogorov equation 
\begin{equation}\label{rwd}
\tilde q(j+1,\omega_n)=\sum_{\{m\in\mathbb Z^d\, :\, |\bar\beta_m|\ge \epsilon_j\}} q(j,\omega_m) p(\omega_m,\omega_n)\,,
\end{equation}
where $p(\omega_m,\omega_n)$ is the transition probability for the random walk from $\omega_m$ to $\omega_n$
and $\epsilon_j$ is the cutoff in iteration $j$.
Before the projection to a lattice point  we have the normally distributed transition probability
\[
\tilde p(\omega_m,\omega)=\frac{e^{-|\omega_m-\omega|^2/(2t)}}{(2\pi t)^{d/2}},
\]
where $t>0$ is the chosen variance, which cannot be much smaller than $(\pi/L)^2$ due to the projection to lattice points. After the projection to a lattice point the transition probability changes slightly 
to
\[
p(\omega_m,\omega_n)= \int_{[-\frac{1}{2},\frac{1}{2}]^d}\tilde p(\omega_m,\omega_n+\frac{\pi}{L}\omega'){\rm d}\omega'\,.
\]
We seek the maximal distance between $\omega_m$, of a resampled frequency in iteration $j$, and a frequency $\omega_n$, obtained in the random walk, that can be resolved by an additional resampling step applied to  the frequencies we have after the random walk. The condition \eqref{C''} for accurate approximation using Theorem \ref{theorem1} including frequency $\omega_n$ becomes
\[
\frac{|\hat f(\omega_n)|}{\tilde q(j+1,\omega_n)}\le C''\,.
\]
We note that $\tilde q(j+1,\omega_n)$ need to be sufficiently large and we have 
\[
\tilde q(j+1,\omega_n)\ge q(j,\omega_m)p(\omega_m,\omega_n)
\]
where the frequency $\omega_m$  has $|\bar \beta_m|\ge \epsilon_j$. A sufficient condition for $\omega_n$ to be in the resolved set is therefore
\begin{equation}\label{p_q}
p(\omega_m,\omega_n)\ge \frac{|\hat f(\omega_n)|}{C''q(j,\omega_m)}=\mathcal O(q_{\epsilon_j}^{\tau})\,,
\end{equation}
for some $\tau>0$ obtained from the reduction $\mathcal O(q_{\epsilon_j}^{1/2})\ll q_{\epsilon_j}^{\tau}$ in the rate constant at each resampling, in Theorem \ref{theorem1} and \eqref{conclude3},  which yields
\[
e^{-|\omega_m-\omega_n|^2/(2t)}\ge \mathcal O(q_{\epsilon_j}^\tau) %\mathcal O(K^{-\delta})
\]
and 
\[
|\omega_m-\omega_n|=\mathcal O\big(\sqrt{\log(q_{\epsilon_j}^{-1})}\big)\,. %\mathcal O(\sqrt{\log K})\,.
\]
That is, we can use a factor $\mathcal O(q_\epsilon^\tau)$ of
the decay $\mathcal O(q_\epsilon^{1/2})$ of the error in rate constant bounds
\eqref{conclude} and \eqref{conclude3} to compensate for the increase of $C''$ in the random walk step.

Consequently at least a spreading of frequencies that differ $\mathcal O\big(\sqrt{\log(q_\epsilon^{-1})}\big)$
%$\mathcal O(\sqrt{\log K})$ 
is possible in the random walk step
and Theorem \ref{theorem1} and \eqref{conclude3} is applicable to recursively resample new frequencies obtained from random walk and cutoff, due to the gain $\mathcal O(q_\epsilon^{1/2})$ in the rate constant from Theorem \ref{theorem1} and \eqref{conclude3}. We conclude that 
the variance of the random walk is related to the error bound of the rate constant.

In the next section we analyze quantitatively the evolution  of the rate constant in successive resampling steps. 

\subsection{Evolution of the rate constant}\label{sec_K}
Theorem \ref{theorem1} and \eqref{conclude3} show that the rate constant for the resampled method becomes
\begin{equation}\label{C''_Eq}
% C_{\bar p}=\big(1+C''\mathcal O(q_\epsilon^{1/2})\big) \mathcal E\big(C_{p_*}, 0\big)\,, %C_{p_*}\,,
C_{\bar p}=\big(1+C''\mathcal O(q_\epsilon^{1/2})\big)C_{p_*}\,,
\end{equation}
%which is accurate when
%$C''\ll q_\epsilon^{-1/2}$. 
provided  $\mathcal E(C_{p},c_K)^{1/2}/\epsilon=\mathcal O(1)$.
%holds initially}
%$C''K^{-(1-d/\ell)/(3-d/\ell)}=\mathcal O(1)$.
If we for instance initially have $p$ standard normal distributed $C''$ is typically unbounded. 
Assume for simplicity that $\hat f$ has compact support on $N_c^d$ frequencies, so that $C''$ is finite. For instance if
the distribution $p$, that samples all frequencies, is uniform we have $C''=\mathcal O(N_c^d)$.
In a setting with infinite number of nodes, $K=\infty$, sufficiently many frequencies exist, so that we can improve the rate constant, $C_{\bar p}$, by iterating the resampling and measure the estimate of the corresponding bound on $C''_j$  by \eqref{C''est} as
\[
C''_{j+1} = 1 + C''_j\underbrace{\mathcal O(q_{\epsilon_j}^{1/2})}_{=:\kappa_j}\,,
\]
using successive splittings $f_{\epsilon_j}$ of $f$
following \eqref{gen_split} without random walk steps.
 We obtain for $j=1,2,3,\ldots$ 
\begin{equation}\label{C''_rate}
C''_{j+1} = 1 +\sum_{i=2}^j\prod_{\ell=i}^j\kappa_\ell +C''_1\prod_{\ell=1}^j\kappa_\ell\,,
%\le \frac{1}{1-O(q_\epsilon^{1/2})} +C''_1\big(O(q_\epsilon^{1/2})\big)^j\,,
\end{equation}
which shows that 
\[
\limsup_{j\to\infty}C''_j=1+\mathcal O(\limsup_{j\to\infty}q_{\epsilon_j}^{1/2})\,,
\]
provided $\kappa_{j-1}\le \kappa_j/c'$ for some constant $c'$. We conclude that the rate constants for the resampling iterates asymptotically have the bound $\big(1+\mathcal O(J^{-(1-d/\ell)/16})\big)C_{p_*}$, with the convergence rate \eqref{C''_rate}
if initially $\mathcal E(C_{p},c_K)^{1/2}/\epsilon_1 =\mathcal O(1)$.

The case of finite number of nodes and several resampling steps requires creation of new frequencies 
%in order  
%to gradually increase the domain of frequencies to sample from
%and to remain 
to have an accurate set of frequencies. % and not depend only on a few starting frequencies.
%In addition there is a complication that frequencies can disappear in the resampling step.
%Neglecting that frequencies can disappear,  
Then the factor $\mathcal O(q_{\epsilon_j}^\tau)$ in \eqref{p_q}, of the decay $\mathcal O(q_{\epsilon_j}^{1/2})$ for the rate constant estimate \eqref{C''_Eq}, is required to compensate for the increase of the $C''_j$ factor in the random walk step. This reasoning motivates
that the random walk steps increase the rate constant bound to
\begin{equation}\label{C''tau}
\limsup_{j\to\infty}C''_j\le 1+\mathcal O(\limsup_{j\to\infty}q_{\epsilon_j}^{1/2-\tau})\,,
\end{equation}
which indicates that, for a suitable random walk parameter $\delta$ 
%using the Fourier transform version Algorithm \ref{alg_FT}
and large amount of data $J$ and nodes $K$, the iterations can converge. In the next section, numerical results show that the algorithm can be unstable for large values of $\delta$. 

\section{Numerical results and adaptive random walk}
\label{sec_num}
This section presents the numerical results based on the resampling/random walk Algorithms~\ref{alg_FT}, \ref{alg1}, and \ref{alg4}, where the least squares problem \eqref{min3} is solved iteratively  using the conjugate gradient method.
%For each resampling the samples are also ordered following the starting order although there is no significant improvement by keeping this order.
The main new contribution compared to the numerical results in \cite{Aku2}
is that here the least squares problem is solved iteratively and adaptive random walks are included, namely %Algorithms \ref{alg1} and \ref{alg2} have initial frequencies equal to zero. %the overparameterized setting is included
Algorithm \ref{alg4} employs adaptive random walk increments based on the estimated covariance of the frequencies
\begin{equation}\label{cov1}
\sum_{n=1}^{\hat N}\sum_{k=1}^K(\omega_k^n-m^n)(\omega_k^n-m^n)^*/(\hat NK)+\hat\epsilon{\rm I}
\end{equation}
%is tested in random walk steps 
to obtain faster convergence of the random walk and resampling iterations. % for $\hat\epsilon>0$.
%and the $d\times d$ identity matrix $\rm I$. 
Here $\omega_k^n\in\rset^d$  is the frequency at iteration $n$ with the empirical mean $m^n:=\sum_{k=1}^K\omega_k^n/K$, and  $\rm I$ denotes the $d\times d$ identity matrix.
%while in the random walk iteration $\hat N$ 
%\[\zeta_k \gets \mbox{random standard normal  in $\rset^{d}$}
%\]  is
%replaced by 
%\[
%\zeta_k\gets \mbox{random normal with mean zero and covariance \eqref{cov1}}\,.
%\]
The positive parameter $\hat\epsilon$ is introduced to ensure positive definite covariance approximations and is fixed as $\hat\epsilon=1/1000$ in the experiments. Numerical tests performed with values of  $\hat{\epsilon}$, ranging from $1/1000$ to $1/10$, exhibit similar performance, indicating robustness with respect to this parameter. The MATLAB codes for numerical implementations of the algorithms described in this section are hosted at the GitHub repository \cite{Code_repository}.
\begin{algorithm}[!ht]
\caption{Adaptive resampling with random walk on lattice}
\label{alg1}
\footnotesize{
\begin{algorithmic}
\State {\bfseries Input:} $\{(x_j, y_j)\}_{j=1}^J$$~\{\textrm{data}\}$
\State {\bfseries Output:} $x\mapsto\sum_{k=1}^K\hat\beta_k e^{{\rm i}\omega_k\cdot x}$
\State Choose a number of resampling iterations $N$, random walk step size $\delta$,  cutoff $\epsilon$,  and Tikhonov parameters $\lambda_1,~\textrm{and}~ \lambda_2$
\State Set $(\omega_1,\ldots,\omega_K) = 0$
%\gets \mbox{standard normal in $\rset^{Kd}$}$ 
%\Comment{zero initial frequencies}
\For{$n = 1$ {\bfseries to} $N$}
\State $\zeta_k \gets \mbox{sample standard normal  in $\rset^{d}$ for $k=1,\ldots,K$}$
\State $\omega_k \gets 
\mbox{project $\omega_k + \delta \zeta_k$ to the periodic lattice 
$\pi L^{-1}\mathbb Z^d$ for $k=1,\ldots,K$}$ $\{\textrm{random walk}\}$
\State $\hat{\beta} \gets \mbox{conjugate gradient approximation to \eqref{min3} given 
$\omega$ and the data}$ $\{\textrm{least squares}\}$
\State $\omega \gets \mbox{sample from $\{\omega_1,\ldots, \omega_K\}$  with replacement $K$  independent frequencies}$ based on Definition \ref{assump1} using multinomial resampling $~\{\textrm{resampling with cutoff}~ \epsilon\}$
\EndFor
\State $\hat{\beta} \gets \mbox{conjugate gradient approximation to \eqref{min3} given 
$\omega$ and the data}$
\State $x\mapsto\sum_{k=1}^K\hat\beta_k e^{{\rm i}\omega_k\cdot x}$
\end{algorithmic}}
\end{algorithm}
\begin{algorithm}[!ht]
\caption{Adaptive resampling with adaptive random walk in $\rset^d$}
\label{alg4}
\footnotesize{
\begin{algorithmic}
\State {\bfseries Input:} $\{(x_j, y_j)\}_{j=1}^J$$\{\textrm{data}\}$
\State {\bfseries Output:} $x\mapsto\sum_{k=1}^K\hat\beta_k e^{{\rm i}\omega_k\cdot x}$
\State Choose a number of resampling iterations $N$, random walk step size $\delta$, cutoff $\epsilon$, Tikhonov parameters $\lambda_1~\text{and}~\lambda_2$ and covariance regularization parameter $\hat\epsilon$
\State Set $(\omega_1,\ldots,\omega_K) =0$
%\gets \mbox{standard normal in $\rset^{Kd}$}$ 
%\Comment{zero initial frequencies}
\State Set $C_1={\rm I}$
$~\{\textrm{initial covariance equal to identity matrix}\}$
\For{$n = 1$ {\bfseries to} $N$}
\State $\zeta_k \gets \mbox{sample multivariate normal in $\rset^{d} $ with mean zero and covariance $C_n+\hat\epsilon {\rm I}$, for $k=1,\ldots,K$}$
\State $\omega_k \gets 
\mbox{ $\omega_k + \delta \zeta_k$ 
 for $k=1,\ldots,K$}$ $\{\textrm{random walk}\}$
\State $\hat{\beta} \gets \mbox{conjugate gradient approximation to \eqref{min3} given 
$\omega$ and the data}$ $\{\textrm{least squares}\}$
\State $\omega \gets$ sample with replacement $K$ times independently from $\{\omega_1,\ldots, \omega_K\}$   with probability $|\hat\beta_k|/\sum_{\{\ell:|\hat\beta_\ell|\ge \epsilon\}}|\hat\beta_\ell|$ for $\omega_k$ using multinomial resampling 
$\{\textrm{resampling with simplified cutoff}\}$
\State $m_n\gets$ average of all resampled frequencies $\omega_k$
\State $\hat C_n=\sum_{k=1}^K(\omega_k-m_n)(\omega_k-m_n)^*/K$ $~\{\textrm{empirical covariance matrix}\}$
\State $C_n\gets$ average of all $\hat C_m$ for $m=1,\ldots,n$ 
\EndFor
\State $\hat{\beta} \gets \mbox{conjugate gradient approximation to \eqref{min3} given 
$\omega$ and the data}$
\State $x\mapsto\sum_{k=1}^K\hat\beta_k e^{{\rm i}\omega_k\cdot x}$
\end{algorithmic}}
\end{algorithm}

\subsection{Numerical implementation: regression task}
\begin{figure}[ht]
\centering
\includegraphics[width=0.4\textwidth]{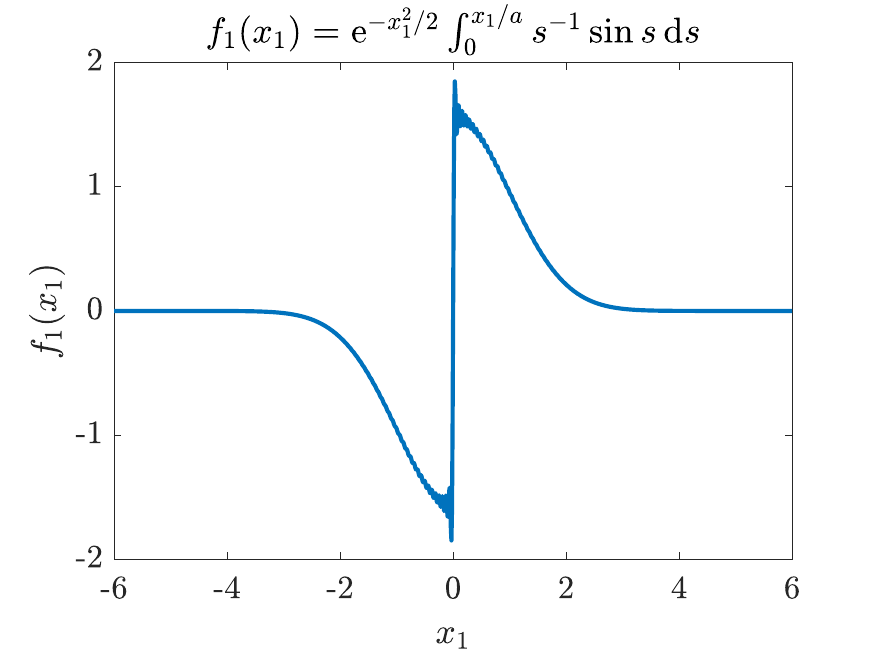}
\includegraphics[width=0.4\textwidth]{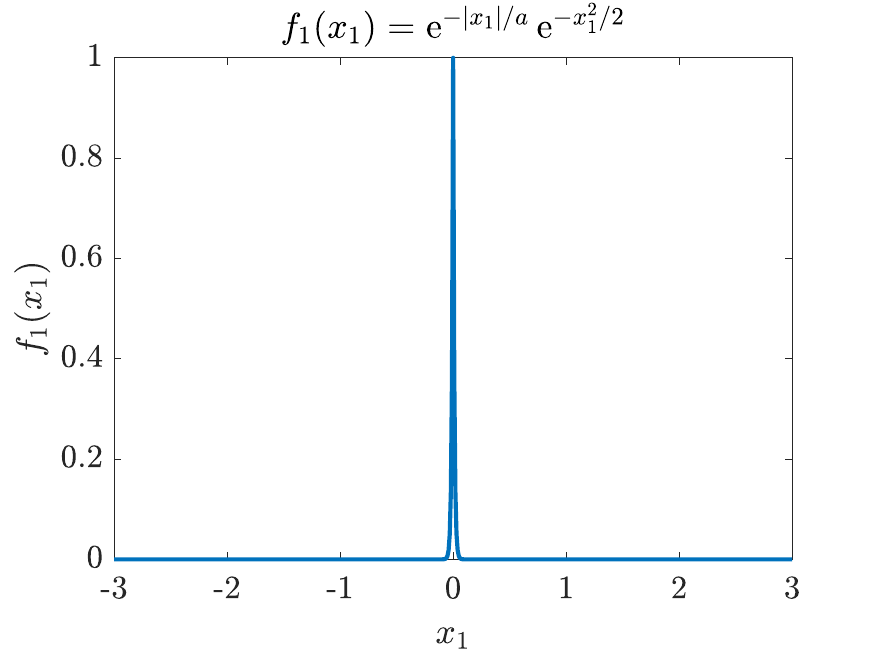}
\caption{ 
For $a=10^{-2}$ the functions $ f_1(x_1):= e^{-x_1^2/2}\int_0^{x_1/a} s^{-1}\sin(s)\mathrm{d}s$ (left) and $f_1(x_1):=e^{-|x_1|/a}e^{-x_1^2/2}$ (right)  are used for the anisotropic part of the functions $f$ in \eqref{eq:f} and \eqref{f_ex_2}.
%with standard normal distributed data points $(x_j,y_j)$ marked red. 
%The inset shows  $|\hat f|$.
%of its Fourier transform
%and the detail of its behaviour at the origin. 
}\label{fig:f_zoomed} 
\end{figure}

% Figures \ref{fig:f11} -- \ref{fig:f23} 
The numerical experiments of this section use independent standard normally distributed data points $x_j\in\rset^d, j=1,\ldots, J,$ and apply Algorithms \ref{alg_FT},~\ref{alg1}, and~\ref{alg4} to the $q$-periodic versions of the functions 
\begin{equation}\label{eq:f}
    f(x) = \int_0^{v\cdot x/a} \frac{\sin(s)}{s}\mathrm{d}s
    \, e^{-|x|^2/2}, \;\; x\in\rset^d,
\end{equation}
and
\begin{equation}\label{f_ex_2}
f(x)=e^{-|v\cdot x|/a}e^{-|x|^2/2}, \quad x\in\rset^d\,,
\end{equation}
for a unit vector
$v=(1,0,0,\ldots)$ in $\rset^d$ 
%(Figures \ref{fig:f11}-\ref{fig:f21}) 
or a random unit vector in $\rset^d$, 
%(Figures \ref{fig:f22}-\ref{fig:f23}), 
%with period $q=12$, 
 see Figure~\ref{fig:f_zoomed}. 
Both functions are anisotropic with slowly decaying Fourier coefficients.
The  Fourier transform of $x_1\mapsto\int_0^{x_1/a} s^{-1}\sin(s)\mathrm{d}s$ is supported in $|\omega_1|\le a^{-1}$ and decays as $|\omega_1|^{-1}$ while the Fourier transform of the function in \eqref{f_ex_2} yields the fat-tailed Cauchy distribution in one direction. This slow decay for both target functions \eqref{eq:f} and \eqref{f_ex_2} implies that sampling the frequencies could be challenging. 

Algorithms \ref{alg1} and \ref{alg4} have five input parameters: the random walk size $\delta$, the number of random walk/resampling iterations $N$, the cutoff threshold $\epsilon$, and the regularization weights $\lambda_1$ and $\lambda_2$. Further hyperparameters for the training of the random Fourier feature model $\beta(x)$ include the number of frequencies $K$,  the training  dataset size $J$,  and the noise standard deviation $s$ in the training  dataset.  This work investigates the effect of varying these hyperparameters in the numerical tests~\ref{tc:test1} to ~\ref{tc:test8}, with the specific settings of the hyperparameters summarized in Table~\ref{Table:numerical_tests}.
\begin{table}[!ht]
\centering
\begin{adjustbox}{max width=\textwidth, center}
\begin{tabular}{ccccccccc}
\toprule
Test case & Applied algorithms & $\delta$ & $K$ & $J$ & $\epsilon$ & $\lambda_1$ & $\lambda_2$ & $s$ \\
\midrule
Test~\ref{tc:test1} & Alg.~\ref{alg1} and~\ref{alg4}  & \textsc{Varied} & $2500$ & $8000$ & $\frac{1}{200}K^{-\frac{1}{2}}$ & $\frac{1}{100}KJ^{-\frac{1}{2}}$ & $0$ & $0$  \\
\midrule
Test~\ref{tc:test2} & Alg.~\ref{alg_FT}, ~\ref{alg1} and~\ref{alg4}  & $0.5$ & \textsc{Varied} & $20000$ & $\frac{1}{200}K^{-\frac{1}{2}}$ & $\frac{1}{20}KJ^{-\frac{1}{2}}$ & $0$ & $0$  \\
\midrule
Test~\ref{tc:test3} & Alg.~\ref{alg_FT} and~\ref{alg4}  & $0.2$ & $5000$ & $20000$ & $\frac{1}{200}K^{-\frac{1}{2}}$ & $\frac{1}{20}KJ^{-\frac{1}{2}}$ & $0$ & $0$ \\
\midrule
Test~\ref{tc:test4} & Alg.~\ref{alg1} & $0.2$ & $2500$ & \textsc{Varied} & $\frac{1}{200}K^{-\frac{1}{2}}$ & $\frac{1}{20}KJ^{-\frac{1}{2}}$ & $0$ & $0$ \\
\midrule
Test~\ref{tc:test5} & Alg.~\ref{alg1} & $0.5$ & $2500$ & $8000$ & \textsc{Varied} & $\frac{1}{20}KJ^{-\frac{1}{2}}$ & $0$ & $0$ \\
\midrule
Test~\ref{tc:test6} & Alg.~\ref{alg4} & $0.5$ & $2500$ & $8000$ & $\frac{1}{200}K^{-\frac{1}{2}}$ & \textsc{Varied} & $0$ & $0$ \\
\midrule
Test~\ref{tc:test7} & Alg.~\ref{alg4} & $0.5$ & $1250$ & $4000$ & $\frac{1}{200}K^{-\frac{1}{2}}$ & $\frac{1}{20}KJ^{-\frac{1}{2}}$ & \textsc{Varied} & $0$ \\
\midrule
Test~\ref{tc:test8} & Alg.~\ref{alg1} & $0.5$ & $10000$ & $50000$ & $\frac{1}{200}K^{-\frac{1}{2}}$ & $\frac{1}{100}KJ^{-\frac{1}{2}}$ & $0$ & \textsc{Varied} \\
\bottomrule
\end{tabular}
\end{adjustbox}
\vspace{1mm}
\caption{Hyperparameter settings for each test case. ``\textsc{Varied}'' indicates the hyperparameter is changed in the corresponding test case.}
\label{Table:numerical_tests}
\end{table}

\begin{figure}[!htbp]
  \centering
  \subfloat[Relative least squares error, Algorithm~\ref{alg1}]{%
    \includegraphics[width=0.49\textwidth]{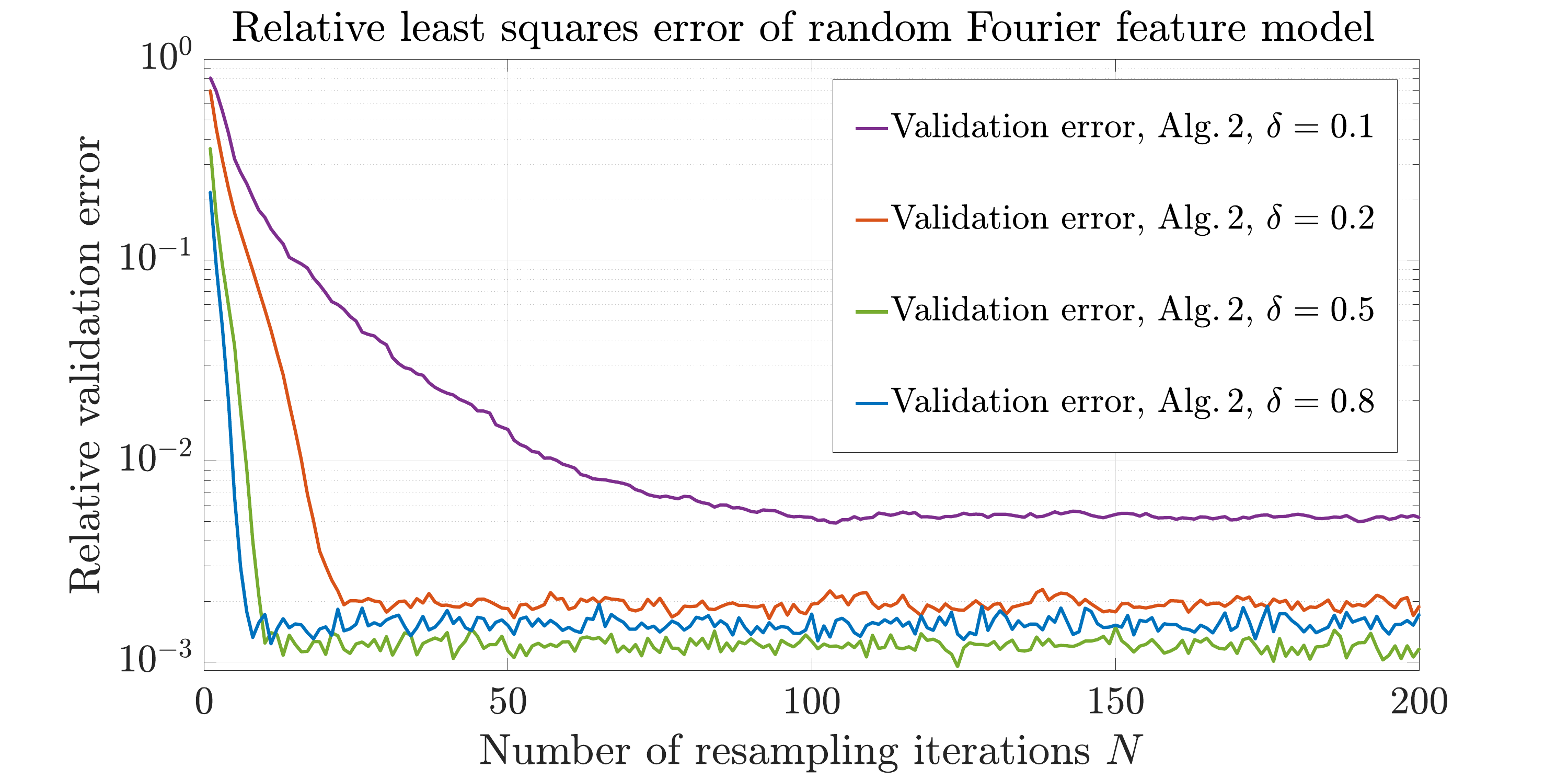}%
  }\hfill
  \subfloat[Relative least squares error, Algorithm~\ref{alg4}]{%
    \includegraphics[width=0.49\textwidth]{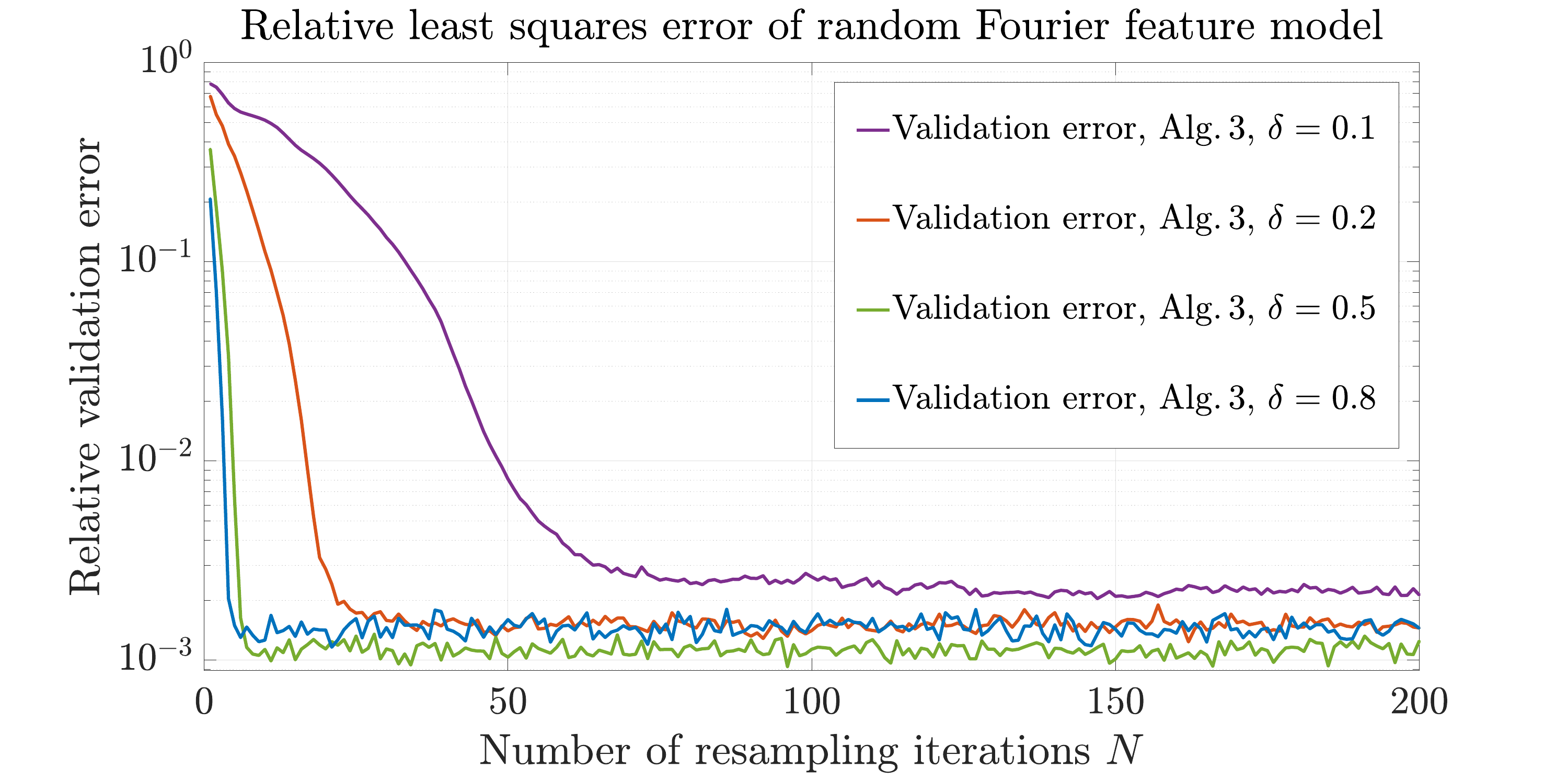}%
  }
  \caption{Relative least squares errors of the trained random Fourier feature model $\beta(x)$ using Algorithms~\ref{alg1} and~\ref{alg4}, with varied random walk step size parameter $\delta$.}
  \label{fig:Fig_a1_vary_delta}
\end{figure}

\Testcase{\label{tc:test1}
\textit{The effect of the random walk step size parameter $\delta$.}
}
As motivated in Section~\ref{sec_iter}, the random walk step size $\delta$  is critical for the performance of the resampling algorithm. If $\delta$ is too large, the resampled distribution becomes smeared out, less accurate, and eventually unstable. Conversely, if $\delta$ is too small, too few frequencies are activated, leading to inefficient frequency sampling and a slow convergence of the resampling iterations, as illustrated in Figure~\ref{fig:Fig_a1_vary_delta}.

We observe that a moderately increased random walk step size $\delta$ improves the convergence of resampling iterations. For this numerical test we use the target function $f(x)$
%$f(x) = \int_0^{v\cdot x/a}\frac{\sin{(s)}}{s}\mathrm{d}s\,\exp({-|x|^2/2)}$ 
given in \eqref{eq:f}, by applying Algorithms~\ref{alg1} and~\ref{alg4} with a period of $q=12$ in dimension $d=4$, with the following hyperparameters: $J=8000$, $K=2500$, $\epsilon=\frac{1}{200}K^{-\frac{1}{2}}$, $\lambda_1=\frac{1}{100}KJ^{-\frac{1}{2}}$, $\lambda_2=0$, and $a=0.1$.

\begin{figure}[!tbp]
    \centering
    \includegraphics[width=0.9\linewidth]{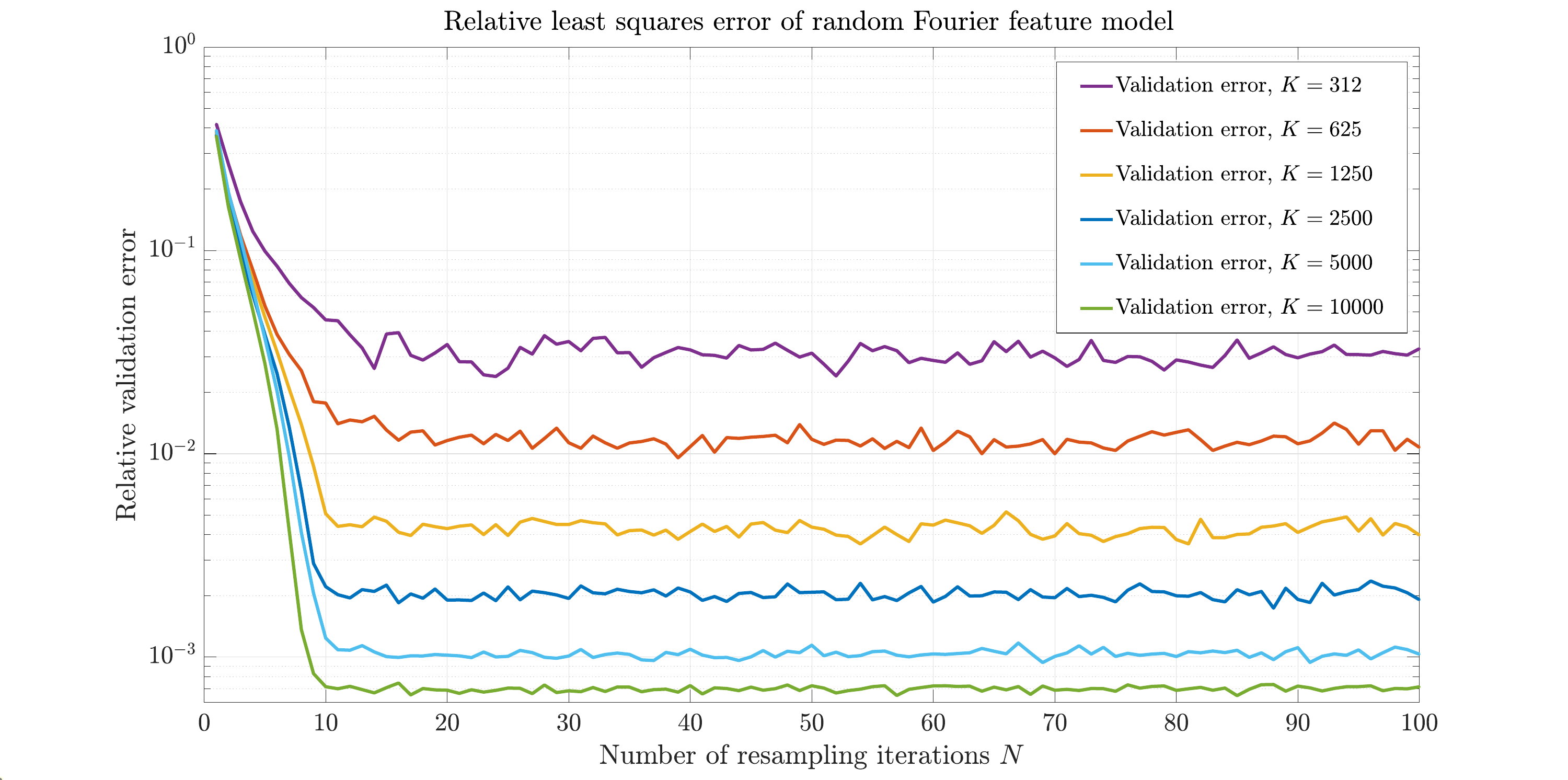}
    \caption{Relative least squares error of the trained random Fourier feature model $\beta(x)$ using Algorithm~\ref{alg1}, with gradually increased number of frequency parameters $K$.}
    \label{fig:Fig_a2_vary_K_1}
\end{figure}

% \subsubsection{The number of frequencies $K$}
\Testcase{\label{tc:test2}
\textit{The effect of the number of nodes $K$.}
}
To  validate the generalization error bound of the random Fourier feature model $\beta(x)$ with respect to the number of frequencies $K$ numerically, we systematically increase $K$ from 312 to 10,000. As illustrated in Figure~\ref{fig:Fig_a2_vary_K_1}, this increase  consistently improves   model accuracy.

  As $K$ increases, the relative least-squares prediction error of the trained random feature model decreases steadily from $3\times 10^{-2}$ to $6\times 10^{-4}$. In Figure~\ref{fig:Fig_a2_vary_K_1} we apply Algorithm~\ref{alg1} for the target function $f(x)$ given in \eqref{eq:f}, with the following hyperparameters: $J=20000$, $\delta=0.5$, $\epsilon=\frac{1}{200}K^{-\frac{1}{2}}$, $\lambda_1=\frac{1}{20}KJ^{-\frac{1}{2}}$, $\lambda_2=0$, $d=4$, $a=0.1$, and period $q=12$.

\begin{figure}[!htbp]
    \centering
    \includegraphics[width=0.9\linewidth]{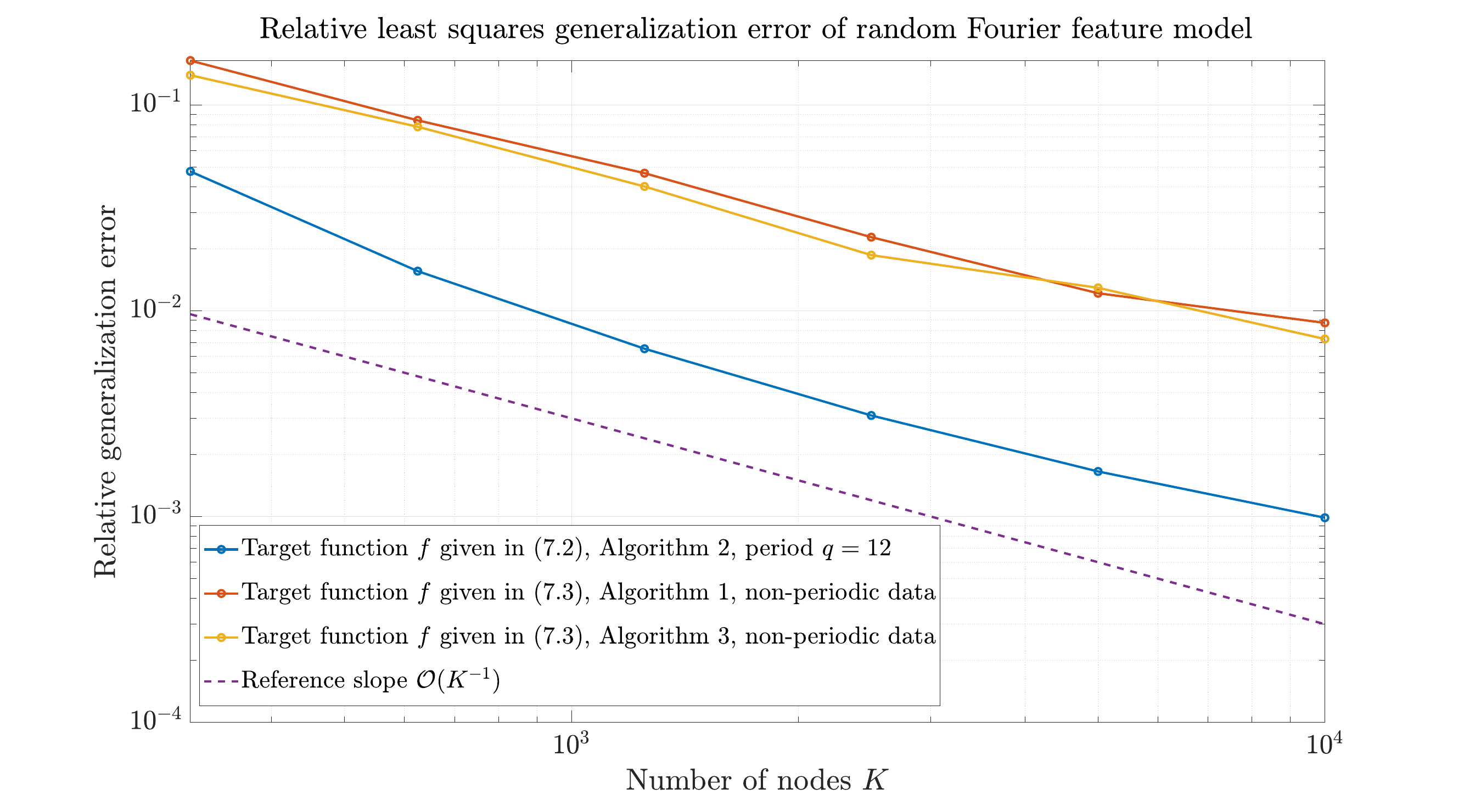}
    \caption{Relative generalization error of the trained random Fourier feature model $\beta(x)$, with an increased number of frequency parameters $K$, applying Algorithms~\ref{alg_FT},~\ref{alg1}, and~\ref{alg4}.}
    \label{fig:Fig_a2_vary_K_2}
\end{figure}

The scalability of the random Fourier feature model is also visually examined in Figure~\ref{fig:Fig_a2_vary_K_2}, where the relative least squares error on the test set shows agreement with the theoretical error bound $\mathcal{O}(K^{-1})$, in all  three test cases using Algorithms~\ref{alg_FT},~\ref{alg1}, and~\ref{alg4}.

\begin{figure}[!htbp]
    \centering
    \includegraphics[width=0.9\linewidth]{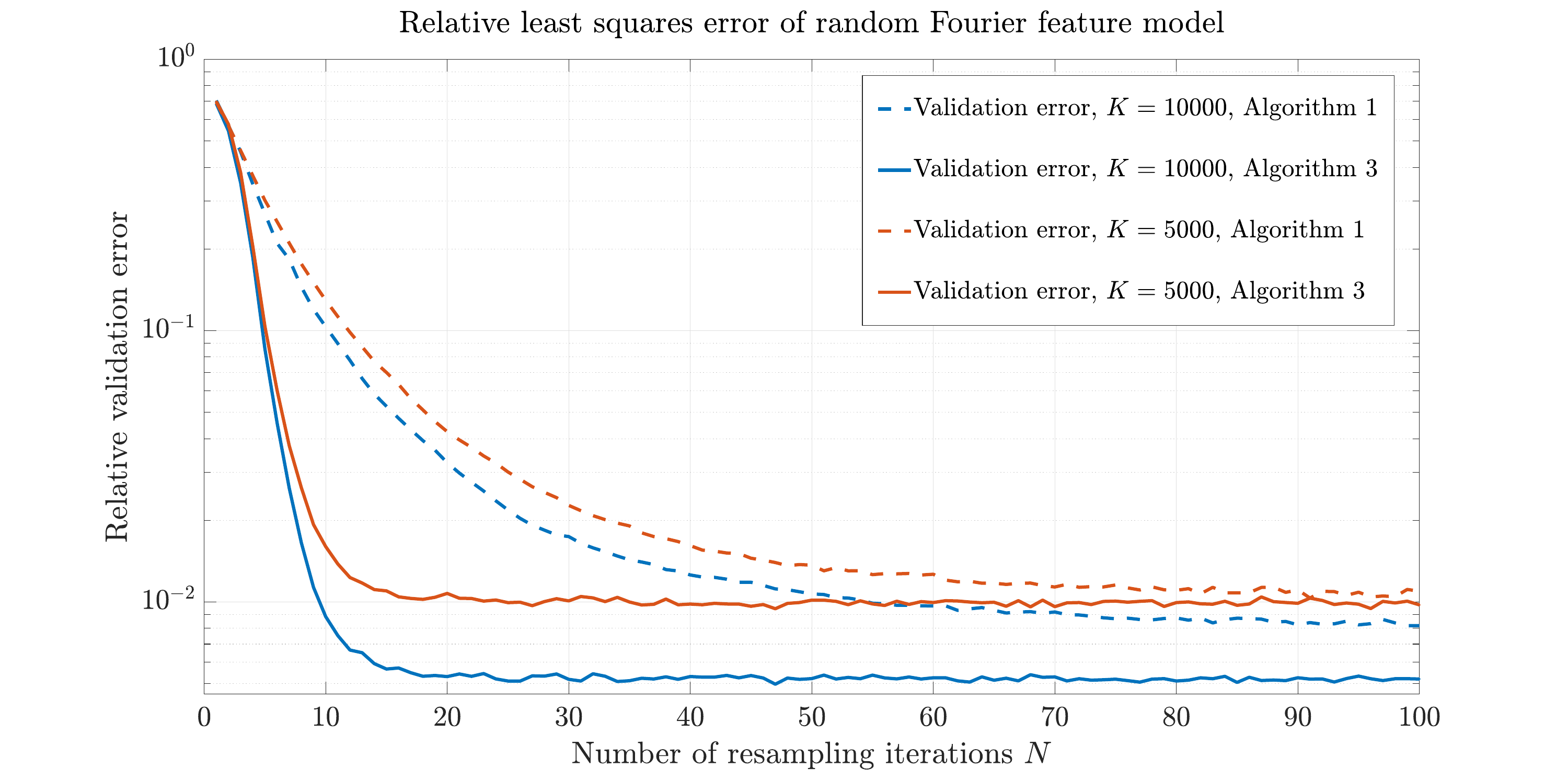}
    \caption{Relative least squares error of the trained random Fourier feature model $\beta(x)$ by applying Algorithm~\ref{alg_FT} (dashed lines) and Algorithm~\ref{alg4} (solid lines).}
    \label{fig:Fig_a2_vary_K_3}
\end{figure}

% \subsubsection{Improvement by adaptive random walk for resampling with Algorithm~\ref{alg1} }
\Testcase{\label{tc:test3}
\textit{Improvement by adaptive random walk for resampling with Algorithm~\ref{alg4}.}
}
 Figure~\ref{fig:Fig_a2_vary_K_3} compares the performance of Algorithms~\ref{alg_FT} and ~\ref{alg4} for training the random Fourier feature model with resampling iterations,  displaying enhanced convergence of Algorithm~\ref{alg4} by employing adaptive random walk using the estimated covariance of the frequency samples.
The random walk increments based on multivariate normal
distribution with estimated covariance \eqref{cov1} can improve the convergence of the random walk/resampling iterations, provided that the  number of nodes $K$ and the training dataset size $J$  are sufficiently large. Figure~\ref{fig:Fig_a2_vary_K_3}  presents an example  in which
Algorithm~\ref{alg_FT} requires a factor of 10 more iterations compared to Algorithm~\ref{alg4}.

The numerical tests applying Algorithms~\ref{alg_FT} and \ref{alg4} in Figures~\ref{fig:Fig_a2_vary_K_2} and~\ref{fig:Fig_a2_vary_K_3} use the target function $f(x)$ given in \eqref{f_ex_2}, with the the following hyperparameters: $J=20000$, $\delta=0.2$, $\epsilon=\frac{1}{200}K^{-\frac{1}{2}}$, $\lambda_1=\frac{1}{20}KJ^{-\frac{1}{2}}$, $\lambda_2=0$, $d=4$, and $a=0.1$.

\begin{figure}[htbp]
  \centering
  \subfloat[Noiseless training data, $s=0$]{%
    \includegraphics[width=0.49\textwidth]{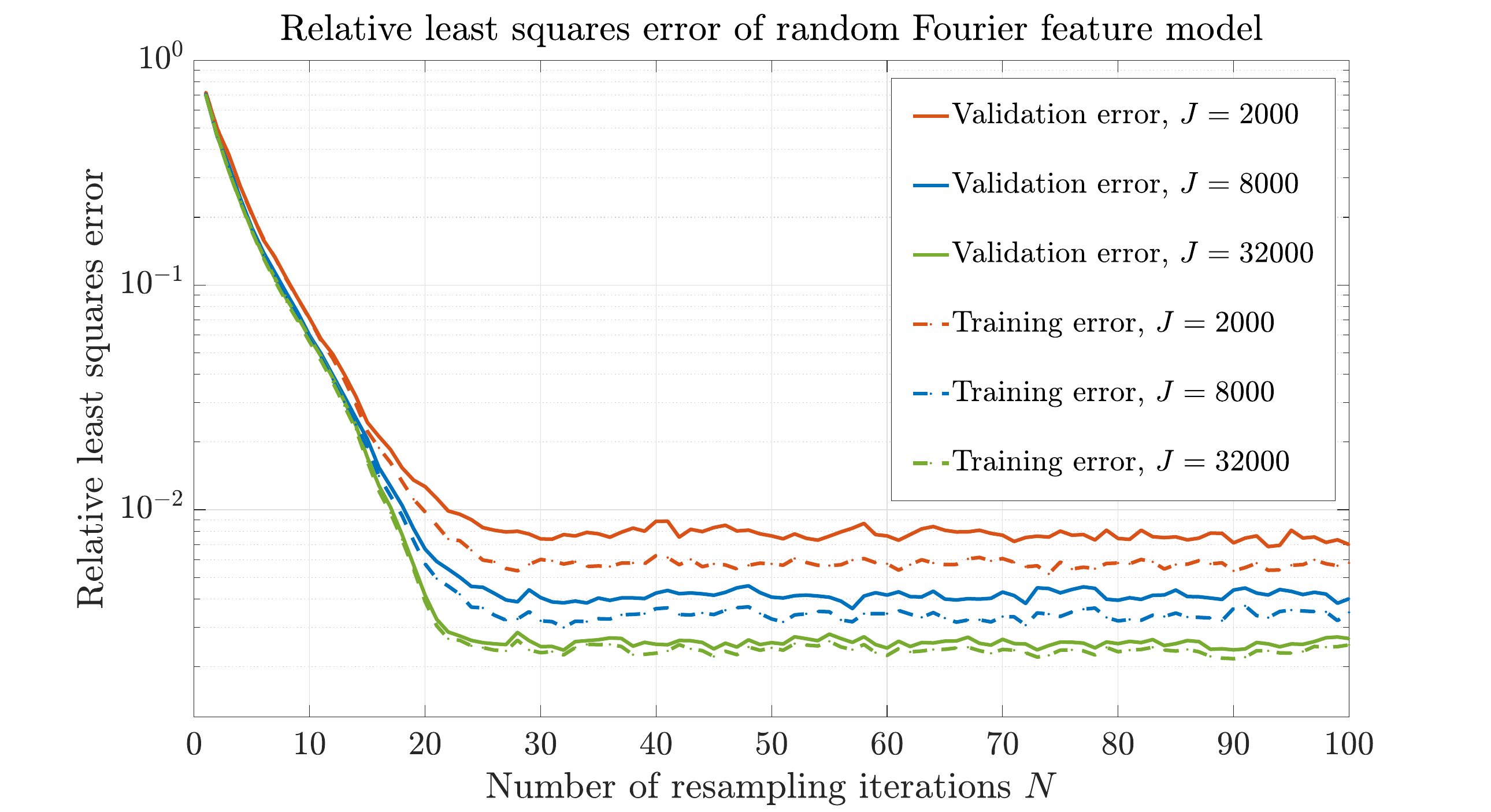}%
    \label{fig:vary_J_noiseless}%
  }\hfill
  \subfloat[Noisy training data, $s=2.5\times 10^{-3}$]{%
    \includegraphics[width=0.49\textwidth]{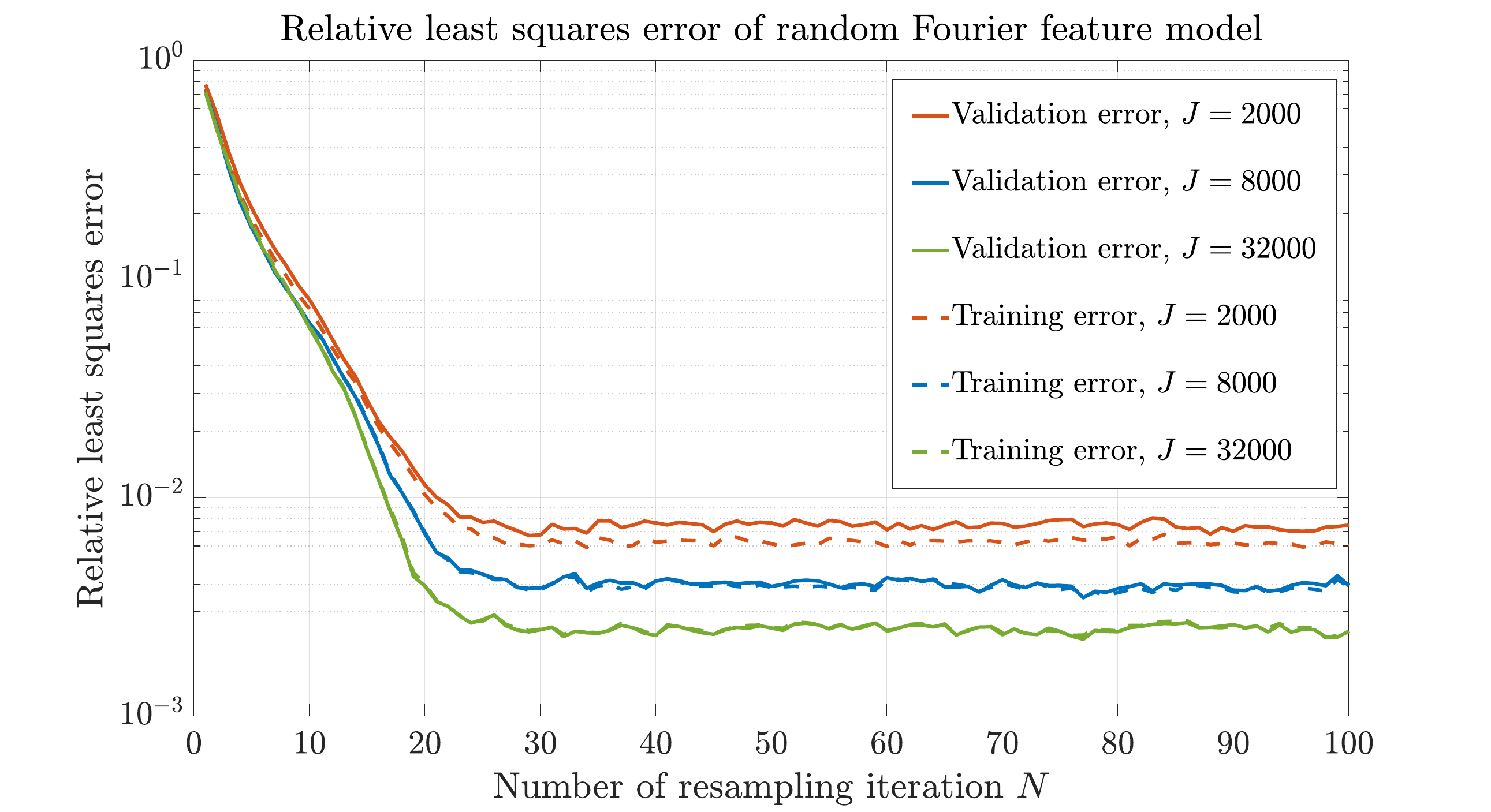}%
    \label{fig:vary_J_noisy}%
  }

  \caption{Relative least squares error of the trained random Fourier feature model $\beta(x)$ with varied training dataset sizes $J$ and fixed network size $K=2500$, using noiseless and noisy training data \(y_j = f(x_j) + \xi_j\), where the noise variables \(\xi_j\) are independently and identically distributed Gaussian random variables with mean zero and variance \(s^2\), with \(s = 0\) and \(s = 2.5 \times 10^{-3}\), respectively.}
  \label{fig:Fig_a3_vary_J}
\end{figure}

% \subsubsection{The training data set size J}
\Testcase{\label{tc:test4}
\textit{The effect of the training dataset size $J$.}
}
In Figure~\ref{fig:Fig_a3_vary_J}, we increase the number of training data points from $J=2000$ to $J=32000$ with a fixed network size $K=2500$ and plot the relative least squares error on the validation set. As expected, the validation error decreases progressively  for larger training set sizes.  Table~\ref{Table:Error_record_vary_J} summarizes the relative least squares errors on the training, validation, and test sets for varying 
$J$. A smaller discrepancy between the training and test errors is observed as $J$ increases, indicating alleviation of overfitting. 

Furthermore, experiments using mildly noisy training data  demonstrate improved consistency across the training, validation, and test errors, as depicted in Figure~\ref{fig:vary_J_noisy}. This result suggests that introducing low-level noise serves as an implicit regularization mechanism, mitigating overfitting by discouraging the model from fitting to spurious data patterns .
\begin{table}[!htbp]
    \centering
    \begin{tabular}{cccc}
        \toprule
        Sample Size \(J\) & Training Error & Validation Error & Test Error \\
        \midrule
        \multicolumn{4}{l}{\textbf{Noiseless data}} \\
        $2000$   & $5.80\times 10^{-3}$ & $7.00\times 10^{-3}$ & $7.02\times 10^{-3}$ \\
        $8000$   & $3.49\times 10^{-3}$ & $4.02\times 10^{-3}$ & $4.23\times 10^{-3}$ \\
        $32000$  & $2.50\times 10^{-3}$ & $2.67\times 10^{-3}$ & $2.53\times 10^{-3}$ \\
        [1.5ex] % extra vertical space
        \multicolumn{4}{l}{\textbf{Noisy data}} \\
        $2000$   & $6.04\times 10^{-3}$ & $7.47\times 10^{-3}$ & $7.61\times 10^{-3}$ \\
        $8000$   & $3.88\times 10^{-3}$ & $3.95\times 10^{-3}$ & $3.98\times 10^{-3}$ \\
        $32000$  & $2.47\times 10^{-3}$ & $2.44\times 10^{-3}$ & $2.54\times 10^{-3}$ \\
        \bottomrule
    \end{tabular}
    \caption{Relative least squares errors of the random Fourier feature model under different training dataset sizes \(J\). The top block uses noiseless data, while the bottom block uses noisy training data $y_j=f(x_j)+\xi_j$ where the noise variables \(\xi_j\) are independently and identically distributed as \(\mathcal{N}(0, s^2)\) with $s=2.5\times 10^{-3}$.}
    \label{Table:Error_record_vary_J}
\end{table}

The hyperparameters used in Table~\ref{Table:Error_record_vary_J} and Figure~\ref{fig:Fig_a3_vary_J} are as follows: $K=2500$, $\delta=0.2$, $\epsilon=\frac{1}{200}K^{-\frac{1}{2}}$, $\lambda_1=\frac{1}{20}KJ^{-\frac{1}{2}}$, $\lambda_2=0$, $d=4$, and $a=0.1$. The target function $f(x)$ is given in \eqref{eq:f} with a period of $q=12$, and the random Fourier feature model $\beta(x)$ is trained by applying Algorithm~\ref{alg1}.

\begin{figure}[!htbp]
    \centering
    \includegraphics[width=0.9\linewidth]{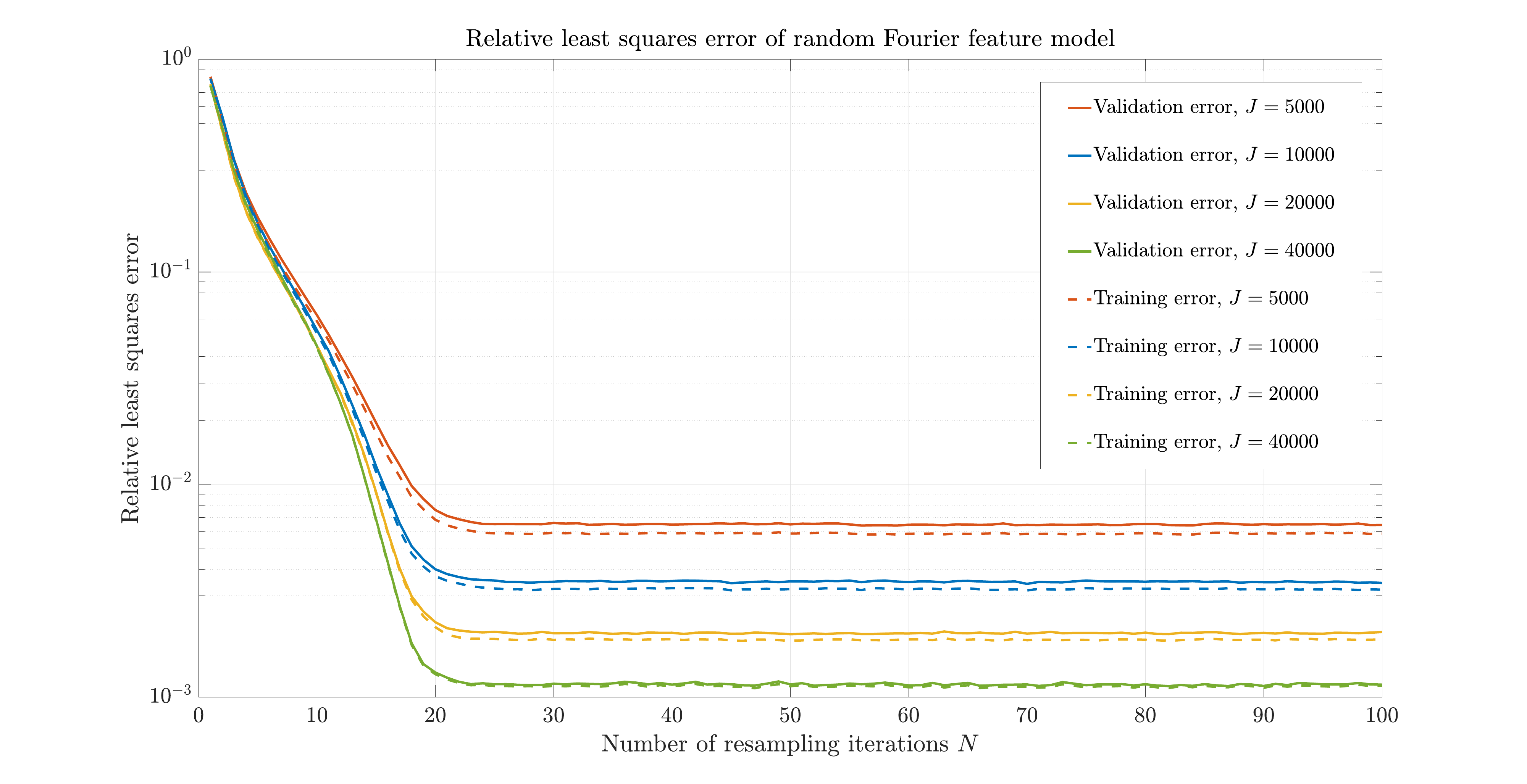}
    \caption{Relative least squares error of the trained random Fourier feature model $\beta(x)$ using varied training data set size $J$ with a fixed network size $K=20000$.}
    \label{fig:Fig_overparametrize_vary_J}
\end{figure}

To evaluate the performance of the proposed training Algorithm~\ref{alg1} in an over-parameterized regime, we apply the random feature model with a large network size $K = 20000$ and vary the noiseless training dataset  size $J$. In particular, we consider the cases where $J \leq K$, using $J = 5000$, $10000$, and $20000$.  Figure~\ref{fig:Fig_overparametrize_vary_J} depicts the training and validation errors of the resulting model $\beta(x)$. In this regime, the generalization error is primarily governed by the size of the training dataset, and Tikhonov regularization is crucial in ensuring numerical stability. Notably, for $J = 5000$ and $10000$, both the training and validation errors decrease steadily over resampling iterations, indicating reliable convergence despite over-parameterization. This behavior is consistent with that observed in the under-parameterized setting ($K < J$), demonstrating the robustness of the proposed algorithm across different parameter regimes.

\begin{figure}[!htbp]
    \centering
    \includegraphics[width=0.9\linewidth]{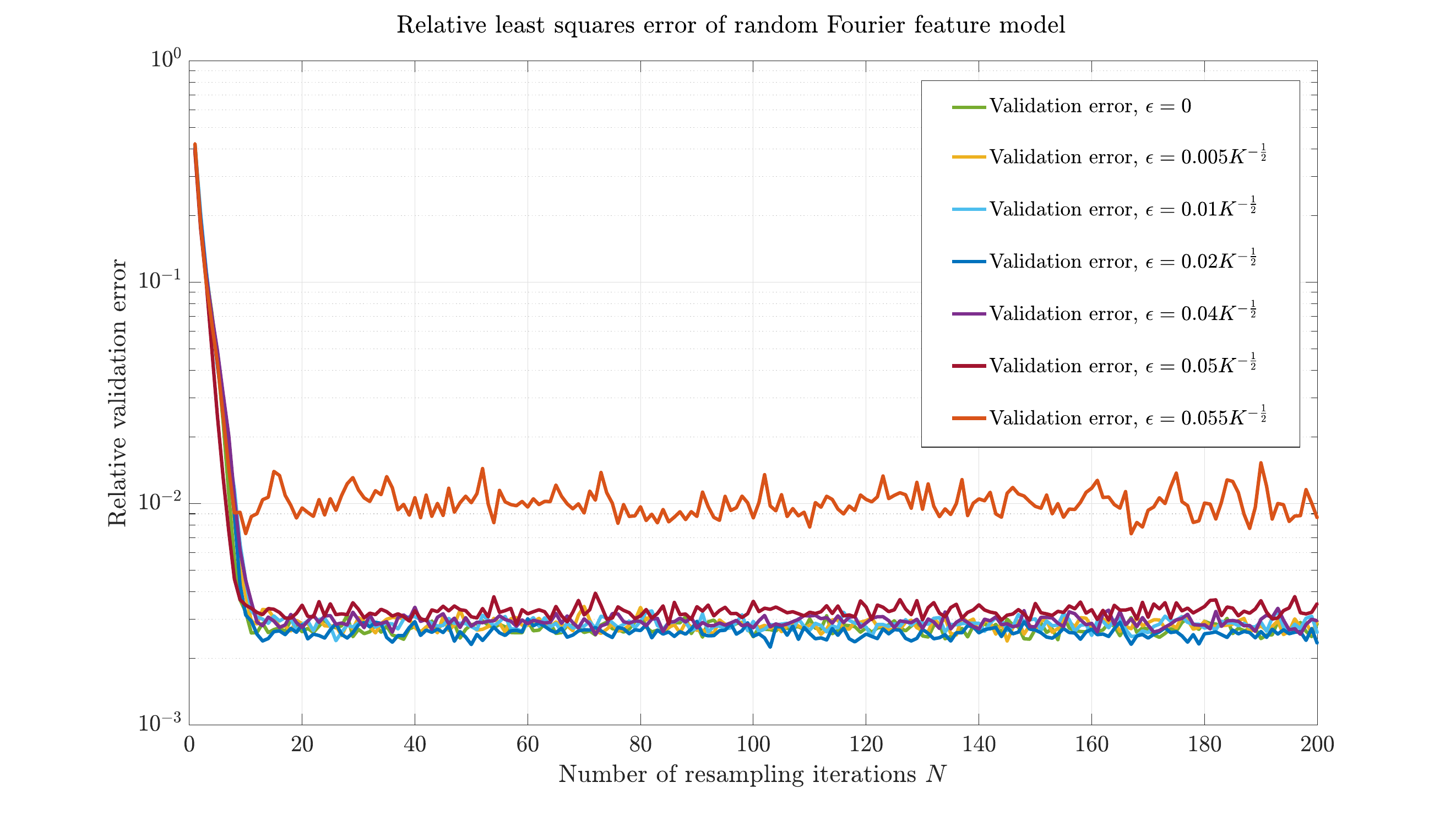}
    \caption{Relative least squares error of the trained random Fourier feature model $\beta(x)$ with varied cutoff parameter $\epsilon$.}
    \label{fig:Fig_a4_vary_epsilon}
\end{figure}

% \subsubsection{The cutoff parameter $\epsilon$}
\Testcase{\label{tc:test5}
\textit{The effect of the cutoff threshold $\epsilon$.}
}
This test case further investigates the effect of the cutoff parameter $\epsilon$ in Algorithm~\ref{alg1} by varying it from $\epsilon = 0$ to $\epsilon = \frac{1}{20}K^{-\frac{1}{2}}$. The results indicate that the algorithm remains effective even without a cutoff during resampling, whereas excessively large values of $\epsilon$ degrade the approximation quality.  Figure~\ref{fig:Fig_a4_vary_epsilon} presents the corresponding validation error decay curves.

The hyperparameters in this test case are $J=8000$, $K=2500$, $\delta=0.5$, $\lambda_1=\frac{1}{20}KJ^{-\frac{1}{2}}$, $\lambda_2=0$, $d=4$, and $a=0.1$. The target function $f(x)$ is given in \eqref{eq:f} with a period of $q=12$.

\begin{figure}[!htbp]
    \centering
    \includegraphics[width=0.9\linewidth]{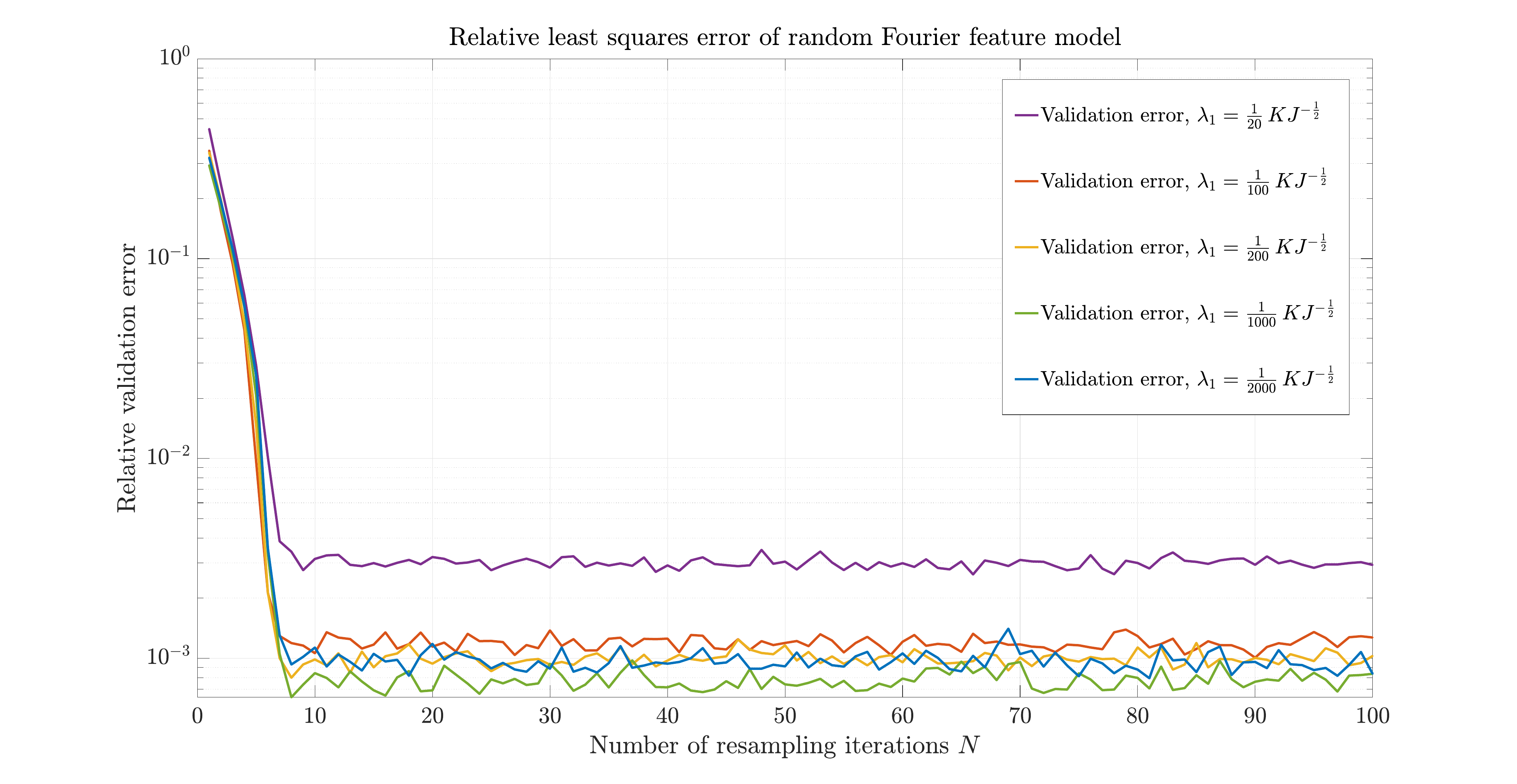}
    \caption{Relative least squares error of the trained random Fourier feature model $\beta(x)$ with varied Tikhonov regularization parameter $\lambda_1$.}
    \label{fig:Fig_vary_lambda1}
\end{figure}

\Testcase{\label{tc:test6}
\textit{The effect of the regularization weight $\lambda_1$.}
}
Tikhonov regularization is employed to alleviate overfitting and suppress the extreme variations in the regressor, with the regularization parameter $\lambda_1$ balancing data fidelity and solution smoothness. %In Figure~\ref{fig:Fig_vary_lambda1} we test with varied $\lambda_1$ regularization and observe a decreasing trend in the prediction error as $\lambda_1$ decreases, whereas too small value of $\lambda_1$ results in unstable and less accurate model performance.
 In Figure~\ref{fig:Fig_vary_lambda1}, we evaluate the model performance under varying values of $\lambda_1$ and observe reduced prediction error as $\lambda_1$ decreases. However, an excessively small value of $\lambda_1$ results in less stable prediction performance and reveals declining model accuracy. For example, the blue curve in Figure~\ref{fig:Fig_vary_lambda1}, corresponding to the smallest tested regularization level, displays a noticeably higher validation error than the green curve, which corresponds to a larger regularization.
 %compared to other tests using moderately small values of $\lambda_1$-penalization.

The hyperparameters in this test case are $J=8000$, $K=2500$, $\delta=0.5$, $\epsilon=\frac{1}{200}K^{-\frac{1}{2}}$, $\lambda_2=0$, $d=4$, and $a=0.1$, for the target function $f$ given in \eqref{eq:f} with a period of $q=12$.

% \subsubsection{The regularization parameter $\lambda_2$}
\Testcase{\label{tc:test7}
\textit{The effect of the regularization weight $\lambda_2$.}
}
With the regularization parameter $\lambda_2>0$, the minimization problem \eqref{min3}
\[
\min_{\hat{\beta}\in\mathbb{C}^K}\big( \frac{1}{J}\sum_{j=1}^J |\beta(x_j)-y_j|^2+\lambda_1\sum_{k=1}^K|\hat{\beta_k}|^2+\lambda_2 (\sum_{k=1}^K|\hat{\beta}_k|^2)^2 \big)
\]
leads to a nonlinear system of equations for the amplitude coefficients $\{\hat{\beta}_k \}_{k=1}^K$, as the derivative of the objective function is nonlinear in $\hat{\beta}_k$. Therefore, we apply Newton's method to solve for the amplitude coefficients by leveraging the convexity of the objective function. Further implementation details on Newton's method for solving the problem \eqref{min3} are available in the GitHub repository \cite{Code_repository} accompanying this manuscript.

\begin{figure}[!htbp]
    \centering
    \includegraphics[width=0.9\linewidth]{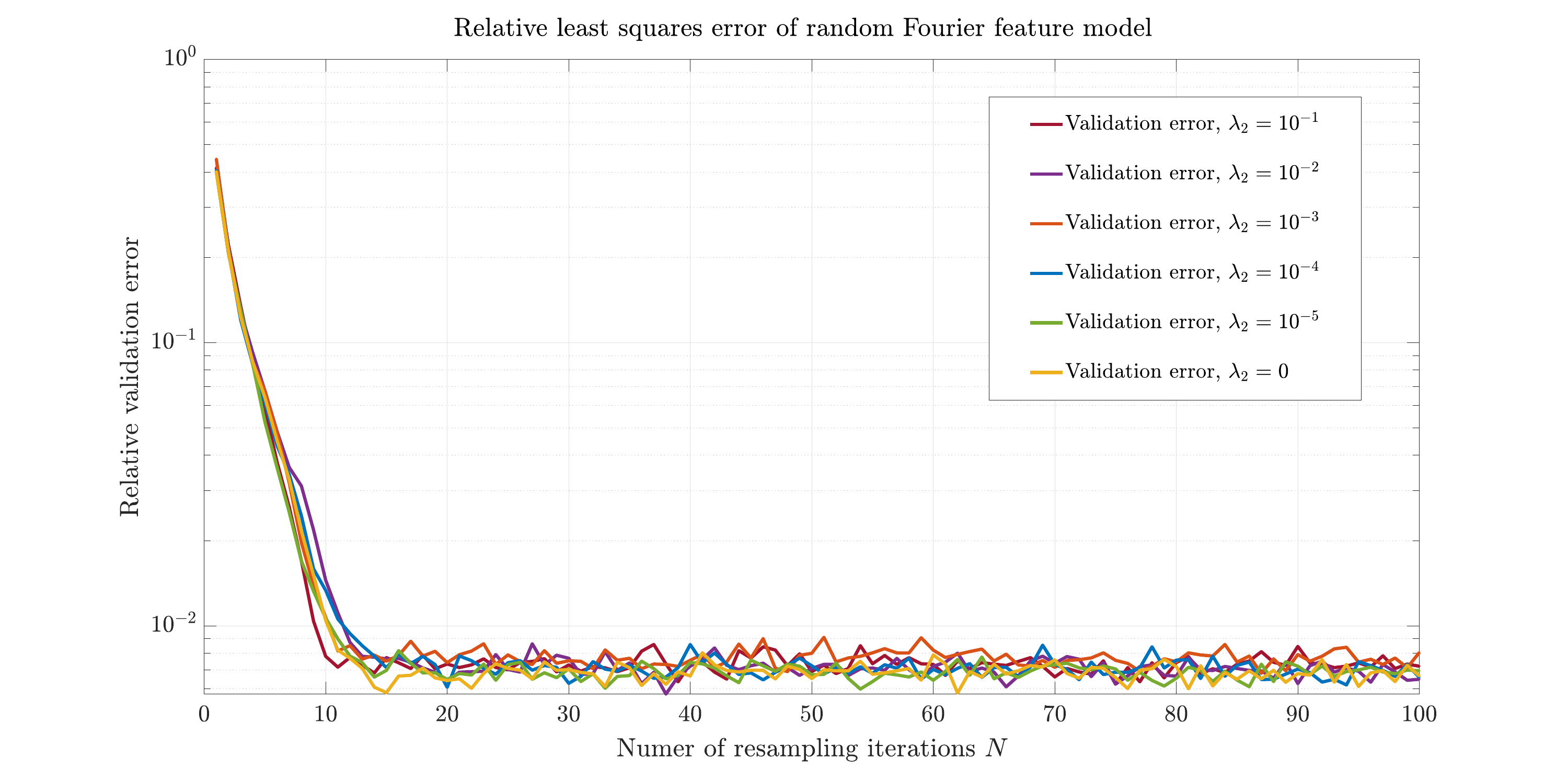}
    \caption{Relative least squares error of the trained random Fourier feature model $\beta(x)$ with varied regularization parameter $\lambda_2$.}
    \label{fig:Fig_a5_vary_lambda2}
\end{figure}

We solve the minimization problem \eqref{min3} with a gradually increasing regularization parameter $\lambda_2$, and observe the stability in the performance of the trained random Fourier feature model $\beta(x)$, as depicted in Figure~\ref{fig:Fig_a5_vary_lambda2}. The numerical test suggests that, in practice, setting $\lambda_2=0$ is sufficient.

The hyperparameters in this test case are $J=4000$, $K=1250$, $\delta=0.5$, $\epsilon=\frac{1}{200}K^{-\frac{1}{2}}$, $\lambda_1=\frac{1}{20}KJ^{-\frac{1}{2}}$, $d=4$, and $a=0.1$. The target function $f$ is given in \eqref{eq:f} with a period of $q=12$. For the implementation of Newton's method on the minimization problem \eqref{min3}, we use the Hessian of the objective function with respect to the real and imaginary part of $\hat{\beta}_k$, respectively, and the solution to the standard least-squares problem with $\lambda_2=0$ is employed as the initial guess for Newton's method.

\Testcase{\label{tc:test8}
\textit{The noise level parameter $s$ in the training data.}
}
In Test~\ref{tc:test4}, we observe that a small amount of noise helps alleviate the overfitting problem for training the random Fourier feature model $\beta(x)$. Following the setup introduced at the beginning of Section~\ref{sec:noise}, in this test case we vary the noise level in the training data  $y_j=f(x_j)+\xi_j$ incorporating independent noise random variables $\xi_j,\ j=1,\ldots, J,$ with mean $\mathbb E[\xi_j]=0$ and variance $\mathbb E[|\xi_j|^2]=s^2$. Specifically, each $\xi_j$ is sampled from a normal distribution $\mathcal{N}(0,s^2)$, where the standard deviation $s$ is gradually increased from $0.025$ to $0.1$.

\begin{figure}[!htbp]
    \centering
\includegraphics[width=0.9\linewidth]{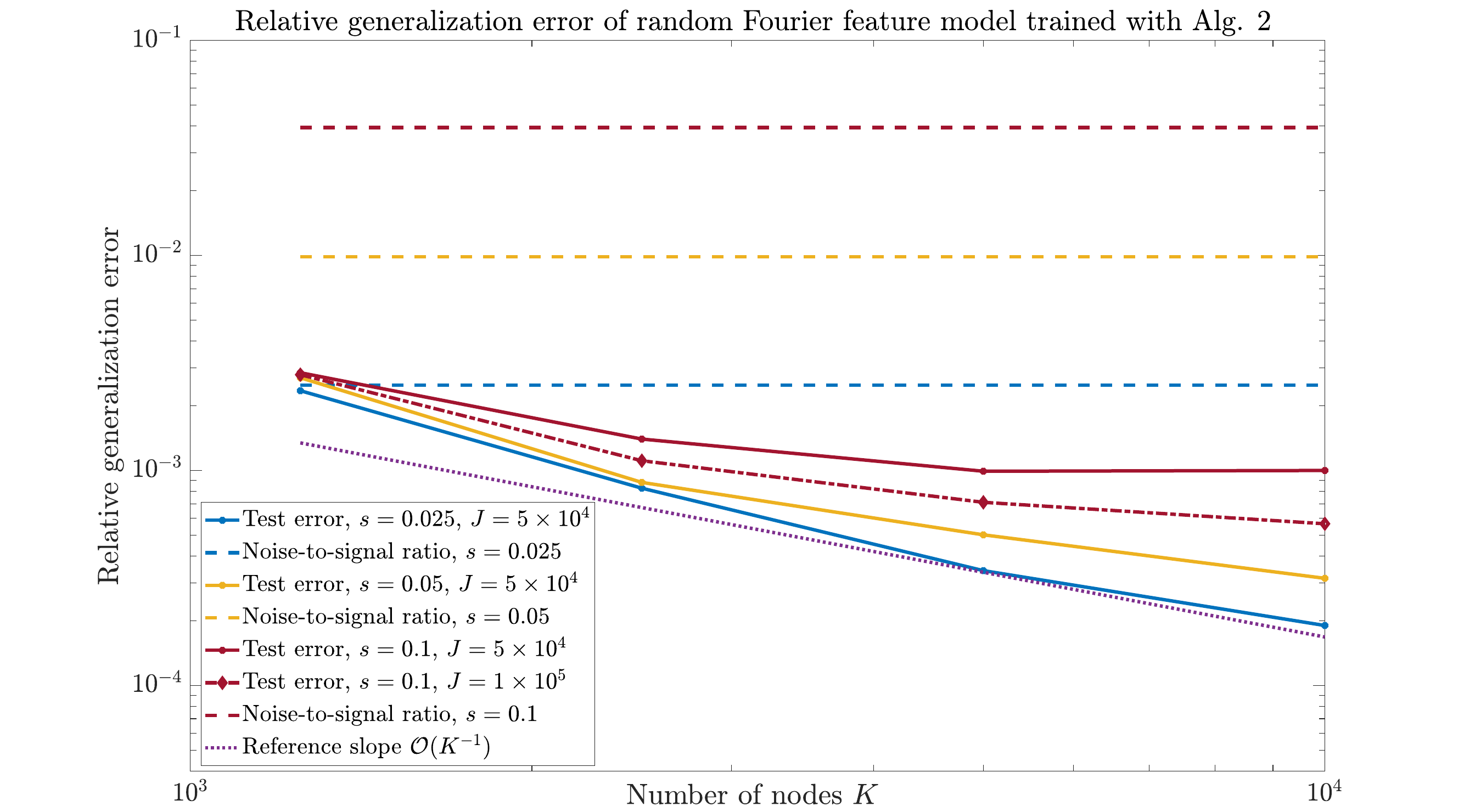}
    \caption{Relative generalization error of the trained random Fourier feature model $\beta(x)$ with varied noise level parameter $s$ and increased network size $K$. The dashed lines corresponding to noise-to-signal ratios under various noise levels are computed using \eqref{eq:def_NSR}. }
    \label{fig:Fig_vary_noise_s}
\end{figure}

Figure~\ref{fig:Fig_vary_noise_s} plots the generalization errors of the trained random Fourier feature model
\[
\frac{\sum_{j=1}^{\tilde{J}} |f(\tilde{x}_j)- \beta(\tilde{x}_j)|^2 }{\sum_{j=1}^{\tilde{J}}|f(\tilde{x}_j)|^2}
\]
on the test dataset $\{\tilde{x}_j,f(\tilde{x}_j)\}_{j=1}^{\tilde{J}}$ with different values of the noise level parameter $s$. We also record the corresponding noise-to-signal ratios (NSR) of the training dataset $\{x_j,y_j\}_{j=1}^J$,
\begin{equation}\label{eq:def_NSR}
\mathrm{NSR}=\frac{\sum_{j=1}^{{J}} |y_j- f({x}_j)|^2 }{\sum_{j=1}^{{J}}|y_j|^2}=\frac{\sum_{j=1}^{{J}}|\xi_j|^2}{\sum_{j=1}^{{J}}|y_j|^2}\,.
\end{equation}
The random Fourier feature model trained with resampling effectively performs denoising across all tested noise-to-signal ratios, $0.25\%$, $1.0\%$, and $4.0\%$, corresponding to standard deviations $s=0.025$, $0.05$, and $0.1$, respectively. The model achieves a relative generalization error of approximately one order of magnitude smaller than the corresponding NSR values. Moreover, the model exhibits good scalability with respect to the number of nodes $K$ under mild noise levels. Under high noise levels, the generalization error continues to decrease as the training dataset size $J$ increases.

In Figure~\ref{fig:Fig_vary_noise_s}, we apply Algorithm~\ref{alg1} for the target function $f$ given in \eqref{eq:f} with a period of $q=12$ and dimension $d=4$, $a=0.1$. We use hyperparameters $\delta=0.5$, $K=10000$, $J=50000$, $\epsilon=\frac{1}{200}K^{-\frac{1}{2}}$, $\lambda_1=\frac{1}{100}KJ^{-\frac{1}{2}}$, and $\lambda_2=0$, with a test dataset size of $\tilde{J}=0.5J$.

\subsection{Numerical implementation: classification on MNIST dataset}
We perform classification experiments on the MNIST handwritten digits dataset using Algorithm \ref{alg_FT} with two modifications. The algorithm runs for 6000 iterations for each of ten neural networks, where each network is dedicated to classifying a single digit.

The original MNIST training dataset of 60,000 images is partitioned into 50,000 training samples and 10,000 validation samples. Evaluation is performed on a separate test set of 10,000 samples. The Tikhonov regularization parameter is set to \( \lambda = 2 \), and the step size is \( \delta = 0.005 \). The conjugate gradient method, with a tolerance of \( 10^{-4} \), is used to solve the linear least squares problems. For the evaluation, we select the neural network weights that correspond to the best validation accuracy across all digits.

Each neural network, \( \beta^n: \mathbb{R}^{784} \to \mathbb{R} \) for \( n = 0,1,\dots,9 \), is defined as follows:
\[
\beta^n(x) = \sum_{k=1}^{K} \left(\hat\beta_k^n\cos(\omega_k^n\cdot x) + \check\beta_k^n\sin(\omega_k^n\cdot x)\right),
\]
where \( \omega_k^n, x \in \mathbb{R}^{784} \), \( \hat\beta_k^n, \check\beta_k^n \in \mathbb{R} \), and \( K = 10000 \). 

The two modifications to Algorithm \ref{alg_FT} are adjusting the network structure as described above and modifying the resampling probability to  
\[
\frac{\|(\hat\beta_k^n, \check\beta_k^n)\|_2}{\sum_{\ell} \|(\hat\beta_\ell^n, \check\beta_\ell^n)\|_2}.
\]

For classification, each handwritten test digit is assigned to class \( n \), where \( n \) corresponds to the neural network \( \beta^n \) that produces the highest output value.

\begin{figure}[!htbp]
    \centering
    \includegraphics[width=0.3\textwidth]{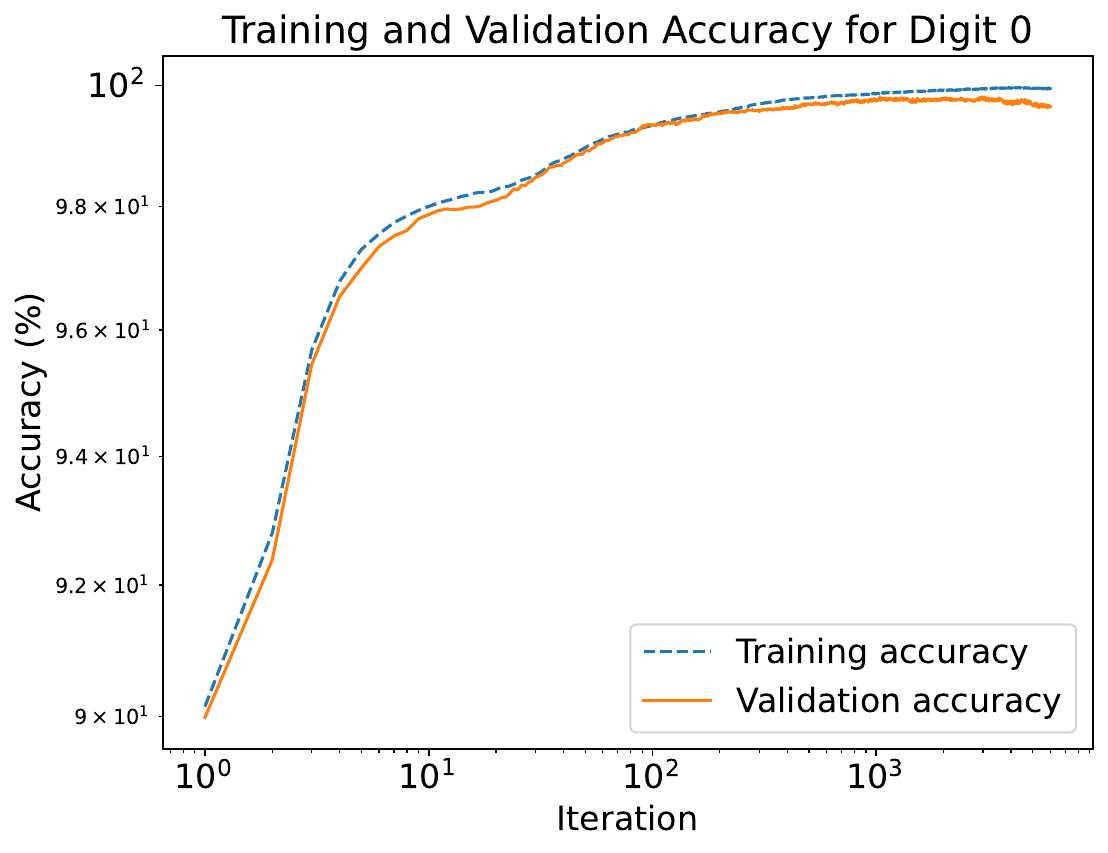}
    \includegraphics[width=0.3\textwidth]{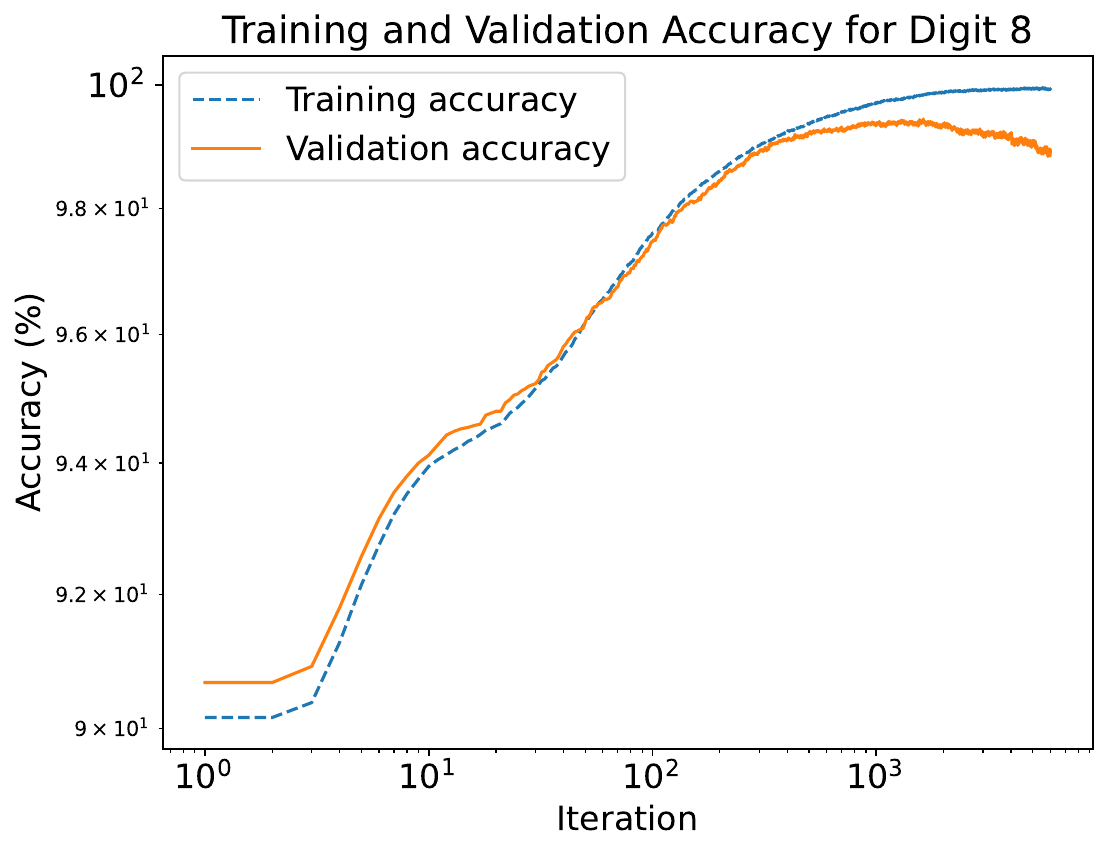}
    \includegraphics[width=0.3\textwidth]{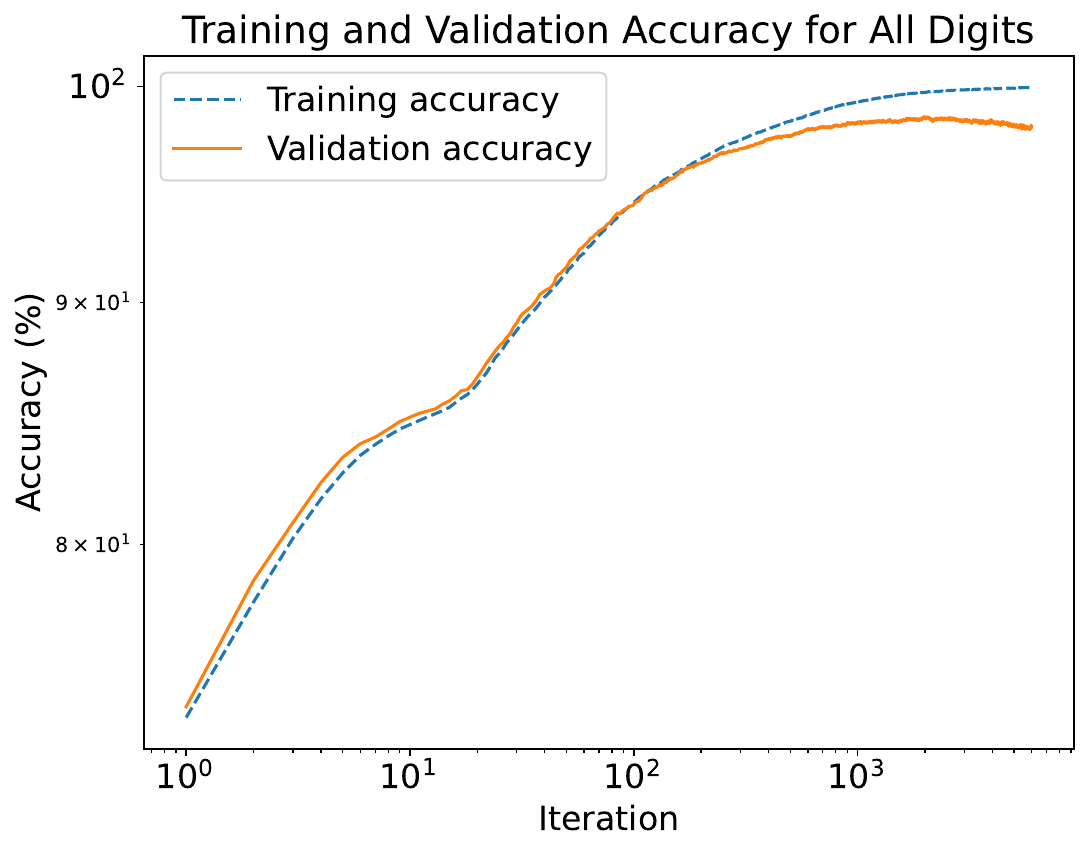}
    \caption{Training and validation accuracy for digit 0, digit 8, and the overall classification accuracy when classifying all digits using the ten neural networks.}
    \label{fig:MNIST_accuracy_no_pca}
\end{figure}

Figure \ref{fig:MNIST_accuracy_no_pca} presents the training and validation accuracy for digit 0 and digit 8, and the overall classification accuracy when using all ten neural networks.

At iteration 1974, the best validation accuracy is achieved, with a training accuracy of 99.71\%, validation accuracy of 98.54\%, and test accuracy of 98.26\%. The neural network with the highest accuracy at this iteration is the one responsible for digit 0, whereas the lowest accuracy is observed for digit 8.

\section{Other decay rates}
\label{other}
For the decay 
\begin{equation}\label{rate1}
|\hat f(\omega_n)|\le (\frac{C}{|n|})^\ell\,,
\end{equation}
we obtained in \eqref{N_*}
\[\sum_{|\hat f(\omega_n)|\ge\epsilon/2} 1 \le 
\sigma_d \frac{N_*^d}{d} =\frac{2^{d/\ell}\sigma_d}{d} \frac{C^d}{\epsilon^{d/\ell}}=\mathcal O(K^{\frac{d/\ell}{3-d/\ell}})\,
\]
In the case that $C$ is large the factor $C^d$ implies that the number of nodes $K$ need to be  large,  depending on $d$, to make 
$\sigma_d \frac{N_*^d}{d} K^{-1/2}$ small, which is required in \eqref{cheby}.

Consider therefore a case where  $C$ is  large only in a few directions $(n_1,\ldots, n_b)$, where the decay rate also can be slow.
%is not large in all directions except one  and in this direction the decay rate is also slow, 
Assume for instance we have the separable setting with 
\[
\omega_n=\frac{\pi}{L}(\underbrace{n_1,\ldots, n_b}_{=:n'}, \underbrace{n_{b+1},\ldots,n_d}_{=:n^\perp})=:(\omega'_n,\omega^\perp_n)\]
with $b\ll d$
and
\begin{equation}\label{f'}
|\hat f(\omega_n)|= |\hat f_\perp(\omega^\perp_n)||\hat f'(\omega'_{n})|
\end{equation}
where $|\hat f_\perp(\omega^\perp_n)|\le  (\frac{C_\perp}{|n^\perp|})^\ell$ decays  fast with the rate $\ell>d$ and a small constant $C_\perp$\,,
while $|\hat f'(\omega'_n)|\le  (\frac{C'}{|n'|})^{\ell'}$ decays slowly with $b\le\ell'\ll d$ and a possibly  large constant $C'$.

The condition $|\hat f(\omega_n)|= |\hat f_\perp(\omega^\perp_n)||\hat f'(\omega'_{n})|>\epsilon$ implies that
\[
\sigma_b |n'|^b < \sigma_b (C' \epsilon^{-1/\ell'})^b (\frac{C_\perp}{|n_\perp|})^{b\ell/\ell'}\,,
\]
and for $b(1+\ell/\ell')>d$ we obtain
\[
\begin{split}
\sum_{|\hat f(\omega_n)|>\epsilon} 1
&\le \sigma_b \sigma_{d-b}(C' \epsilon^{-1/\ell'})^b \int_1^\infty (\frac{C_\perp}{|n_\perp|})^{b\ell/\ell'}|n_\perp|^{d-b-1}{\rm d}|n_\perp|\\
&=\sigma_b \sigma_{d-b}C'^b \epsilon^{-b/\ell'} \frac{C_\perp^{b\ell/\ell'}}{b(1+\frac{\ell}{\ell'})-d}\,,
%
%\le \sum_{|\hat f_\perp(\omega^\perp_n)|>\epsilon/\min_{n'}|\hat f_1(\omega'_{n})|}1
%+  \sum_{|\hat f_1(\omega'_{n})|>\epsilon/\min_{n^\perp}|\hat f_\perp(\omega^\perp_{n})|}1\\
%&\le \sigma_{d-b} C_\perp^{d-b} \epsilon^{-(d-b)/\ell} + \sigma_b C'^{b}  \epsilon^{-b/\ell'}\\
%&\le \sigma_{d-b} C_\perp^{d-b} K^{(d-b)/(2\ell)} + \sigma_b C'^b K^{b/(2\ell')}\,,
\end{split}
\]
%and the resolving condition becomes
%\[
%\sigma_{d-b} C_\perp^{d-b} K^{(d-b)/(2\ell)} + \sigma_b C'^b K^{b/(2\ell')}\le \gamma^{-1/4}\mathcal O(K^{1/2})\,,
%\]
which does not include a large factor $C^d$ for $d\gg 1$ when $b\ll d$.
%
%Therefore the rate constants $C_p$ during the iterations are bounded and
%the function $f_\epsilon$ is  approximated with asymptotically optimal rate constant.

The remaining cutoff error  has with the slow decay \eqref{f'}, in a few directions $b$, for $\ell'>b$ and $\ell> d-b$ the bound
\[
\begin{split}
\sum_{|\hat f(\omega_n|\le \epsilon} |  \hat f(\omega_n| &\le \sum_{|\hat f(\omega_n|\le \epsilon}  (\frac{C'}{|n'|})^{\ell'}
 (\frac{C_\perp}{|n_\perp|})^\ell\\
 &\le  \sum_{(\frac{C'}{|n'|})^{\ell'}C_\perp^\ell \le \epsilon} (\frac{C'}{|n'|})^{\ell'} \sum_{n_\perp\in\mathbb Z^{d-b}\setminus 0}(\frac{C_\perp}{|n_\perp|})^\ell
 + \sum_{ (\frac{C_\perp}{|n_\perp|})^\ell C'^{\ell'}\le \epsilon}(\frac{C_\perp}{|n_\perp|})^\ell \sum_{n'\in\mathbb Z^b\setminus 0} (\frac{C'}{|n'|})^{\ell'}\\
&=\mathcal O(\epsilon^{1-b/\ell'} + \epsilon^{1-(d-b)/\ell})\,,
%K^{-1}\mathcal O(K^{b/(2\ell')} C'^b+K^{(d-b)/(2\ell)} C_\perp^{d-b} )\,.
\end{split}
\]
which also avoids constant large factors $C^d$.

The remainder factor in the proof of Theorem \ref{theorem1} becomes
\[
\mathbb E_p[\|\eta\|] \mathcal O(\epsilon^{-1} +\epsilon^{-b/\ell'}) + 
q_\epsilon^{-1}\mathcal O(\epsilon^{1-b/\ell'} + \epsilon^{1-(d-b)/\ell} )+\mathcal O(q_\epsilon)\,,
\]
so that  the estimate in \eqref{conclude} also holds for the anisotropic decay in \eqref{f'} by replacing $d/\ell$ in Theorem \eqref{theorem1}  by $b/\ell'$, in the case   $b/\ell'>(d-b)/\ell$.

In conclusion, the conditions for convergence of the resampling in Theorem \ref{theorem1}
can be applied to a setting with slow decay in a few directions \eqref{f'} and fast decay in the remaining several directions.
%
%holds for the slower decay rate \eqref{f'} provided the actual function is cut off, that is by approximating $f_\epsilon$ instead of $f$.
%In particular the resolving conditions \eqref{N_bound} then holds for $K$ independent of $d$ 
%for some functions that decay slowly in a few directions
%and fast in the remaining several directions. 
In some sense, this is the setting where we can expect neural network approximation to  work well because slow decay in the majority of directions leads to a large optimal rate constant $C_{p_*}=(\sum_{n\in\mathbb Z^d}|\hat f(\omega_n)|)^2=\mathcal O(C^{2\ell})$ for $\ell>d\gg 1$ and $C\gg 1$.

\subsection{Smaller generalization error}\label{small_GE}
The decay in \eqref{decay} implies a smaller generalization error  $\min_{(\hat\beta,\nu)\in\mathbb C^K\times \frac{\pi}{L}
\mathbb Z^{dK} }\mathbb E_x[|f(x)-\beta(x)|^2]= \mathcal O(K^{-2\ell/d+1})$, 
for the case with an infinite amount of data $J=\infty$,
%of the generalization error 
by sampling the $K$ frequencies with the largest amplitudes as follows.
By choosing $\sum_{|\hat f(\omega_n)|\ge \epsilon} 1=K$ we obtain by \eqref{N_*}  that $\epsilon=\mathcal O(K^{-\ell/d})$.
Consequently  $\beta$ becomes equal to the cutoff function $f_\epsilon$.  Assume that 
$\|\rho\|_{L^\infty(\mathbb T^d)}$ is bounded; thus, by the orthogonality of the 
Fourier basis, the remaining error is 
\begin{equation}\label{fast}
\begin{split}
\min_{(\hat\beta,\nu)\in\mathbb C^K\times \frac{\pi}{L}
\mathbb Z^{dK} }\mathbb E_x[|f(x)-\beta(x)|^2] &\le \|\rho\|_{L^\infty}\min_{(\hat\beta,\omega)\in\mathbb C^K\times\mathbb R^K} \int_{\mathbb T^d} |f(x)-\beta(x)|^2{\rm d} x\\
 &=\|\rho\|_{L^\infty} \int_{\mathbb T^d} |f(x)-f_\epsilon(x)|^2{\rm d} x\\
 &=(2L)^d\|\rho\|_{L^\infty}\sum_{|\hat f(\omega_n)|\le \epsilon}|\hat f(\omega_n)|^2\\
&= \mathcal O(\epsilon^{2-d/\ell})\\
&= \mathcal O(K^{1-2\ell/d})\,,
\end{split}
\end{equation}
which is smaller than the generalization error estimate $\mathcal O(K^{-1})$ obtained in Theorem \ref{theorem1} for $\ell>d$.
Therefore, the error estimate in Theorem \ref{theorem1} is not sharp for functions with the decay in \eqref{decay}. 
An algorithm providing this smaller generalization error \eqref{fast}, for $\ell>d$, needs to find the $K$ frequencies with the largest amplitudes, which would require a different analysis of the non-convex problem for the frequencies.  The computational work for an algorithm to achieve this small minimal generalization error can typically be large.

\section{The generalization error for random features with finite amount of data and general activation functions}
\label{sec_general}
In this section, we derive the error estimates for the generalization error in the case of finite amount of data and for general activation functions.
Let for $x\in\rset^d$
\[
\tilde x:= \left[\begin{array}{c}
x\\
1\\
\end{array}
\right]
\]
and assume a  function $f\in L^\infty(\rset^d)$ has the following representation:
\begin{equation}\label{f_repr}
f(x)=\int_{\rset^{d+1}}a(\omega)\sigma(\omega\cdot\tilde x){\rm d}\omega\,,
\end{equation}
for a function $a:\rset^{d+1}\to\rset$ and a given nonpolynomial  activation function $\sigma:\rset\to\rset$.
For example, the activation function $\sigma(\omega\cdot\tilde x):=\cos(\omega\cdot\tilde x)$ yields by the Fourier transform $\hat f:\rset^d\to\rset$, defined by $\hat f(\omega):=(2\pi)^{-d/2}\int_{\rset^d}f(x)e^{-{\rm i}\omega\cdot x}{\rm d}x$ for $f\in L^1(\rset^d)$,
\[
\begin{split}
f(x)&= {\rm Re}\int_{\rset^d} \hat f(\omega)e^{{\rm i}\omega\cdot x}{\rm d}\omega\\
&=\int_{\rset^d} |\hat f(\omega)| \cos\big(\omega\cdot x +{\rm arg}\hat f(\omega)\big){\rm d}\omega\\
&=\int_\rset\int_{\rset^d} |\hat f(\omega)| \cos\big(\omega\cdot x +\omega_0\big)
\delta\big(\omega_0-{\rm arg}\hat f(\omega)\big){\rm d}\omega{\rm d}\omega_0\,.\\
\end{split}
\]
For other activation functions, we do not have explicit formulations for which functions have a representation \eqref{f_repr}, but the Universal Approximation Theorem \cite{activation} proves that continuous functions can be approximated in maximum norm on compact sets by neural networks based on any nonpolynomial activation function to any accuracy.

In the next section, we consider, for  given nonnegative regularization parameters $\lambda_1$ and $\lambda_2$, the generalization error for the optimization problem
\begin{equation}\label{min_J0}
\begin{split}
&\min_{\hat\beta\in \rset^K } %\lim_{J\to\infty}
\Big(J^{-1}\sum_{j=1}^J|f(x_j)-\beta(x_j)|^2 +\lambda_1\sum_{k=1}^K|\hat\beta_k|^2  +\lambda_2(\sum_{k=1}^K|\hat\beta_k|^2)^2\Big)
\,,\\
\end{split}
\end{equation}
based on the random feature approximation
\begin{equation}\label{bar_f}
\beta(x):=\sum_{k=1}^K\hat\beta_k \sigma(\omega_k\cdot\tilde x)
\end{equation}
where $\{\omega_k\in\rset^{d+1} \ |\ k=1\ldots K\}$ are independent random samples from a probability distribution
$p:\rset^{d+1}\to[0,\infty)$.
 The given data $\{x_j\in\rset^d\ |\ j=1,\ldots J\}$ are independent random samples from some
unknown probability distribution $\mu:\rset^d\to[0,\infty)$. We use the notation $\mathbb E_\omega$  for the expected value with respect to $\omega$  and $\mathbb E_{\{x_j\}}$ for the expected value with respect to the data together with $\mathbb E_{\mathbf x}$ for the expectation with respect to
the random variable $\mathbf x$ with probability distribution $\mu$.

Section \ref{RFnoise} includes independent noise $\xi_j$ in the random feature training \[
\min_{\hat\beta\in \rset^K } %\lim_{J\to\infty}
\Big(J^{-1}\sum_{j=1}^J|f(x_j)+\xi_j-\beta(x_j)|^2  +\lambda_1\sum_{k=1}^K|\hat\beta_k|^2
 +\lambda_2(\sum_{k=1}^K|\hat\beta_k|^2)^2\Big)\]
and establishes a corresponding generalization error bound.
\subsection{The generalization error for random features without noise}\label{RF}
This section derives the following estimate of the generalization error for the random  feature optimization \eqref{min_J0} with a finite amount of data.
\begin{theorem}[Generalization error without noise]\label{thm:generalization} Assume  
that $f:\rset^d\to\rset$ has the representation \eqref{f_repr}, where
\begin{equation}\label{2moment}
\begin{split}
&\int_{\rset^{d+1}}|a(\omega)|\big(1+(\mathbb E_{\mathbf x}[\sigma^2(\omega\cdot\mathbf x)])^{1/2}
\big){\rm d}\omega =\mathcal O(1)\,,\,,\\
%&\mathbb E_\omega\big[\mathbb E_{\mathbf x}[\sigma^2(\omega\cdot \mathbf x) +\sigma^4(\omega\cdot \mathbf x)]\big] =\mathcal O(1)\,,\\
C_1 &:=
\mathbb E_\omega 
\big[\mathbb E_{\{x_j\}}[\big(
\sigma(\omega_k\cdot x_j)f(x_j)-\mathbb E_{\mathbf x}[\sigma(\omega_k\cdot \mathbf x)f(\mathbf x)]\big)^2]\big]=\mathcal O(1)\,,\\
E_{k,\ell} &:=\mathbb E_{\mathbf x}[\sigma(\omega_k\cdot \tilde x)\sigma(\omega_\ell\cdot \tilde x)]=\mathcal O(1)\,,\\
C_2' &:=\mathbb E_\omega\big[\mathbb E_{\mathbf x}[\big( \sigma(\omega_k\cdot \tilde x)\sigma(\omega_\ell\cdot \tilde x) - E_{k,\ell}\big)^2]\big]=\mathcal O(1)\,,\quad k\ne\ell\,,\\
C_2'' &:= \mathbb E_\omega\big[
\mathbb E_{\mathbf x}[\big( \sigma^2(\omega_k\cdot \tilde x) - E_{k,k}\big)^2]\big]=\mathcal O(1)\,,\\
%
%\big(1+\mathbb E_{\mathbf x}[ \sigma^4(\omega\cdot \mathbf x)] 
\end{split}
\end{equation}
and that the sampling distribution $p:\rset^{d+1}\to[0,\infty)$ satisfies
\begin{equation}\label{ap_bound}
\begin{split}
%&\|A\|_{L^\infty} + \|A'\|_{L^\infty} + \|A''\|_{L^\infty} +\|B\|_{L^\infty} =\mathcal O(1)\,,\\
&\sup_{\omega\in\rset^{d+1}}\frac{|a(\omega)|\big(1+(\mathbb E_{\mathbf x}[\sigma^2(\omega\cdot\mathbf x)])^{1/2}\big)}{p(\omega)}=C''=\mathcal O(1) \,.
\end{split}
\end{equation}% assumptions in Theorem \ref{} holds.
Then, the random feature optimization in \eqref{min_J0} and \eqref{bar_f}
has the expected generalization error estimate
\begin{equation}\label{gen_O3}
\begin{split}
\mathbb E_{\{x_j\}}\big[\mathbb E_\omega [\int_{\rset^{d}} |f(x)-\beta(x)|^2\mu(x){\rm d}x ] \big]
&\le (1+\delta) \frac{C_p}{K}(1+\lambda_1)+\frac{\lambda_2(C_p^2+K^{-1}C_p(C'')^2)}{K^2} \\
&\quad +  \frac{1}{2\delta\lambda_2}( \frac{K^2-K}{J}C_2' + \frac{K}{J}C_2'')
+\frac{KC_1}{2\lambda_1 \delta J}\\
%&=\mathcal O\big(\frac{1}{K}+\lambda +(1+\frac{1}{\lambda K})(\sqrt{\frac{\log\frac{K}{\delta}}{J}})+\delta\big)\\
&=\mathcal O\big(\frac{1}{K}+\frac{\lambda_1}{K} +
\frac{K}{\lambda_1 J}+\frac{\lambda_2}{K^2}+\frac{K^2}{\lambda_2 J}\big)\,.
\end{split}
\end{equation}
\end{theorem}

We note that the choice of $\lambda_2= CK^2 J^{-1/2}$ and $\lambda_1= CKJ^{-1/2}$, for some positive constant $C$, yields
\[
\mathbb E_{\{x_j\}}\big[\mathbb E_\omega [%\min_{\eta\in \mathcal C^K } %\lim_{J\to\infty}
\int_{\rset^{d}}|f(x)-\beta(x)|^2 \mu(x){\rm d}x]\big]
=\mathcal O\big(\frac{1}{K}+\frac{1}{\sqrt J}\big)\,,
\]
and the constant in $\mathcal O$ in \eqref{gen_O3} depends only on the $\mathcal O$ constants in \eqref{2moment} and \eqref{ap_bound}.
\begin{proof}[Proof of the generalization error] The proof has four steps:
\begin{itemize}
\item[Step 0.] Formulates the training optimization and introduces the notation.
\item[Step 1.] Estimates the training error and formulates the generalization error.
\item[Step 2.] Derives a representation of the generalization error, 
%that separates the data error and the trained amplitudes $\hat\beta$, which depend on the data, 
using the regularization terms and the independence of the data points $x_j$ and the $\omega_k$ samples.
% Relates the generalization error to the regularization terms and applies independence of the data points to estimate the remainder. %Uses Hoeffding's inequality combined with the regularization term to estimate the data error.
\item[Step 3.] Combines Steps 1 and 2 to obtain  \eqref{gen_O3}.
%\item[Step 5.] Combines Steps 1-4 to obtain  \eqref{gen_O3}.
\end{itemize}

{\it Step 0.}
We have a finite number of data points
\[\{\big(x_j,f(x_j)\big)\ |\ j=1,\ldots,J\}\,,\]
where $x_j$ are independent samples from $\mu$, 
and  $\beta(x)=\sum_{k=1}^K\hat\beta_k\sigma(\omega_k\cdot \tilde x)$ is a solution to the optimization problem \eqref{min_J0}.
  The loss function
\[
\begin{split}
L\big(\beta(x),f(x)\big):=|\beta(x)-f(x)|^2 \,,
\end{split}
\]
satisfies
\begin{equation}\label{LL}
\begin{split}
&L\big(\beta(x),f(x)\big)=\sum_{k=1}^K\sum_{\ell=1}^K\hat\beta_k\hat\beta_\ell\sigma(\omega_k\cdot \tilde x)\sigma(\omega_\ell\cdot \tilde x)
-2\sum_{k=1}^K\hat\beta_k\sigma(\omega_k\cdot \tilde x)f(x)+|f(x)|^2\,.\\
\end{split}
\end{equation}
%Assume that $\mathbf x\in\rset^d$ is a random variable with probability distribution $\mu$ and 
We introduce the following notation for the mean and empirical mean with respect to the data
\[
\begin{split}
\mathbb E_{\mathbf x}[L\big(\beta(\mathbf x),f(\mathbf x)\big)]&:=\int_{\rset^{d}}L\big(\beta(x),f(x)\big)\mu(x){\rm d}x\,,\\
\hat{\mathbb E}_J[L\big(\beta(x_j),f(x_j)\big)]&:=\frac{1}{J}\sum_{j=1}^JL\big(\beta(x_j),f(x_j)\big)\,.\\
\end{split}
\]
The regularization terms will provide bounds for $\sum_{k=1}^K|\hat\beta_k|^2$ and $(\sum_{k=1}^K|\hat\beta_k|^2)^2$, so that we can separate the data approximation error
\[|\hat{\mathbb E}_J[\sigma(\omega_k\cdot \tilde x_j)\sigma(\omega_\ell\cdot \tilde x_j)]
- \mathbb E_{\mathbf x}[\sigma(\omega_k\cdot \tilde{\mathbf x})\sigma(\omega_\ell\cdot \tilde{\mathbf x})]|\]
 in \eqref{LL}.
%In Step 2 we rewrite  the second right hand side term in \eqref{LL} to also also relate to the least squares Gram matrix $J^{-1}\sum_{j=1}^J\sigma(\omega_k\cdot \tilde x_j)\sigma(\omega_\ell\cdot \tilde x_j)$.

{\it Step 1.} % (Neural network Fourier representation and training error).} 
Based on the representation \eqref{f_repr}, the amplitude coefficients
\begin{equation}\label{eta_F}
\hat a_k := \frac{a(\omega_k)}{Kp(\omega_k)}\,,\ k=1,\ldots, K\,,
\end{equation}
define the $f(x)$-%minimal variance 
approximation
\[
\Bw(x) := \sum_{k=1}^K \hat a_k \sigma(\omega_k\cdot \tilde x)\,,
\]
that satisfies the ortogonality relation
\begin{equation}\label{zero_mean}
\mathbb E_\omega[\Bw(x)-f(x)]=0,
\end{equation}
which is a consequence of  \eqref{eta_F} since
\begin{equation}\label{mean_ok}
\begin{split}
\mathbb E_\omega[\Bw(x)]&= \mathbb E_\omega[\sum_{k=1}^K   
\frac{a(\omega_k)\sigma(\omega_k\cdot \tilde x)}{Kp(\omega_k)}] \\
&=\int_{\mathbb R^{d+1}}a(\omega)  \sigma(\omega\cdot \tilde x){\rm d}\omega\\ 
&=f(x)\,.
\end{split}
\end{equation}

The mean \eqref{mean_ok} and the independence of $\omega_k$ imply that
\begin{equation*}\label{var_f}
\begin{split}
    \mathbb E_\omega [|\Bw(x)-f(x)|^2] 
    %&= 
    %\mathbb E_\omega [\sum_{k=1}^K\sum_{\ell=1}^K 
    %\big(\hat\eta_ke^{{\rm i} \omega_k\cdot x} - v_i(x)\big)
    %\big(\hat\eta_\ell e^{{\rm i} \omega_\ell\cdot x} - v_i(x)\big)^*]\\
    &=\frac{1}{K^2}\sum_{k=1}^K\sum_{\ell=1}^K
    \mathbb E_\omega [
    \big(\frac{a(\omega_k)}{p(\omega_k)}\sigma(\omega_k\cdot \tilde x)-f(x)\big)
     \big(\frac{a(\omega_\ell)}{p(\omega_\ell)}\sigma(\omega_\ell\cdot \tilde x)-f(x)\big)]\\
    %\big(\frac{\widehat v_i(\omega_\ell)}{p_n(\omega_\ell)}e^{{\rm i} \omega_\ell\cdot x}-v_i(x)\big)^*]\\
    &=\frac{1}{K^2}\sum_{k=1}^K
    \mathbb E_\omega [
    |\frac{a(\omega_k)}{p(\omega_k)}\sigma(\omega_k\cdot \tilde x)-f(x)|^2]\\
    &=\frac{1}{K}\big(\mathbb E_\omega [
    |\frac{a(\omega_k)\sigma(\omega_k\cdot \tilde x)}{p(\omega_k)}|^2]-|f(x)|^2\big)\,,\\
%    &\le \frac{1}{K}\big(\mathbb E_\omega [
%    |\frac{a(\omega_k)c}{p(\omega_k)}|^2]-|f(x)|^2\big)\\
%      &\le \frac{1}{K}\big(\mathcal O(\|a\|^2_{L^1(\rset^d)})-|f(x)|^2\big)
%    \,,\\
    \end{split}
\end{equation*}
and by  \eqref{ap_bound} and \eqref{2moment},  we obtain 
 \begin{equation}\label{gen_error2}
 \begin{split}
\mathbb E_{\mathbf x}\big[\mathbb E_\omega[|\Bw(x)-f(x)|^2]\big]
&= \frac{1}{K}\mathbb E_{\mathbf x}\big[ \int_{\mathbb R^{d+1}}\frac{|a(\omega)\sigma(\omega\cdot\tilde x)|^2}{p(\omega)}{\rm d}\omega -|f(x)|^2\big]\\
&=\frac{1}{K} \big(\int_{\mathbb R^{d+1}}\frac{a^2(\omega)\mathbb E_{\mathbf x}[\sigma^2(\omega\cdot\tilde x)]}{p(\omega)}{\rm d}\omega 
-\mathbb E_{\mathbf x}[|f(\mathbf x)|^2]\big)\\
& =:\frac{1}{K}\big(C_p-\mathbb E_{\mathbf x}[|f(\mathbf x)|^2]\big)= \mathcal O(\frac{1}{K})\,.\\
%\mathbb E_\omega[L\big(\Bw_i(x),v_i(x)\big)] 
%&= \frac{1}{K}\big( \int_{\mathbb R^{3n}}\frac{|\widehat v_i(\omega)|^2(1+|\omega|^2)}{p_n(\omega)}{\rm d}\omega -|v_i(x)|^2 -|v'_i(x)|^2\big)\,,
\end{split}
\end{equation}

Assumption \eqref{ap_bound} 
also yields
\begin{equation}\label{C''2}
|\hat a_k | \le \frac{C''}{K}=\mathcal O(K^{-1})\,.
\end{equation}
We note that the bound of the variance of $\Bw$ in  \eqref{gen_error2} is minimized by using the  probability distribution 
%\begin{equation}\label{p_opt}
%%p(\omega)=
%\frac{|a(\omega)|}{\|a(\omega)\|_{L^1}}\,.
%\end{equation}
%and  we assume
$%\begin{equation}\label{p_opt}
\omega\mapsto \frac{|a(\omega)|(\mathbb E_{\mathbf x}[\sigma^2(\omega\cdot x)])^{1/2}}{\|a(\omega)(\mathbb E_{\mathbf x}[\sigma^2(\omega\cdot x)])^{1/2}\|_{L^1}}\,.
$%\end{equation}

Let $\beta$ be a solution to the optimization problem \eqref{min_J0}, %with positive $\alpha_j$, 
then \eqref{C''2} and \eqref{gen_error2} imply
that the expected training error satisfies 
\begin{equation}\label{tr}
\begin{split}
&\mathbb E_{\{x_j\}}\Big[\mathbb E_\omega\big[\min_{\hat\beta\in \mathbb R^K }\Big( %\lim_{J\to\infty}
\hat{\mathbb E}_J[|\beta(x)-f(x)|^2] +\lambda_1\sum_{k=1}^K |\hat\beta_k|^2+\lambda_2(\sum_{k=1}^K |\hat\beta_k|^2)^2\Big)\big]\Big]\\
&\le \mathbb E_{\{x_j\}}\Big[\mathbb E_\omega\big[\hat{\mathbb E}_J[|\bar f(x)-f(x)|^2]+\lambda_1\sum_{k=1}^K|\hat a_k|^2
+\lambda_2(\sum_{k=1}^K|\hat a_k|^2)^2\big]\Big]\\
&\le \frac{C_p}{K}(1+\lambda_1)
+ \frac{\lambda_2}{K^4}\sum_{k=1}^K\sum_{\ell=1}^K\mathbb E_\omega[
\frac{|a(\omega_k)|^2}{p^2(\omega_k)}
\frac{|a(\omega_\ell)|^2}{p^2(\omega_\ell)}]\\
&\le \frac{C_p}{K}(1+\lambda_1)
+ \frac{\lambda_2}{K^4}\big(\sum_{k\ne \ell}\mathbb E_\omega[
\frac{|a(\omega)|^2}{p^2(\omega)}]
\mathbb E_\omega[\frac{|a(\omega)|^2}{p^2(\omega)}] + \sum_{k=\ell}\mathbb E_\omega[
\frac{|a(\omega)|^4}{p^4(\omega)}]\big)\\
&\le \frac{C_p}{K}(1+\lambda_1)+\frac{\lambda_2\big(C_p^2+K^{-1}C_p(C'')^2\big)}{K^2}=:\mathcal E_t(C_p,C'')\\
&=\mathcal O(\frac{1}{K} +\frac{\lambda_1}{K} +\frac{\lambda_2}{K^2})
\,.\\
\end{split}
\end{equation}

The expected generalization error (i.e., the test error) can by \eqref{tr} be written as follows:
\begin{equation}\label{LLE}
\begin{split}
&\mathbb E_{\{x_j\}}\Big[\mathbb E_\omega\big[\mathbb E_{\mathbf x}[L\big(\beta(\mathbf x),f(\mathbf x)\big)]\big]\Big]\\
&=\mathbb E_{\{x_j\}}\Big[\mathbb E_\omega\big[\hat{\mathbb E}_J[L\big(\beta(x_j),f(x_j)\big)]\big]\Big]
\\&\quad
 +\mathbb E_{\{x_j\}}\Big[\mathbb E_\omega\big[\mathbb E_{\mathbf x}[L\big(\beta(\mathbf x),f(\mathbf x)\big)]\big]\Big]
-\mathbb E_{\{x_j\}}\Big[\mathbb E_\omega\big[\hat{\mathbb E}_J[L\big(\beta(x_j),f(x_j)\big)]\big]\Big]\\
&=\mathcal E_t  %}_{=\mathcal O(\frac{1}{K}+\frac{\lambda_1}{K} +\frac{\lambda_2}{K^2})}\\
%&\quad
+\mathbb E_{\{x_j\}}\Big[\mathbb E_\omega\big[\mathbb E_{\mathbf x}[L\big(\beta(\mathbf x),f(\mathbf x)\big)]\big]\Big]
-\mathbb E_{\{x_j\}}\Big[\mathbb E_\omega\big[\hat{\mathbb E}_J[L\big(\beta(x_j),f(x_j)\big)]\big]\Big]\,.
\end{split}
\end{equation}
It remains to estimate the right-hand side using  \eqref{LL}. 
%and
%the least squares Gram matrix $J^{-1}\sum_{j=1}^J\sigma(\omega_k\cdot \tilde x_j)\sigma(\omega_\ell\cdot \tilde x_j)$.

{\it Step 2.}
In this step, we will estimate the last two terms in \eqref{LLE}, %We have
\begin{equation}\label{ELL}
\begin{split}
&\mathbb E_{\{x_j\}}\Big[\mathbb E_\omega\big[{\mathbb E}_{\mathbf x}[L\big(\beta(\mathbf x),f(\mathbf x)\big)]-\hat{\mathbb E}_J[L\big(\beta(x_j),f(x_j)\big)]\big]\Big]\\
&=:\mathbb E_{\{x_j\}}\Big[\mathbb E_\omega\big[({\mathbb E}_{\mathbf x}-\hat{\mathbb E}_J)[L\big(\beta(x),f(x)\big)]\big]\Big]\\
&=\mathbb E_{\{x_j\}}\Big[\mathbb E_\omega\big[({\mathbb E}_{\mathbf x}-\hat{\mathbb E}_J)[
\sum_{k=1}^K\sum_{\ell=1}^K \hat\beta_k\hat\beta_\ell
\sigma(\omega_k\cdot \tilde x)\sigma(\omega_\ell\cdot \tilde x)
%\eta_k\eta_\ell^*e^{{\rm i}(\omega_k-\omega_\ell)\cdot x}
]\big]\Big]\\
&\quad -2\mathbb E_{\{x_j\}}\Big[\mathbb E_\omega\big[({\mathbb E}_{\mathbf x}-\hat{\mathbb E}_J)[\sum_{k=1}^K
\hat\beta_k\sigma(\omega_k\cdot \tilde x)f(x)]\big]\Big]
%
%(\eta_k\hat\eta_\ell^*e^{{\rm i}(\omega_k-\omega_\ell)\cdot x}+
%\eta_\ell^*\hat\eta_ke^{{\rm i}(\omega_k-\omega_\ell)\cdot x})]\big]
+\mathbb E_{\{x_j\}}\Big[\mathbb E_\omega
\big[({\mathbb E}_{\mathbf x}-\hat{\mathbb E}_J)|f(x)|^2]\big]\Big]\,.\\
\end{split}
\end{equation}

%{\it Step 3.}
Cauchy's inequality implies that, for any $\delta>0$,
 \begin{equation}\label{Ho}
\begin{split}
&|\mathbb E_\omega\big[\sum_{k=1}^K\sum_{\ell=1}^K\hat\beta_k\hat\beta_\ell \big(\hat{\mathbb E}_J[
\sigma(\omega_k\cdot \tilde x)\sigma(\omega_\ell\cdot \tilde x)]
-\mathbb E_{\mathbf x}[\sigma(\omega_k\cdot \tilde x)\sigma(\omega_\ell\cdot \tilde x)]
\big)\big]|\\
&\le\Big( \mathbb E_\omega\big[\sum_{k=1}^K\sum_{\ell=1}^K(\hat\beta_k\hat\beta_\ell )^2\big]\Big)^{1/2}
\Big( \mathbb E_\omega\big[\sum_{k=1}^K\sum_{\ell=1}^K\big((\hat{\mathbb E}_J-\mathbb E_{\mathbf x})[
\sigma(\omega_k\cdot \tilde x)\sigma(\omega_\ell\cdot \tilde x)]\big)^2\big]\Big)^{1/2}\\
&\le \frac{\lambda_2\delta }{2}\mathbb E_\omega\big[\big(\sum_{k=1}^K(\hat\beta_k)^2\big)^2\big]
+\frac{1}{2\lambda_2\delta} \mathbb E_\omega\big[\sum_{k=1}^K\sum_{\ell=1}^K\big((\hat{\mathbb E}_J-\mathbb E_{\mathbf x})[
\sigma(\omega_k\cdot \tilde x)\sigma(\omega_\ell\cdot \tilde x)]\big)^2\big]
\,,
\end{split}
\end{equation}
where the last inequality used that two real numbers $a$ and $b$, for any positive $\gamma$, satisfy $2|ab|\le \frac{a^2}{\gamma} + \gamma b^2$.
By \eqref{tr}, we obtain
\begin{equation}\label{Ew1}
\frac{\lambda_2\delta}{2}\mathbb E_{\{x_j\}}\Big[\mathbb E_\omega\big[\big(\sum_{k=1}^K\sum_{\ell=1}^K(\hat\beta_k)^2\big)^2\big]\Big]\le\frac{\delta}{2} \mathcal E_t=\mathcal O(\frac{1}{K} +\frac{\lambda_1}{K} +\frac{\lambda_2}{K^2})\,.
\end{equation}

We use the independence of all $x_j$ and $\omega_k$ together with the zero expected value \[
\mathbb E_{\{x_j\}}[ \sigma(\omega_k\cdot \tilde x_j)\sigma(\omega_\ell\cdot \tilde x_j)]-
\mathbb E_{\mathbf x}[
\sigma(\omega_k\cdot \tilde x)\sigma(\omega_\ell\cdot \tilde x)]=0\]
 to estimate the last term in \eqref{Ho}% as follows
\begin{equation}\label{ExW}
\begin{split}
&\mathbb E_{\{x_j\}}\Big[\mathbb E_\omega\big[\sum_{k=1}^K\sum_{\ell=1}^K\big((\hat{\mathbb E}_J-\mathbb E_{\mathbf x})[
\sigma(\omega_k\cdot \tilde x)\sigma(\omega_\ell\cdot \tilde x)]\big)^2\big]\Big]\\
&=\sum_{k=1}^K\sum_{\ell=1}^K\mathbb E_\omega\Big[\mathbb E_{\{x_j\}}\big[\big(\hat{\mathbb E}_J[\sigma(\omega_k\cdot \tilde x)\sigma(\omega_\ell\cdot \tilde x)]-
\underbrace{\mathbb E_{\mathbf x}[
\sigma(\omega_k\cdot \tilde x)\sigma(\omega_\ell\cdot \tilde x)]}_{=:E_{k,\ell}}\big)^2\big]\Big]\\
&=\sum_{k=1}^K\sum_{\ell=1}^K\mathbb E_\omega\big[\mathbb E_{\{x_j\}}[
\sum_{j=1}^J\sum_{i=1}^J
\frac{\big( \sigma(\omega_k\cdot \tilde x_j)\sigma(\omega_\ell\cdot \tilde x_j) - E_{k,\ell}\big)}{J}
\frac{\big( \sigma(\omega_k\cdot \tilde x_i)\sigma(\omega_\ell\cdot \tilde x_i) - E_{k,\ell}\big)}{J}
]\big]\\
&=\sum_{k=1}^K\sum_{\ell=1}^K\mathbb E_\omega\big[\mathbb E_{\{x_j\}}[
J^{-2}\sum_{j=1}^J\big( \sigma(\omega_k\cdot \tilde x_j)\sigma(\omega_\ell\cdot \tilde x_j) - E_{k,\ell}\big)^2]\big]\\
&=J^{-1}\sum_{k=1}^K\sum_{\ell=1}^K\mathbb E_\omega\big[\mathbb E_{\{x_j\}}[\big( \sigma(\omega_k\cdot \tilde x_j)\sigma(\omega_\ell\cdot \tilde x_j) - E_{k,\ell}\big)^2]\big]\\
&=\frac{1}{J}\Big(\sum_{k\ne\ell}\underbrace{\mathbb E_\omega\big[\mathbb E_{\mathbf x}[\big( \sigma(\omega_k\cdot \tilde x)\sigma(\omega_\ell\cdot \tilde x) - E_{k,\ell}\big)^2]\big]}_{=:C_2'}
+\sum_{k=\ell}\underbrace{\mathbb E_\omega\big[
\mathbb E_{\mathbf x}[\big( \sigma^2(\omega_k\cdot \tilde x) - E_{k,k}\big)^2]\big]}_{=:C_2''}\Big)\\
%&=\frac{1}{J}\sum_{k=1}^K\sum_{\ell=1}^K
%\underbrace{\mathbb E_\omega\big[\mathbb E_{\mathbf x}[\big( \sigma(\omega_k\cdot \tilde x)\sigma(\omega_\ell\cdot \tilde x) - E_{k,\ell}\big)^2]\big]}_{=:C_2}\\
&= \frac{K^2-K}{J}C_2' + \frac{K}{J}C_2''\\
&=\mathcal O(\frac{K^2}{J})\,,
\end{split}
\end{equation}
with \eqref{2moment} used in the last step,
and together with \eqref{Ew1} estimate \eqref{ExW} implies for any $\delta>0$
\begin{equation}\label{Exws}
\begin{split}
&|\mathbb E_{\{x_j\}}\Big[\mathbb E_\omega\big[\sum_{k=1}^K\sum_{\ell=1}^K\hat\beta_k\hat\beta_\ell \big(\hat{\mathbb E}_J[
\sigma(\omega_k\cdot \tilde x)\sigma(\omega_\ell\cdot \tilde x)]
-\mathbb E_{\mathbf x}[\sigma(\omega_k\cdot \tilde x)\sigma(\omega_\ell\cdot \tilde x)]
\big)\big]\Big]|\\
&\le \frac{\delta}{2}\mathcal E_t + \frac{1}{2\delta\lambda_2}( \frac{K^2-K}{J}C_2' + \frac{K}{J}C_2'')\\
&=
\mathcal O(\frac{1}{K} +\frac{\lambda_1}{K} +\frac{\lambda_2}{K^2}+ \frac{K^2}{\lambda_2 J})\,.
\end{split}
\end{equation}
We note that the crucial cancellation in \eqref{ExW} for $i\ne j$
uses that $\omega_k$ is independent of the data $\{x_j\}$ which holds
for random feature approximations but not for neural networks based on
minimization over both $\hat\beta$ and $\omega$.

The second term on the right-hand side of \eqref{ELL} can be estimated similarly by using Cauchy's inequality as
\[
\begin{split}
&|\mathbb E_{\{x_j\}} \Big[\mathbb E_\omega\big[(\hat{\mathbb E}_J-\mathbb E_{\mathbf x})[\sum_{k=1}^K\hat{\beta}_k\sigma(\omega_k\cdot x_j)f(x_j)]\big]\Big] |\\&=
|\mathbb E_{\{x_j\}} \Big[\mathbb E_\omega\big[\sum_{k=1}^K\hat{\beta}_k(\hat{\mathbb E}_J-\mathbb E_{\mathbf x})[\sigma(\omega_k\cdot x_j)f(x_j)]\big]\Big]|\\
&\le
\mathbb E_{\{x_j\}} \Big[\big(\mathbb E_{\omega}[\sum_{k=1}^K\hat{\beta}_k^2]\big)^{1/2}
\Big(\mathbb E_\omega\big[\sum_{k=1}^K\big((\hat{\mathbb E}_J-\mathbb E_{\mathbf x})
[\sigma(\omega_k\cdot x_j)f(x_j)]\big)^2\big]\Big)^{1/2}\Big]\\
&\le \frac{\lambda_1\delta}{2}\mathbb E_{\{x_j\}} \big[\mathbb E_\omega[\sum_{k=1}^K\hat{\beta}_k^2]\big]
+\frac{1}{2\lambda_1\delta}\mathbb E_{\{x_j\}} \Big[\mathbb E_\omega\big[\sum_{k=1}^K\big(
(\hat{\mathbb E}_J-\mathbb E_{\mathbf x})[\sigma(\omega_k\cdot x_j)f(x_j)]\big)^2\big]\Big]\,,\\
%&= \frac{\lambda_1}{2}\mathbb E_{\{x_j,\xi_j\}} \big[\mathbb E_\omega[\sum_{k=1}^K\hat{\beta}_k^2]\big]
%+\frac{1}{2\lambda_1}\mathbb E_\omega \big[\mathbb E_{\{x_j,\xi_j\}}[\sum_{k=1}^K\big(\hat{\mathbb E}_J[\sigma(\omega_k\cdot x_j)\xi_j]\big)^2]\big]\\
\end{split}
\]
and the independence of $x_j$ and $\omega_k$ along with the assumption that  $\|f\|_{L^\infty}\le C$, imply as in \eqref{ExW}  that 
%and the independence of $x_j$ and $\omega_k$ implies as in \eqref{ExW}  that 
\[
\begin{split}
&\frac{1}{2\lambda_1\delta}\mathbb E_{\{x_j\}} \Big[\mathbb E_\omega\big[\sum_{k=1}^K\big(
(\hat{\mathbb E}_J-\mathbb E_{\mathbf x})[\sigma(\omega_k\cdot x_j)f(x_j)]\big)^2\big]\Big]\\
&= \frac{1}{2\lambda_1\delta}\mathbb E_\omega \Big[\mathbb E_{\{x_j\}}\big[\sum_{k=1}^K\big(
(\hat{\mathbb E}_J-\mathbb E_{\mathbf x})[\sigma(\omega_k\cdot x_j)f(x_j)]\big)^2\big]\Big]\\
&=\frac{1}{2\lambda_1\delta}\mathbb E_\omega \Big[\sum_{k=1}^K\mathbb E_{\{x_j\}}\big[
J^{-2}\sum_{j=1}^J\sum_{i=1}^J\big(
\sigma(\omega_k\cdot x_j)f(x_j)-\mathbb E_{\mathbf x}[\sigma(\omega_k\cdot \mathbf x)f(\mathbf x)]\big)\times\\
&\qquad\times\big(
\sigma(\omega_k\cdot x_i)f(x_i)-\mathbb E_{\mathbf x}[\sigma(\omega_k\cdot \mathbf x)f(\mathbf x)]\big)\big]\Big]\\
&=\frac{1}{2\lambda_1\delta}\mathbb E_\omega \big[\sum_{k=1}^K\mathbb E_{\{x_j\}}[J^{-2}\sum_{j=1}^J
\big(\sigma(\omega_k\cdot x_j)f(x_j)-\mathbb E_{\mathbf x}[\sigma(\omega_k\cdot \mathbf x)f(\mathbf x)]\big)^2]\big]\\
&=\frac{K}{2\lambda_1 \delta J}\underbrace{
\mathbb E_\omega 
\big[\mathbb E_{\{x_j\}}[\big(
\sigma(\omega_k\cdot x_j)f(x_j)-\mathbb E_{\mathbf x}[\sigma(\omega_k\cdot \mathbf x)f(\mathbf x)]\big)^2]\big]}_{=:C_1}\\
&=\mathcal O(\frac{K}{\lambda_1 J})\,.
\end{split}
\]

We also have the zero expected value
\begin{equation}\label{H2}
\begin{split}
\mathbb E_{\{x_j\}}\big[\mathbb E_{\mathbf x}[|f(x)|^2]-\hat{\mathbb E}_J[|f(x)|^2]\big]=0\,.
\end{split}
\end{equation}

{\it Step 3.}
By combining \eqref{LLE}, \eqref{ELL}, %\eqref{HL},  
\eqref{Ho},  %\eqref{Ew1}, 
\eqref{Exws}, \eqref{H2},  and \eqref{tr} we obtain 
%for $\delta:=2\delta_1+\delta_2$ 
\begin{equation}\label{Exw22}
\begin{split}
\mathbb E_{\{x_j\}}\Big[\mathbb E_\omega\big[\mathbb E_{\mathbf x}[L\big(\beta(x),f(x)\big)]\big] \Big]
&\le \mathcal E_t(1+\delta) +  \frac{1}{2\delta\lambda_2}( \frac{K^2-K}{J}C_2' + \frac{K}{J}C_2'')
+\frac{KC_1}{2\lambda_1 \delta J}\\
%&=: \mathcal E\\
&=
%\frac{1}{2}\mathbb E_{\{x_j\}}\Big[\mathbb E_\omega\big[\mathbb E_{\mathbf x}[L\big(\beta(x),f(x)\big)]\big]\Big]\\
% &\quad+
\mathcal O(\frac{1}{K} +\frac{\lambda_1}{K}+
\frac{K}{\lambda_1 J} +\frac{\lambda_2}{K^2}+ \frac{K^2}{\lambda_2 J}  )\,,
\end{split}
\end{equation}
which proves \eqref{gen_O3}.

\end{proof}

\subsection{The generalization error for random features with noisy data}\label{RFnoise}
The generalization error with data polluted by noise can be analyzed as follows. Assume that for the data points $x_j\in \rset^d, \ j=1,\ldots,J\,,$ we observe
\begin{equation}\label{noise1}
y_j=f(x_j)+\xi_j\,,\quad j=1,\ldots,J\,,
\end{equation}
where $\xi_j\in \rset$ are independent random samples, independent also of $x_j$ and $\omega_j$, with mean zero, $\mathbb E_\xi[\xi_j]=0$, and bounded
variance, $\mathbb E_\xi[\xi_j^2]=s^2$. Let $\mathbb E_ {\{x_j,\xi_j\}}$ denote the expected value with respect to the training points $x_j$ and noise $\xi_j$, for $j=1,\ldots, J$.

\begin{theorem}[Generalization error with noise]\label{thm_noise}
Suppose 
\[
C_3:=\mathbb E_\omega \big[\mathbb E_{\{x_j,\xi_j\}}[
\sigma^2(\omega_k\cdot x_j)]\big]=\mathcal O(1)
\]
and the assumptions in Theorem \ref{thm:generalization} hold with the difference that the optimization problem is
\begin{equation}\label{min_JN}
\begin{split}
&\min_{\hat\beta\in \rset^K } %\lim_{J\to\infty}
\Big(J^{-1}\sum_{j=1}^J|f(x_j)+\xi_j-\beta(x_j)|^2  +\lambda_1\sum_{k=1}^K|\hat\beta_k|^2
 +\lambda_2(\sum_{k=1}^K|\hat\beta_k|^2)^2\Big)
\,,\\
\end{split}
\end{equation}
then for any $\delta>0$, 
\begin{equation}\label{gen_ON}
\begin{split}
\mathbb E_ {\{x_j,\xi_j\}}\big[\mathbb E_\omega [\int_{\rset^{d}} |f(x)-\beta(x)|^2\mu(x){\rm d}x ] \big]
&\le (1+\frac{3\delta}{2}) \frac{C_p}{K}(1+\lambda_1)+\frac{\lambda_2(C_p^2+K^{-1}C_p(C'')^2)}{K^2} \\
&\quad +  \frac{1}{2\delta\lambda_2}( \frac{K^2-K}{J}C_2' + \frac{K}{J}C_2'')
+\frac{KC_1}{2\lambda_1 \delta J}
+ \frac{KC_3 s^2}{2\delta\lambda_1 J}\\
%&=\mathcal O\big(\frac{1}{K}+\lambda +(1+\frac{1}{\lambda K})(\sqrt{\frac{\log\frac{K}{\delta}}{J}})+\delta\big)\\
&=\mathcal O\big(\frac{1}{K}+\frac{\lambda_1}{K} +\frac{K(1+s^2)}{\lambda_1 J}
+\frac{\lambda_2}{K^2}+\frac{K^2}{\lambda_2 J}\big)\,,
\end{split}
\end{equation}
which for $\lambda_1=K\frac{\sqrt{1+s^2}}{\sqrt J}$ and $\lambda_2=\frac{K^2}{\sqrt J}$ becomes
\begin{equation*}\label{gen_Jxi22}
\begin{split}
\mathbb E_ {\{x_j,\xi_j\}}\big[\mathbb E_\omega [\int_{\rset^{d}} |f(x)-\beta(x)|^2\mu(x){\rm d}x ] \big]
%&=\mathcal O\big(\frac{1}{K}+\lambda +(1+\frac{1}{\lambda K})(\sqrt{\frac{\log\frac{K}{\delta}}{J}})+\delta\big)\\
&=\mathcal O\big(\frac{1}{K}+\frac{1+|s|}{\sqrt J} \big)\,,
\end{split}
\end{equation*}
and in the case $\sigma(z)=e^{{\rm i}z}$ there holds
\[
\begin{split}
C_1 &\le \mathbb E_x[ |f(x)|^2] \,,\\
C_2' &\le 1\,,\\
C_2'' &=0\,,\\
C_3 &=1\,,\\
\end{split}
\]
which implies
\[
\begin{split}
&\mathbb E_ {\{x_j,\xi_j\}}\big[\mathbb E_\omega [\int_{\rset^{d}} |f(x)-\beta(x)|^2\mu(x){\rm d}x ] \big]\\
&\le \frac{(1+3\delta/2)C_p}{K} + (\frac{(1+3\delta/2) C_1C_p}{J\delta})^{1/2}+(\frac{(1+3\delta/2)(1+K^{-1})(C_p^2+\frac{(C'')^2C_p}{K})}{J\delta})^{1/2}
\end{split}
\]
for $\lambda_1=K(\frac{(1+3\delta/2)C_pc_1}{J\delta})^{1/2}$ and $\lambda_2=K(\frac{(1+3\delta/2)(1+K^{-1})(C_p^2+(C'')^2C_pK^{-1})}{J\delta})^{-1/2}$. 
\end{theorem}

We note that the effect of the noise in fact vanishes as the number of data points $J$ tends to infinity,
which requires the no bias condition $\mathbb E_\xi[\xi_j]=0$. More precisely, for the optimal penalty parameter choice  $\lambda_1=K/\sqrt J$ and $\lambda_2=K^2/\sqrt J$ and the amount of data $J$ proportional to $K^\alpha$, the generalization error bound is  
$\mathcal O(K^{-1})$ for $\alpha\ge 2$
and $\mathcal O(K^{-\alpha/2})$ for $\alpha<2$. 
%In the case without noise, i.e. $s=0$, the generalization error is $\mathcal O(K^{-1})$ for $\alpha\ge 2$  and $\mathcal O(K^{-\alpha/2})$ for $\alpha<2$. 
This indicates that neural network approximation with high accuracy in high dimension, $d\gg  1$, is possible but requires many data points. %in particular for noisy data.

\begin{proof}
The theorem is an extension of Theorem \ref{thm:generalization}, and the proof has two additional steps:

\begin{itemize}
\item[Step 1] estimates the training error, including noisy data, and
\item[Step 2] applies Step 1 combined with  \eqref{ELL},  \eqref{Ho}, and  \eqref{H2} to prove \eqref{gen_ON}.
\end{itemize}

{\it Step 1.} Assume that $\beta$ is a solution to the optimization \eqref{min_JN}. 
We have by the optimization \eqref{min_JN}, \eqref{mean_ok}, \eqref{gen_error2}, and \eqref{C''2} together with the independence of $x_j$ and $\xi_j$
\[
\begin{split}
&\mathbb E_{\{x_j,\xi_j\},\omega}\big[\hat{\mathbb E}_J[|f(x_j)-\beta(x_j)+\xi_j|^2] +\lambda_1\sum_{k=1}^K|\hat\beta_k|^2
+\lambda_2(\sum_{k=1}^K|\hat\beta_k|^2)^2-\hat{\mathbb E}_J[\xi_j^2]\big]\\
&\le
\mathbb E_{\{x_j,\xi_j\},\omega}\big[\hat{\mathbb E}_J[|f(x_j)-\bar f(x_j)+\xi_j|^2]
+\lambda_1\sum_{k=1}^K|\hat a_k|^2
+\lambda_2(\sum_{k=1}^K|\hat a_k|^2)^2 
-\hat{\mathbb E}_J[\xi_j^2]\big]\\
&=\mathbb E_{\{x_j,\xi_j\},\omega}\big[\hat{\mathbb E}_J[|f(x_j)-\bar f(x_j)|^2] +\lambda_1\sum_{k=1}^K|\hat a_k|^2
+\lambda_2(\sum_{k=1}^K|\hat a_k|^2)^2 \big]\\
&\quad-2\mathbb E_{\{x_j\},\omega}\Big[\hat{\mathbb E}_J\big[\underbrace{{\mathbb E}_{\{\xi_j\}}[\big(f(x_j)
-\bar f(x_j)\xi_j]}_{=0}\big]\Big] +\mathbb E_{\{x_j,\xi_j\},\omega}\big[\hat{\mathbb E}_J[\xi_j^2]-\hat{\mathbb E}_J[\xi_j^2]\big]\\
&=\mathbb E_{\{x_j,\xi_j\},\omega}\big[\hat{\mathbb E}_J[|f(x_j)-\bar f(x_j)|^2\big] +\lambda_1\sum_{k=1}^K|\hat a_k|^2
+\lambda_2(\sum_{k=1}^K|\hat a_k|^2)^2 \big] \\
&\le\mathcal E_t\\
&=\mathcal O(\frac{1}{K} + \frac{\lambda_1}{K}+ \frac{\lambda_2}{K^2})
\end{split}
\]
and the rewriting
\[
\begin{split}
&\mathbb E_\omega\big[\hat{\mathbb E}_J[|f(x_j)-\beta(x_j)+\xi_j|^2] +\lambda_1\sum_{k=1}^K|\hat\beta_k|^2
+\lambda_2(\sum_{k=1}^K|\hat\beta_k|^2)^2\big]-\hat{\mathbb E}_J[\xi_j^2]\\
%&=\mathbb E_\omega\big[\hat{\mathbb E}_J[|f(x_j)-\beta(x_j)|^2] \big]
%-2\mathbb E_\omega \big[\hat{\mathbb E}_J[\big(f(x_j)-\beta(x_j)\big)\xi_j]+\lambda(\sum_{k=1}^K|\beta_k|)^2\big]\\
&=\mathbb E_\omega\big[\hat{\mathbb E}_J[|f(x_j)-\beta(x_j)|^2] +\lambda_1\sum_{k=1}^K|\hat\beta_k|^2
+\lambda_2(\sum_{k=1}^K|\hat\beta_k|^2)^2\big]\\
&\quad -2\mathbb E_\omega\big[\hat{\mathbb E}_J[\big(f(x_j)-\beta(x_j)\big)\xi_j]\big]\,,\\
\end{split}
\]
implies the training error estimate
\[
\begin{split}
&\mathbb E_{\{x_j,\xi_j\}}\Big[\mathbb E_\omega\big[\hat{\mathbb E}_J[|f(x_j)-\beta(x_j)|^2] +\lambda_1\sum_{k=1}^K|\hat\beta_k|^2
+\lambda_2(\sum_{k=1}^K|\hat\beta_k|^2)^2\big]\Big]\\
&\le \mathcal E_t %\mathcal O(\frac{1}{K} + \frac{\lambda_1}{K}+ \frac{\lambda_2}{K^2})%\\&\quad 
+2\mathbb E_{\{x_j,\xi_j\}}\Big[\mathbb E_\omega\big[\hat{\mathbb E}_J[\big(f(x_j)-\beta(x_j)\big)\xi_j]\big]\Big]\,.
\end{split}
\]

We seek a bound on the expected training error
\[
\mathbb E_{\{x_j,\xi_j\}}\Big[\mathbb E_\omega\big[\hat{\mathbb E}_J[|f(x_j)-\beta(x_j)|^2] +\lambda_1\sum_{k=1}^K|\hat\beta_k|^2
+\lambda_2(\sum_{k=1}^K|\hat\beta_k|^2)^2\big]\Big]
\]
and it remains to estimate  the following: % We start with the last one:
\[
\begin{split}
\mathbb E_{\{x_j,\xi_j\}}\Big[\mathbb E_\omega\big[\hat{\mathbb E}_J[\big(f(x_j)-\beta(x_j)\big)\xi_j]\big]\Big]
&=\mathbb E_\omega\Big[\mathbb E_{\{x_j,\xi_j\}}\big[\hat{\mathbb E}_J[f(x_j)\xi_j]\big]\Big]\\
&\quad
-\mathbb E_{\{x_j,\xi_j\}} \Big[\mathbb E_\omega\big[\hat{\mathbb E}_J[\sum_{k=1}^K\hat{\beta}_k\sigma(\omega_k\cdot x_j)\xi_j]\big]\Big]\,.
%=:\hat{\mathbb E}_J\big[\mathbb E_\xi[\sum_{k=1}^K\hat{\beta}_k\eta_j \big]\,.
\end{split}
\]
The independence of $x_j$ and $\xi_j$ implies that the first term on the right-hand side vanishes
\[
\mathbb E_{\{x_j,\xi_j\}}\big[\hat{\mathbb E}_J[f(x_j)\xi_j]\big]
=\hat{\mathbb E}_J\big[\mathbb E_{\mathbf x}\big[f(x_j)\underbrace{\mathbb E_{\xi_j}[\xi_j}_{=0}]\big]=0\,.
\]
%and we start with the term based on $\eta_j$.
The second term on the right-hand side can be estimated, as in \eqref{Ho}:
\[
\begin{split}
&|\mathbb E_{\{x_j,\xi_j\}} \Big[\mathbb E_\omega\big[\hat{\mathbb E}_J[\sum_{k=1}^K\hat{\beta}_k\sigma(\omega_k\cdot x_j)\xi_j]\big]\Big] |\\&=
|\mathbb E_{\{x_j,\xi_j\}} \Big[\mathbb E_\omega\big[\sum_{k=1}^K\hat{\beta}_k\hat{\mathbb E}_J[\sigma(\omega_k\cdot x_j)\xi_j]\big]\Big]|\\
&\le
\mathbb E_{\{x_j,\xi_j\}} \Big[\big(\mathbb E_\omega[\sum_{k=1}^K\hat{\beta}_k^2]\big)^{1/2}
\Big(\mathbb E_\omega\big[\sum_{k=1}^K\big(\hat{\mathbb E}_J[\sigma(\omega_k\cdot x_j)\xi_j]\big)^2\big]\Big)^{1/2}\Big]\\
&\le \frac{\lambda_1\delta}{2}\mathbb E_{\{x_j,\xi_j\}} \big[\mathbb E_\omega[\sum_{k=1}^K\hat{\beta}_k^2]\big]
+\frac{1}{2\lambda_1\delta}\mathbb E_{\{x_j,\xi_j\}} \Big[\mathbb E_\omega\big[\sum_{k=1}^K\big(\hat{\mathbb E}_J[\sigma(\omega_k\cdot x_j)\xi_j]\big)^2\big]\Big]\,.\\
%&= \frac{\lambda_1}{2}\mathbb E_{\{x_j,\xi_j\}} \big[\mathbb E_\omega[\sum_{k=1}^K\hat{\beta}_k^2]\big]
%+\frac{1}{2\lambda_1}\mathbb E_\omega \big[\mathbb E_{\{x_j,\xi_j\}}[\sum_{k=1}^K\big(\hat{\mathbb E}_J[\sigma(\omega_k\cdot x_j)\xi_j]\big)^2]\big]\\
\end{split}
\]
The independence of $\xi_j$ and $x_j$ together with $\mathbb E_{\xi_j}[\xi_j]=0$ imply that the last term above has the following bound:
\[
\begin{split}
&\frac{1}{2\lambda_1\delta}\mathbb E_{\{x_j,\xi_j\}} \Big[\mathbb E_\omega\big[\sum_{k=1}^K\big(\hat{\mathbb E}_J[\sigma(\omega_k\cdot x_j)\xi_j]\big)^2\big]\Big]\\
&= \frac{1}{2\lambda_1\delta}\mathbb E_\omega \Big[\mathbb E_{\{x_j,\xi_j\}}\big[\sum_{k=1}^K\big(\hat{\mathbb E}_J[\sigma(\omega_k\cdot x_j)\xi_j]\big)^2\big]\Big]\\
&=\frac{1}{2\lambda_1\delta}\mathbb E_\omega \big[\sum_{k=1}^K\mathbb E_{\{x_j,\xi_j\}}[J^{-2}\sum_{j=1}^J\sum_{i=1}^J
\sigma(\omega_k\cdot x_j)\xi_j  \sigma(\omega_k\cdot x_i)\xi_i]\big]\\
&=\frac{1}{2\lambda_1\delta}\mathbb E_\omega \big[\sum_{k=1}^K\mathbb E_{\{x_j,\xi_j\}}[J^{-2}\sum_{j=1}^J
\sigma^2(\omega_k\cdot x_j)\xi_j^2]\big]\\
&=\frac{K}{2\lambda_1\delta J}\mathbb E_\omega \big[\mathbb E_{\{x_j,\xi_j\}}[
\sigma^2(\omega_k\cdot x_j)\xi_j^2]\big]\\
&=\frac{K}{2\lambda_1\delta J}\underbrace{\mathbb E_\omega \big[\mathbb E_{\{x_j,\xi_j\}}[
\sigma^2(\omega_k\cdot x_j)]\big]}_{=:C_3}s^2\\
&=\mathcal O(\frac{Ks^2}{\lambda_1 J})\,.
\end{split}
\]

In conclusion, the expected training error bound is
\begin{equation}\label{train}
\begin{split}
&\mathbb E_{\{x_j,\xi_j\}}\Big[\mathbb E_\omega\big[\hat{\mathbb E}_J[|f(x_j)-\beta(x_j)|^2] +\lambda_1\sum_{k=1}^K|\hat\beta_k|^2
+\lambda_2(\sum_{k=1}^K|\hat\beta_k|^2)^2\big]\Big]\\
&\le (1+\frac{3\delta}{2})\mathcal E_t + \frac{KC_3 s^2}{2\delta\lambda_1 J}\\
&=\mathcal O(\frac{1}{K} + \frac{\lambda_1}{K} + \frac{Ks^2}{\lambda_1 J} + \frac{\lambda_2}{K^2})\,.
\end{split}
\end{equation}

{\it Step 2.} The expected generalization error, corresponding to \eqref{LLE},
but now with data including noise, becomes
\begin{equation}\label{LLEJ}
\begin{split}
&\mathbb E_{\{x_j,\xi_j\}}\Big[\mathbb E_\omega\big[\mathbb E_{\mathbf x}[L\big(\beta(\mathbf x),f(\mathbf x)\big)]\big]\Big]\\
&=\mathbb E_{\{x_j,\xi_j\}}\Big[\mathbb E_\omega\big[\hat{\mathbb E}_J[L\big(\beta(x_j),f(x_j)\big)]\big]\Big]
\\&\quad
 +\mathbb E_{\{x_j,\xi_j\}}\Big[\mathbb E_\omega\big[\mathbb E_{\mathbf x}[L\big(\beta(\mathbf x),f(\mathbf x)\big)]\big]
-\mathbb E_\omega\big[\hat{\mathbb E}_J[L\big(\beta(x_j),f(x_j)\big)]\big]\Big]\\
&= (1+\frac{3\delta}{2})\mathcal E_t + \frac{KC_3 s^2}{2\delta\lambda_1 J}\\
%\mathcal O(\frac{1}{K} + \frac{\lambda_1}{K} + \frac{Ks^2}{\lambda_1 J} + \frac{\lambda_2}{K^2})
&\quad
+\mathbb E_{\{x_j,\xi_j\}}\Big[\mathbb E_\omega\big[\mathbb E_{\mathbf x}[L\big(\beta(\mathbf x),f(\mathbf x)\big)]\big]
-\mathbb E_\omega\big[\hat{\mathbb E}_J[L\big(\beta(x_j),f(x_j)\big)]\big]\Big]\,.
\end{split}
\end{equation}
By combining \eqref{LLEJ}, \eqref{ELL}, %\eqref{HL},  
\eqref{Ho},  %\eqref{Ew1}, 
\eqref{Exws} and \eqref{H2} we obtain as in \eqref{Exw22}
%\eqref{LLEJ}, \eqref{ELL}, \eqref{HL},  \eqref{Ho},  \eqref{H2},  and \eqref{trJ} 
%we obtain 
\[
\begin{split}
&\mathbb E_{\{x_j,\xi_j\}}\Big[\mathbb E_\omega\big[\mathbb E_{\mathbf x}[L\big(\beta(\mathbf x),f(\mathbf x)\big)]\big]\Big]\\
&\le  (1+\frac{3\delta}{2})\mathcal E_t 
+  \frac{1}{2\delta\lambda_2}( \frac{K^2-K}{J}C_2' + \frac{K}{J}C_2'')
+\frac{KC_1}{2\lambda_1 \delta J}
+ \frac{KC_3 s^2}{2\delta\lambda_1 J}\\
&=
\mathcal O(\frac{1}{K} + \frac{\lambda_1}{K} + \frac{K(1+s^2)}{\lambda_1 J} + \frac{\lambda_2}{K^2} +\frac{K^2}{\lambda_2 J})
\end{split}
\]
which proves the Theorem.
\end{proof}

\section*{Acknowledgments}
This research was supported by the 
Swedish Research Council grant 2019-03725.
The computations were enabled by resources provided by the National Academic Infrastructure for Supercomputing in Sweden (NAISS) at the PDC Center for High Performance Computing, KTH Royal Institute of Technology, partially funded by the Swedish Research Council through grant agreement no. 2022-06725. The work of Xin Huang was supported in part by the Kempe Stiftelserna project JCSMK23-0168.
We also acknowledge the financial support from the King Abdullah University of Science and Technology (KAUST) Office 
of Sponsored Research (OSR) under Award No. OSR-2019-CRG8-4033, and from the Alexander von Humboldt Foundation.

% %%%%% Anamika editing
\section*{Funding}
Funding information here.

%\bibliographystyle{plain}
%\bibliography{reference}

%USE THE BELOW OPTIONS IN CASE YOU NEED AUTHOR YEAR FORMAT.
%\bibliographystyle{abbrvnat}
%\bibliography{reference}

\end{document}